\documentclass[12pt,reqno]{amsart} 

\usepackage{mathrsfs,amsfonts,amsthm,amsmath,amssymb,thmtools} 
\usepackage[dvipsnames]{xcolor} 

\usepackage{enumitem}                
\usepackage{lineno}                 

\usepackage[numbers,sort&compress]{natbib}  

\usepackage{graphics}
\usepackage{graphicx}   
\graphicspath{{images/}} 
\usepackage{subfig}     
\usepackage{rotating}

\usepackage[hypertexnames=false, hyperindex, pdfencoding=auto, psdextra, bookmarksdepth=4]{hyperref}
\hypersetup{
    bookmarksopen      = true,
    bookmarksopenlevel = 1,
    bookmarksnumbered  = true,
    pdftitle    = {},
    pdfauthor   = {},
    linktoc     = page,
    anchorcolor = green,
    colorlinks  = true,        
    linkcolor   = black,
    citecolor   = blue,
    filecolor   = blue,
    urlcolor    = magenta,
}

\usepackage[capitalise, noabbrev, nameinlink]{cleveref}    

\makeatletter
\DeclareRobustCommand{\Crefinorder}[1]{%
  {\@cref@sortfalse\@cref@compressfalse\Cref{#1}}}
\makeatother

\newtheorem{innercustomthm}{Theorem}
\newenvironment{customthm}[1]
    {\renewcommand\theinnercustomthm{#1}\innercustomthm}
    {\endinnercustomthm}

\theoremstyle{plain}

\newtheorem{theorem}{Theorem}[section]
\newtheorem{lemma}[theorem]{Lemma}

\newtheorem{proposition}[theorem]{Proposition}
\newtheorem{corollary}[theorem]{Corollary}

\theoremstyle{definition}   
\newtheorem{definition}[theorem]{Definition}
\newtheorem{remark}[theorem]{Remark}
\newtheorem{example}[theorem]{Example}

\makeatletter
    \def\subsection{\@startsection{subsection}{2}%
    \z@{.5\linespacing\@plus.7\linespacing}{.3\linespacing}%
    {\normalfont\bfseries}}
\makeatother 

\makeatletter                                         
    \newcommand{\LeftEqNo}{\let\veqno\@@leqno}        
\makeatother

\numberwithin{equation}{section}                      

\usepackage{etoolbox}

\renewcommand{\contentsnamefont}{\normalfont\bfseries}

\makeatletter

\patchcmd{\@starttoc}
  {\centering\contentsnamefont}
  {\raggedright\contentsnamefont}
  {}{}

\patchcmd{\@tocline}
  {\hfil}
  {\dotfill}
  {}{}

\renewcommand{\l@section}
  {\@tocline{1}{5pt}{0pt}{2em}{}}
\renewcommand{\l@subsection}
  {\@tocline{2}{0pt}{1.5em}{2.8em}{}}

\makeatother

\begin{document}
\raggedbottom
\

\vspace{-2cm}

\title[Homogeneous operators in the hyperfinite type $\mathrm{II}_1$ factor]
{A formal representation of shift structures in the hyperfinite type $\mathrm{II}_1$ factor}

\author{Junsheng Fang}
\address{Junsheng Fang, School of Mathematical Sciences, Hebei Normal University, Shijiazhuang, 050024, China}
\email{jfang@hebtu.edu.cn}

\author{Chunlan Jiang}
\address{Chunlan Jiang, School of Mathematical Sciences, Hebei Normal University, Shijiazhuang, 050024, China}
\email{cljiang@hebtu.edu.cn}

\author{Minghui Ma}
\address{Minghui Ma, School of Mathematical Sciences, Dalian University of Technology, Dalian, 116024, China}
\email{minghuima@dlut.edu.cn}

\author{Rui Shi}
\address{Rui Shi, School of Mathematical Sciences, Dalian University of Technology, Dalian, 116024, China}
\email{ruishi@dlut.edu.cn, ruishi.math@gmail.com}

\keywords{Weighted shift model, $K_0$-group, strongly irreducible operator, strong generator, Jones index}

\makeatletter
\@namedef{subjclassname@2020}{\textup{2020} Mathematics Subject Classification}
\makeatother
\subjclass[2020]{Primary 46L10; Secondary 47C15, 46L80}

\begin{abstract}
In this paper, we establish a connection between classical weighted shift operators and homogeneous operators associated with irrational rotation algebras in the hyperfinite type $\mathrm{II}_1$ factor $\mathcal R$, focusing on the structure of their relative commutants.
For injective homogeneous operators of degree one, we use bilateral weighted shift operators to model their relative commutants and obtain unilateral models in the noninvertible case.
We compute $K_0$-groups and provide splitting criteria for direct sums.
In particular, we construct strongly irreducible operators in $\mathcal R$ whose relative commutants have $K_0$-groups isomorphic to $\mathbb Z$ and $\mathbb Z^2$.
We realize every Jones index of an irreducible subfactor of $\mathcal R$ as $[\mathcal R:W^*(T)]$ for some strongly irreducible operator $T$ relative to $\mathcal R$.
Furthermore, we construct strong generators of $\mathcal R$ and strongly irreducible operators in every type $\mathrm{II}_1$ factor with separable predual.
\end{abstract}

\maketitle

\section{Introduction}

While concepts and techniques from operator algebras are used to solve problems in operator theory \cite{BDF-73,Voiculescu-76}, one of many approaches to understanding the structure of von Neumann algebras is to extend ideas and methods from classical matrix theory and Hilbert space operator theory to factors.
Representative examples include the dimension theory of Murray and von Neumann \cite{MvN-36}, the determinant theory of Fuglede and Kadison \cite{FK-52}, Voiculescu's random matrix models and free probability theory \cite{Voiculescu-91}, the work of Haagerup and Schultz on Brown measures and invariant subspaces \cite{HS-07,HS-09}, and the Schur-type decomposition of Dykema, Sukochev, and Zanin \cite{DSZ-15}.
In this paper, we pursue this theme by seeking counterparts of the classical theory of weighted shift operators in the hyperfinite type \(\mathrm{II}_1\) factor.

Motivated by the theory of the unweighted unilateral and bilateral shift operators, Shields \cite{Shields-74} developed a systematic account of weighted shift operators and their interplay with analytic function theory. 
A central idea in his work is to realize an injective weighted shift as multiplication by the coordinate function on a suitably weighted Hilbert space of formal power series or Laurent series. 
This approach identifies the commutant of a weighted shift with a corresponding multiplier algebra and provides a powerful tool for investigating its structure.
In particular, it yields what is now called the \emph{strong irreducibility} of injective unilateral weighted shift operators \cite[Corollary 2, p.~63]{Shields-74}.
Inspired by this point of view, we develop weighted shift models for \emph{homogeneous operators} arising from irrational rotation algebras, with the aim of bringing the classical weighted shift theory to bear on the structure of relative commutants in the hyperfinite type $\mathrm{II}_1$ factor.
Along this thread, we investigate the strong irreducibility of homogeneous operators, study the $K$-theory of their relative commutants, and develop new operator-theoretic approaches to realizing irreducible Jones index values.

Let $\mathcal A_\theta$ denote the irrational rotation $C^*$-algebra \cite{EE-93} generated by two unitary operators $u$ and $v$ satisfying
\begin{equation*}
	vu=\alpha uv,\qquad \alpha=e^{2\pi i\theta},\qquad \theta\in\mathbb R\setminus\mathbb Q.
\end{equation*}
Let $\tau$ be the unique tracial state on $\mathcal A_\theta$ and $(\pi_\tau,\mathcal H_\tau,\xi_\tau)$ the corresponding GNS triple.
Then $\pi_\tau(\mathcal A_\theta)''$ is isomorphic to the hyperfinite type $\mathrm{II}_1$ factor $\mathcal R$.
We identify $\mathcal A_\theta$ with its image under $\pi_\tau$ in $\mathcal R$ and continue to denote its generators by $u$ and $v$.
With respect to these generators, a nonzero measurable operator $T$ affiliated with $\mathcal R$ is called \emph{homogeneous} of \emph{degree} $m\in\mathbb Z$ if $T=u^mf(v)$ for some measurable function $f$ on the unit circle $\mathbb T$ equipped with normalized Haar measure $\mu$, as in \Cref{def homogeneous}.
We call $f$ the \emph{weight function} of $T$, in analogy with weighted shift operators.
The role of the scalar weight sequence in the classical theory is played here by a measurable weight function sampled along an irrational rotation.
For an injective homogeneous operator $T=uf(v)$ of degree one in $\mathcal R$, we compare the $L^2$-Fourier expansions of operators in $\{T\}'\cap\mathcal R$ with formal Laurent series to obtain a bilateral weighted shift model.
When $T$ is noninvertible, restricting the weight sequence to nonnegative indices gives a unilateral weighted shift model.

\begin{customthm}{A}\label{thm A}
Let $T=uf(v)\in\mathcal R$, where $f\ne 0$ almost everywhere.
Then there exist an injective bilateral weighted shift operator $B$ and an isometric algebra isomorphism from $\{T\}'\cap\mathcal R$ onto $\{B\}'$ that sends $T$ to $B$.
If $T$ is noninvertible, the same conclusion holds with an injective unilateral weighted shift operator in place of $B$.

Conversely, suppose that $B$ is a noninvertible injective bilateral weighted shift operator (or an injective unilateral weighted shift operator).
Then there exist an injective noninvertible operator $T=uf(v)\in\mathcal R$ and a Banach algebra isomorphism from $\{T\}'\cap\mathcal R$ onto $\{B\}'$ that sends $T$ to $B$.
\end{customthm}

The bilateral and unilateral weighted shift models follow from \Cref{thm bilateral-model,thm unilateral-model}, respectively.
The converse follows from \Cref{prop noninvertible-bilateral-realization} in the bilateral case and from \Cref{prop unilateral-realization} in the unilateral case.
As an application of the unilateral weighted shift model, we show in \Cref{cor commutant-Hinfty} that if $T=uf(v)$ is injective and noninvertible and satisfies $r(T)=\|T\|$, then $\{T\}'\cap\mathcal R$ is isometrically isomorphic to $H^\infty(\mathbb D)$, the algebra of bounded holomorphic functions on the unit disk.

To establish the weighted shift models in \Cref{thm A}, we develop Fourier methods in \Cref{sec pre} and use them in \Cref{sec commutant} to determine the relative commutants in $\mathcal R$ of injective homogeneous operators of nonzero degree.
These methods also allow us to study \emph{irreducibility} in \Cref{sec IR} and \emph{strong irreducibility} in \Cref{sec SI}.
The study of irreducible operators on separable Hilbert spaces dates back to Halmos \cite{Halmos-68}.
Since then, a great deal of fruitful work has been devoted to the structure of operators and $C^*$-algebras \cite{BDF-73,Voiculescu-76}.
Irreducibility is preserved by unitary equivalence but need not be preserved by similarity.
This motivates the notion of strong irreducibility, introduced by Gilfeather \cite{Gilfeather-72}, which requires every operator in the similarity orbit to be irreducible.
In finite dimensions, this concept characterizes matrices similar to a single Jordan block.
Jiang \cite{JC-23} viewed strongly irreducible operators as infinite-dimensional analogues of Jordan blocks.
The theory of strongly irreducible operators is developed systematically in \cite{JW-book-98,JW-book-06}.
The notion of strongly irreducible operators extends naturally to factors.
An operator $T$ in a factor $\mathcal M$ is called \emph{strongly irreducible} relative to $\mathcal M$ if the relative commutant $\{T\}'\cap\mathcal M$ contains no idempotents other than $0$ and $1$; equivalently, every operator in the similarity orbit of $T$ in $\mathcal M$ generates an irreducible subfactor of $\mathcal M$.
The authors in \cite[Theorem 8.2]{FJLX-13} proved that $u+v$ is strongly irreducible relative to $\mathcal R$.
This investigation was subsequently extended in \cite[Theorem~4.4]{ZFS-17}, which treats the strong irreducibility of homogeneous operators of degree one with continuous weights.
In \Cref{thm SI-noninvertible}, we show that every injective noninvertible homogeneous operator of nonzero degree in $\mathcal R$ is strongly irreducible relative to $\mathcal R$.
As an application, we show in \Cref{cor SI-sum-power} that $u^m+v^n$ is strongly irreducible relative to $\mathcal R$ for all nonzero integers $m$ and $n$.

The $K_0$-group is a fundamental invariant in the structure and classification theory of $C^*$-algebras \cite{Blackadar-book-98}, as illustrated by Elliott's classification of separable unital AF algebras \cite{Elliott-76}.
In operator theory, $K_0$-groups of commutants have been used to study finite strongly irreducible decompositions \cite{CFJ-02} and the similarity classification of Cowen--Douglas operators \cite{Jiang-04,JGJ-05}.
For every unital complex algebra $\mathcal A$, let $V(\mathcal A)$ denote the monoid of algebraic equivalence classes of idempotents in matrix algebras over $\mathcal A$, with addition given by $[e]+[f]=[e\oplus f]$.
We denote by $K_0(\mathcal A)$ its Grothendieck group.

\begin{customthm}{B}\label{thm B}
Let $T=uf(v)\in\mathcal R$, where $f\ne 0$ almost everywhere.
\begin{enumerate}
\item[$(1)$] If $f\in C(\mathbb T)$, then $K_0(\{T\}'\cap\mathcal A_\theta)\cong\mathbb Z$, where $\mathcal A_\theta=C^*(u,v)$.
\item[$(2)$] If $T$ is noninvertible and $r(T)=\|T\|$, then $K_0(\{T\}'\cap\mathcal R)\cong\mathbb Z$.
\item[$(3)$] Suppose that $T$ is invertible and
\begin{equation*}
	r(T)=\|T\|=b,\qquad r(T^{-1})=\|T^{-1}\|=a^{-1},\qquad 0<a<b.
\end{equation*}
Then $\{T\}'\cap\mathcal R$ is isomorphic to $H^\infty(\Omega)$, where $\Omega=\{z\in\mathbb C\colon a<|z|<b\}$.
In particular, $K_0(\{T\}'\cap\mathcal R)\cong\mathbb Z$.
\end{enumerate}
\end{customthm}

Parts $(1)$--$(3)$ of \Cref{thm B} follow from \Crefinorder{cor rotation-commutant-K0,cor H-infty-K0=Z,cor invertible-commutant-Hinfty}, respectively.
In \Cref{thm SI-K0-Z2}, we construct an invertible homogeneous operator $T=uf(v)$ such that for any $c_1,c_2\in\mathbb C\setminus\{0\}$, the operator $T_{c_1,c_2}=c_1(T\otimes 1)+c_2(1\otimes T)$ is strongly irreducible relative to $\mathcal R^{\otimes2}=\mathcal R\bar\otimes\mathcal R$ and satisfies
\begin{equation*}
    V(\{T_{c_1,c_2}\}'\cap\mathcal R^{\otimes2})
    \cong V(C(\mathbb T^2)),
    \qquad
    K_0(\{T_{c_1,c_2}\}'\cap\mathcal R^{\otimes2})
    \cong K_0(C(\mathbb T^2))
    \cong\mathbb Z^2.
\end{equation*}
Furthermore, for operators of the form $u^m+v^n$ with $m,n\in\mathbb Z\setminus\{0\}$, we study their strong irreducibility and classify them up to similarity in \Cref{subsec um-vn}, and prove a splitting result for the $K_0$-group in \Cref{cor monomial-sum-K0-split}.

Jones \cite{Jones-83} introduced the index as a fundamental invariant of inclusions of type $\mathrm{II}_1$ factors.
The authors in \cite{FJLX-13} identified the motivation of their work with the connection between strong irreducibility relative to type $\mathrm{II}_1$ factors and irreducible subfactors.
If $T\in\mathcal R$ is strongly irreducible relative to $\mathcal R$, then $W^*(T)$ is an irreducible subfactor of $\mathcal R$, and we call the index $[\mathcal R:W^*(T)]$ a \emph{strongly irreducible index} of $\mathcal R$, as defined in \Cref{def SI-index}.
Moreover, an operator $T\in\mathcal R$ is a \emph{strong generator} of $\mathcal R$ if every operator in its similarity orbit in $\mathcal R$ generates $\mathcal R$ as a von Neumann algebra, as defined in \Cref{def SG}.

\begin{customthm}{C}\label{thm C}
\begin{enumerate}
\item[$(1)$] There exists a homogeneous operator $T=uf(v)\in\mathcal R$ that is a strong generator of $\mathcal R$ and is strongly irreducible relative to $\mathcal M$ for every irreducible inclusion $\mathcal R\subseteq\mathcal M$ of type $\mathrm{II}_1$ factors.
\item[$(2)$] Every type $\mathrm{II}_1$ factor $\mathcal M$ with separable predual contains an operator that is strongly irreducible relative to $\mathcal M$.
\item[$(3)$] Every Jones index of an irreducible subfactor of $\mathcal R$ is a strongly irreducible index of $\mathcal R$.
\end{enumerate}
\end{customthm}

The strong generator in part $(1)$ is constructed in \Cref{thm SG}.
Parts $(2)$ and $(3)$ are applications of this construction and are proved in \Crefinorder{cor SG,cor SI-index}, respectively.
Moreover, we show in \Cref{cor SI-Boolean-lattice} that $\operatorname{Proj}(L^\infty(\mathbb T,\mu))$ can be realized as the relative invariant subspace lattice of the strong generator in part $(1)$.

\emph{Organization of the paper.}
In \Cref{sec pre}, we develop the Fourier methods used throughout the paper.
As an application, \Cref{rem Kaplansky} provides an affirmative answer to the Ulam--Kaplansky problem up to similarity for injective homogeneous operators of nonzero degree.
In \Cref{sec IR}, we study the irreducibility and reducible approximation of homogeneous operators.
In \Cref{sec commutant}, by introducing the notion of homogeneously power-free operators, we obtain the commutant structure of injective homogeneous operators in \Cref{thm commutant}.
In \Cref{sec SI}, we establish strong irreducibility in \Cref{thm SI-noninvertible}, describe roots of homogeneous operators in \Cref{prop SI-root}, and calculate the $K_0$-group of commutants relative to the irrational rotation algebra in \Cref{cor rotation-commutant-K0}.
In \Cref{sec bilateral}, we develop the bilateral weighted shift model.
The spectrum and point spectrum of a homogeneous operator are computed in \Cref{subsec spectrum}.
In \Cref{sec unilateral}, we develop the corresponding unilateral weighted shift model.
We also prove a version of Carleson's corona theorem in \Cref{prop corona}.
In \Cref{sec K0}, we study $K_0$-groups of relative commutants and splitting criteria for direct sums.
In \Cref{sec SG}, we construct strong generators and derive the existence and index results for strongly irreducible operators.
In \Cref{sec K0=Z2}, we also construct strongly irreducible operators whose relative commutants have $K_0$-groups isomorphic to $\mathbb Z^2$.
Finally, we provide a $V$-monoid criterion for the existence and uniqueness of finite strongly irreducible decompositions in \Cref{app CFJ}, with applications in \Cref{eg CFJ-applications}.

\tableofcontents

\section{Preliminaries}\label{sec pre}

Throughout, let $(\mathcal R,\tau)$ denote the hyperfinite type $\mathrm{II}_1$ factor generated by two unitary operators $u$ and $v$ satisfying
\begin{equation*}
	vu=\alpha uv,\qquad \alpha=e^{2\pi i\theta},\qquad \theta\in\mathbb R\setminus\mathbb Q.
\end{equation*}
Every operator $x$ in $\mathcal R$ can be viewed as a vector in $L^2(\mathcal R,\tau)$.
Write $L^0(\mathcal R,\tau)$ for the *-algebra of measurable operators affiliated with $\mathcal R$ and denote by $\mu$ the normalized Haar measure on the unit circle $\mathbb T$.
Note that a measurable operator in $L^0(\mathcal R,\tau)$ is invertible if and only if it has trivial kernel.

\subsection{\texorpdfstring{$L^2$}{L2}-Fourier expansions and homogeneous operators}

Let $\mathcal V=W^*(v)$ be the maximal abelian self-adjoint subalgebra of $\mathcal R$ generated by $v$, $e_{\mathcal V}$ the orthogonal projection from $L^2(\mathcal R,\tau)$ onto $L^2(\mathcal V,\tau)$, and $E_{\mathcal V}=e_{\mathcal V}|_{\mathcal R}$ the trace-preserving conditional expectation from $\mathcal R$ onto $\mathcal V$.
With respect to the orthogonal decomposition
\begin{equation*}
	L^2(\mathcal R,\tau)=\bigoplus_{n\in\mathbb Z}u^nL^2(\mathcal V,\tau),
\end{equation*}
every operator $x$ in $\mathcal R$ has a unique $L^2$-Fourier expansion
\begin{equation}\label{equ Fourier}
	x=\sum_{n\in\mathbb Z}u^nf_n(v),\qquad f_n(v)=E_{\mathcal V}(u^{-n}x),\qquad \|f_n\|_\infty\leqslant\|x\|,
\end{equation}
where $f_n\in L^\infty(\mathbb T,\mu)$ for every $n\in\mathbb Z$.
Moreover, if
\begin{equation*}
	x=\sum_{n\in\mathbb Z}u^nf_n(v),\qquad y=\sum_{n\in\mathbb Z}u^ng_n(v),
\end{equation*}
then the equality $h(v)u^n=u^nh(\alpha^nv)$ implies that
\begin{equation*}
	xy=\sum_{n\in\mathbb Z}u^n\sum_{k\in\mathbb Z}f_{n-k}(\alpha^kv)g_k(v),\qquad x^*=\sum_{n\in\mathbb Z}u^n\overline{f_{-n}}(\alpha^nv).
\end{equation*}
We introduce the concept of homogeneous operators in $L^0(\mathcal R,\tau)$.

\begin{definition}\label{def homogeneous}
A nonzero measurable operator $T\in L^0(\mathcal R,\tau)$ is \emph{$m$-homogeneous} if it is of the form $T=u^mf(v)$, where $m\in\mathbb Z$ and $f$ is a measurable function on $\mathbb T$.
We call $m$ the \emph{degree} of $T$ and $f$ its \emph{weight function}.
\end{definition}

In the following lemma, we provide an abstract characterization of homogeneous operators in $L^0(\mathcal R,\tau)$.

\begin{lemma}\label{lem abstract-homogeneous}
Let $\lambda\in\mathbb C$ and $T\in L^0(\mathcal R,\tau)$ such that
\begin{equation*}
	vTv^*=\lambda T.
\end{equation*}
If $\lambda=\alpha^m$ for some $m\in\mathbb Z$, then $T=u^mf(v)$ for some measurable function $f$ on $\mathbb T$.
Otherwise, $T=0$.
\end{lemma}

\begin{proof}
If $\lambda=\alpha^m$ and $T=u^mf(v)$ for some $m\in\mathbb Z$, then clearly $vTv^*=\lambda T$.
Suppose that $vTv^*=\lambda T$ and $T\ne0$.
Then $\lambda\ne0$.
Let $T=W|T|$ be the polar decomposition of $T$ and let $c=|\lambda|$.
By the uniqueness of the polar decomposition, we have
\begin{equation*}
    vWv^*=\frac{\lambda}{c}W,
    \qquad
    v|T|v^*=c|T|.
\end{equation*}
Let $H=|T|(I+|T|)^{-1}\in\mathcal R$.
Since $v|T|v^*=c|T|$, we have $vHv^*=c|T|(I+c|T|)^{-1}$ and
\begin{equation*}
    vHv^*-H
    =(c-1)|T|(I+c|T|)^{-1}(I+|T|)^{-1}.
\end{equation*}
It follows that $c=1$ by taking the trace on both sides.
Thus, we have
\begin{equation*}
    vWv^*=\lambda W,
    \qquad
    vHv^*=H.
\end{equation*}
Write $W=\sum_{n\in\mathbb Z}u^ng_n(v)\in\mathcal R$.
Then $(\alpha^n-\lambda)g_n=0$ for every $n\in\mathbb Z$.
Since $W\ne0$, we have $g_m\ne0$ for some $m\in\mathbb Z$, and hence $\lambda=\alpha^m$.
Consequently, we have
\begin{equation*}
    v(u^{-m}W)v^*=u^{-m}W,
    \qquad
    vHv^*=H.
\end{equation*}
Since $u^{-m}W,H\in\{v\}'\cap\mathcal R=W^*(v)$, we can write $u^{-m}W=g(v)$ and $H=h(v)$.
Consequently, $T=u^mf(v)$, where $f=gh(1-h)^{-1}$.
This completes the proof.
\end{proof}

We now give another description of \eqref{equ Fourier}.
For $z\in\mathbb T$, let $\pi_z$ be the automorphism of $\mathcal R$ determined by
\begin{equation}\label{equ pi-z}
    \pi_z(u)=zu,\qquad \pi_z(v)=v.
\end{equation}
It extends to a unitary operator on $L^2(\mathcal R,\tau)$, which we still denote by $\pi_z$.
For every $x=\sum_{n\in\mathbb Z}u^nf_n(v)\in L^2(\mathcal R,\tau)$ and $n\in\mathbb Z$, let
\begin{equation}\label{equ Pi-n}
    \Pi_n(x)=\int_{\mathbb T}z^{-n}\pi_z(x)\,d\mu(z).
\end{equation}
It is clear that $\Pi_n(x)=u^nf_n(v)$.
Therefore, the $L^2$-Fourier expansion \eqref{equ Fourier} for $x$ in $\mathcal R$ becomes
\begin{equation*}
	x=\sum_{n\in\mathbb Z}\Pi_n(x),\qquad \Pi_n(x)\in u^nW^*(v),\qquad
	\|\Pi_n(x)\|\leqslant\|x\|.
\end{equation*}
The following intertwining equation for two homogeneous operators in $L^0(\mathcal R,\tau)$ will be used frequently in this paper.

\begin{lemma}\label{lem frequently}
Let $T_j=u^{m_j}f_j(v)\in L^0(\mathcal R,\tau)$ for $j=1,2$.
Then for any $x\in L^2(\mathcal R,\tau)$, we have $T_1x=xT_2$ if and only if
\begin{equation*}
	T_1x_n=x_{n+m_1-m_2}T_2\qquad\text{for every}~n\in\mathbb Z,
\end{equation*}
where $x_n=\Pi_n(x)$ is the $n$-th Fourier component of $x$ given by \eqref{equ Pi-n}.
\end{lemma}

\begin{proof}
Let $A=(1+|f_2(v)|)^{-1}$ and $B=u^{m_1}(1+|f_1(v)|)^{-1}u^{-m_1}$.
Then
\begin{equation*}
	A,B\in W^*(v)\qquad\text{and}\qquad BT_1,T_2A\in\mathcal R.
\end{equation*}
Moreover, both $A$ and $u^{-m_1}B$ are invertible in $L^0(\mathcal R,\tau)$.
Since
\begin{equation*}
	u^{-m_1}B(T_1x_n-x_{n+m_1-m_2}T_2)A\in u^nL^2(W^*(v),\tau)
\end{equation*}
and we have the following $L^2$-Fourier expansion
\begin{equation*}
	u^{-m_1}B(T_1x-xT_2)A=\sum_{n\in\mathbb Z}u^{-m_1}B(T_1x_n-x_{n+m_1-m_2}T_2)A,
\end{equation*}
it follows that $T_1x=xT_2$ if and only if $T_1x_{n}=x_{n+m_1-m_2}T_2$ for every $n\in\mathbb Z$.
\end{proof}

Next we determine the norm continuity of the map $z\mapsto\pi_z(x)$, where $\pi_z$ is defined in \eqref{equ pi-z} and $x\in\mathcal R$.

\begin{lemma}\label{lem pi-z-continuous}
Let $x\in\mathcal R$.
Then $x\in C^*(u,W^*(v))$ if and only if the map $z\mapsto\pi_z(x)$ is norm continuous.
\end{lemma}

\begin{proof}
Let $\mathfrak A=C^*(u,W^*(v))$.
Since $\pi_z(u)=zu$ and $\pi_z$ fixes $W^*(v)$ pointwise, the map $z\mapsto\pi_z(x)$ is norm continuous for every operator in the *-algebra generated by $u$ and $W^*(v)$.
Because each $\pi_z$ is isometric, this continuity extends to every operator in $\mathfrak A$.
Conversely, suppose that $z\mapsto\pi_z(x)$ is norm continuous.
For every $\varepsilon>0$, there exists an open neighborhood $U$ of $1$ in $\mathbb T$ such that $\|\pi_z(x)-x\|<\varepsilon$ for every $z\in U$.
Let $h$ be a nonnegative continuous function on $\mathbb T$ such that
\begin{equation*}
    \operatorname{supp}(h)\subseteq U,
    \qquad
    \int_{\mathbb T}h(z)\,d\mu(z)=1.
\end{equation*}
The norm continuity of $z\mapsto\pi_z(x)$ ensures that $y_1=\int_{\mathbb T}h(z)\pi_z(x)\,d\mu(z)$ is a Bochner integral in $\mathcal R$ and
\begin{equation*}
    \|y_1-x\|
    \leqslant
    \int_{\mathbb T}h(z)\|\pi_z(x)-x\|\,d\mu(z)
    <\varepsilon.
\end{equation*}
By the uniform density of trigonometric polynomials in $C(\mathbb T)$, we can choose
\begin{equation*}
    p(z)=\sum_{n=-N}^{N}c_nz^n,\qquad
    \|p-h\|_\infty<\frac{\varepsilon}{1+\|x\|}.
\end{equation*}
Then
\begin{equation*}
    y_2=\int_{\mathbb T}p(z)\pi_z(x)\,d\mu(z)
      =\sum_{n=-N}^{N}c_n\Pi_{-n}(x)
      \in\mathfrak A.
\end{equation*}
Moreover, we have
\begin{equation*}
    \|y_1-y_2\|
    \leqslant\|x\|\int_{\mathbb T}|p(z)-h(z)|\,d\mu(z)
    <\varepsilon.
\end{equation*}
It follows that $\|x-y_2\|<2\varepsilon$.
Since $\varepsilon>0$ is arbitrary and $\mathfrak A$ is norm closed, we conclude that $x\in\mathfrak A$.
\end{proof}

\subsection{Fej\'er and Abel--Poisson approximations}\label{subsec 2-approximations}

Inspired by Shields' survey \cite[p.~90]{Shields-74}, for every vector $x\in L^2(\mathcal R,\tau)$ and every integer $N\geqslant 0$, we define
\begin{equation*}
    s_N(x)=\sum_{|n|\leqslant N}\Pi_n(x),
    \qquad
    \sigma_N(x)=\frac{s_0(x)+s_1(x)+\cdots+s_N(x)}{N+1},
\end{equation*}
where $\Pi_n(x)=\int_{\mathbb T}z^{-n}\pi_z(x)\,d\mu(z)$ as defined in \eqref{equ Pi-n}.
Then $\{s_N(x)\}_{N=0}^{\infty}$ are the \emph{partial sums} and $\{\sigma_N(x)\}_{N=0}^{\infty}$ are the \emph{Fej\'er sums} (also known as the \emph{Ces\`aro sums}) of $x$.
It is straightforward to verify that
\begin{equation}\label{equ Fejer}
    \sigma_N(x)
    =\sum_{|n|\leqslant N}
      \left(1-\frac{|n|}{N+1}\right)\Pi_n(x).
\end{equation}
Recall that the classical \emph{Fej\'er kernel} is defined by
\begin{equation}\label{equ Fejer-kernel}
    K_N(z)
    =\sum_{|n|\leqslant N}
      \left(1-\frac{|n|}{N+1}\right)z^n
    =\frac{1}{N+1}\left|\sum_{k=0}^{N}z^k\right|^2,
    \qquad z\in\mathbb T.
\end{equation}
By definition, we have $K_N(z)\geqslant 0$ and $\int_{\mathbb T}K_N(z)d\mu(z)=1$.
We now establish the following Fej\'er approximation result.

\begin{proposition}\label{prop Fejer}
For every integer $N\geqslant0$, the restriction of $\sigma_N$ to $\mathcal R$ is a normal unital trace-preserving completely positive map from $\mathcal R$ into itself, and
\begin{enumerate}
\item[$(1)$] $\|\sigma_N(x)\|\leqslant\|x\|$ for every $x\in\mathcal R$;
\item[$(2)$] $\|\sigma_N(x)\|_2\leqslant\|x\|_2$ and $\|\sigma_N(x)-x\|_2\to 0$ for every $x\in L^2(\mathcal R,\tau)$;
\item[$(3)$] $\sigma_N(x)-x\to 0$ in the $*$-strong-operator topology for every $x\in\mathcal R$;
\item[$(4)$] for every $x\in\mathcal R$, $\|\sigma_N(x)-x\|\to 0$ if and only if $x\in C^*(u,W^*(v))$.
\end{enumerate}
All limits are taken as $N\to\infty$.
\end{proposition}

\begin{proof}
By \eqref{equ Pi-n} and \eqref{equ Fejer-kernel}, the Fej\'er sum \eqref{equ Fejer} becomes
\begin{equation*}
	\sigma_N(x)=\int_{\mathbb T}K_N(z)\pi_z(x)d\mu(z).
\end{equation*}
Thus, $\sigma_N$ is normal unital trace-preserving completely positive on $\mathcal R$.
Moreover,
\begin{equation*}
    \|\sigma_N(x)\|\leqslant\int_{\mathbb T}K_N(z)\|\pi_z(x)\|d\mu(z)
    =\|x\|\qquad\text{for every}~x\in\mathcal R.
\end{equation*}
Now let $x\in L^2(\mathcal R,\tau)$.
Similarly, we have $\|\sigma_N(x)\|_2\leqslant\|x\|_2$.
By orthogonality of its homogeneous components, we have
\begin{equation*}
    \|s_N(x)-x\|_2^2=\sum_{|n|>N}\|\Pi_n(x)\|_2^2\to 0,
\end{equation*}
and
\begin{equation*}
    \|\sigma_N(x)-x\|_2\leqslant\frac{1}{N+1}\sum_{k=0}^{N}\|s_k(x)-x\|_2\to 0.
\end{equation*}
If $x\in\mathcal R$, then $\|\sigma_N(x)\|\leqslant\|x\|$ for every $N\geqslant 0$.
Since the $L^2$-topology coincides with the $*$-strong-operator topology on bounded subsets of $\mathcal R$, we see that $\sigma_N(x)\to x$ in the $*$-strong-operator topology.

We prove (4) as follows.
Let $\mathfrak A=C^*(u,W^*(v))$.
It is clear that $\sigma_N(x)\in\mathfrak A$ for every $x\in\mathcal R$.
Thus, $\|\sigma_N(x)-x\|\to 0$ implies that $x\in\mathfrak A$.
Conversely, for every $x\in\mathfrak A$ and $\varepsilon>0$, there exists a finite sum $y=\sum_{|n|\leqslant m}\Pi_n(y)$ such that $\|x-y\|<\varepsilon$.
Moreover, $\|\Pi_n(y)\|\leqslant\|y\|$ for each $|n|\leqslant m$.
For $N\geqslant m$, we have
\begin{equation*}
    \sigma_N(y)-y
    =-\frac{1}{N+1}\sum_{|n|\leqslant m}|n|\Pi_n(y).
\end{equation*}
It follows that
\begin{equation*}
    \|\sigma_N(x)-x\|\leqslant\|\sigma_N(x-y)\|+\|\sigma_N(y)-y\|+\|x-y\|
    \leqslant 2\varepsilon+\frac{m(m+1)}{N+1}\|y\|.
\end{equation*}
Thus, $\limsup_{N\to\infty}\|\sigma_N(x)-x\|\leqslant2\varepsilon$ for every $\varepsilon>0$, and hence $\|\sigma_N(x)-x\|\to 0$.
\end{proof}

\begin{remark}
As a direct consequence of \Cref{prop Fejer}, we have $\|\sigma_N(x)-x\|\to 0$ for every $x\in C^*(u,v)$.
\end{remark}

More generally, for every $1\leqslant p<\infty$ and $x\in L^p(\mathcal R,\tau)$, we define the $n$-th Fourier component $\Pi_n(x)$ by \eqref{equ Pi-n}.
Then $\Pi_n(x)\in u^nL^p(\mathcal V,\tau)$ for every $n\in\mathbb Z$, and we write $x\sim\sum_{n\in\mathbb Z}\Pi_n(x)$ for its formal Fourier expansion.
The corresponding Fej\'er sums satisfy $\|\sigma_N(x)-x\|_p\to0$ as $N\to\infty$.
We omit further details.

For every vector $x\in L^2(\mathcal R,\tau)$, its \emph{Abel--Poisson sum} is defined by
\begin{equation}\label{equ Abel-Poisson}
	\Phi_r(x)=\sum_{n\in\mathbb Z}r^{|n|}\Pi_n(x),\qquad 0\leqslant r<1.
\end{equation}
In particular, $\Phi_0(x)=\Pi_0(x)$.
Moreover, for $x\in\mathcal R$, the series in \eqref{equ Abel-Poisson} converges absolutely in operator norm.
Recall that the classical \emph{Poisson kernel} is defined by
\begin{equation}\label{equ Poisson-kernel}
	P_r(z)=\sum_{n\in\mathbb Z}r^{|n|}z^n
	=\mathrm{Re}\left(\frac{1+rz}{1-rz}\right)
	=\frac{1-r^2}{|1-rz|^2}.
\end{equation}
By definition, we have $P_r(z)\geqslant 0$ and $\int_{\mathbb T}P_r(z)d\mu(z)=1$.
Similar to \Cref{prop Fejer}, we have the following Abel--Poisson approximation result.

\begin{proposition}\label{prop Abel-Poisson}
For every $0\leqslant r<1$, the restriction of $\Phi_r$ to $\mathcal R$ is a normal unital trace-preserving completely positive map from $\mathcal R$ into itself, and
\begin{enumerate}
\item[$(1)$] $\|\Phi_r(x)\|\leqslant\|x\|$ for every $x\in\mathcal R$;
\item[$(2)$] $\|\Phi_r(x)\|_2\leqslant\|x\|_2$ and $\|\Phi_r(x)-x\|_2\to 0$ for every $x\in L^2(\mathcal R,\tau)$;
\item[$(3)$] $\Phi_r(x)-x\to 0$ in the $*$-strong-operator topology for every $x\in\mathcal R$;
\item[$(4)$] for every $x\in\mathcal R$, $\|\Phi_r(x)-x\|\to 0$ if and only if $x\in C^*(u,W^*(v))$.
\end{enumerate}
All limits are taken as $r\to 1$.
\end{proposition}

\begin{proof}
By \eqref{equ Pi-n} and \eqref{equ Poisson-kernel}, the Abel--Poisson sum \eqref{equ Abel-Poisson} becomes
\begin{equation*}
	\Phi_r(x)=\int_{\mathbb T}P_r(z)\pi_z(x)\,d\mu(z),\qquad 0\leqslant r<1.
\end{equation*}
The remainder of the proof is similar to that of \Cref{prop Fejer}.
\end{proof}

Suppose that $x,y\in\mathcal R$ such that $\Pi_n(x)=\Pi_n(y)=0$ for every $n<0$.
Then
\begin{equation}\label{equ Pi-n-xy}
	\Pi_n(xy)=
	\begin{cases}
		\sum_{k=0}^n\Pi_{n-k}(x)\Pi_k(y), & n\geqslant 0,\\
		0, & n<0.
	\end{cases}
\end{equation}
It follows that $\Phi_r(xy)=\Phi_r(x)\Phi_r(y)$.
In particular, if $x$ is an idempotent in $\mathcal R$ with $\Pi_n(x)=0$ for every $n<0$, then $\Phi_r(x)$ is also an idempotent for every $0\leqslant r<1$.
This observation will be used in the calculation of $K_0$-groups.

\subsection{Power boundedness and similarity}

Recall that an operator $T$ is said to be \emph{power bounded} if $\sup_{n\geqslant 0}\|T^n\|<\infty$.
We say that $T$ is \emph{bilaterally power bounded} if $T$ is invertible and both $T$ and $T^{-1}$ are power bounded.

\begin{proposition}\label{prop power-bounded}
Let $T$ be a homogeneous operator in $\mathcal R$.
\begin{enumerate}
\item[$(1)$] $T$ is power bounded if and only if $T$ is similar to a contraction in $\mathcal R$, implemented by an invertible positive operator in $W^*(v)$;
\item[$(2)$] $T$ is bilaterally power bounded if and only if $T$ is similar to a unitary operator in $\mathcal{R}$, implemented by an invertible positive operator in $W^*(v)$.
\end{enumerate}
\end{proposition}

\begin{proof}
Let $T=u^mf(v)\in\mathcal{R}$ with degree $m\in\mathbb Z$ and weight $f\in L^\infty(\mathbb T,\mu)$.

(1) It is clear that every operator similar to a contraction is power bounded.
Suppose that $C=\sup_{n\geqslant 0}\|T^n\|<\infty$.
Let $F_0=1$ and
\begin{equation*}
    F_n(z)=\prod_{j=0}^{n-1}f(\alpha^{jm}z)\qquad\text{for}~n\geqslant 1.
\end{equation*}
Since $T^n=u^{mn}F_n(v)$, we have $\|F_n\|_\infty=\|T^n\|\leqslant C$.
Let
\begin{equation*}
    g(z)=\sup_{n\geqslant 0}|F_n(z)|.
\end{equation*}
Then $1\leqslant g(z)\leqslant C$.
In particular, $g(v)$ is a positive invertible operator in $W^*(v)$.
It is clear that
\begin{equation*}
    |f(z)|g(\alpha^mz)
    =\sup_{n\geqslant 0}|f(z)F_n(\alpha^mz)|
    =\sup_{n\geqslant 0}|F_{n+1}(z)|
    \leqslant g(z).
\end{equation*}
Consequently, $g(v)Tg(v)^{-1}=u^mf(v)g(\alpha^mv)g(v)^{-1}$ is a contraction in $\mathcal R$.

(2) The sufficiency is clear.
Suppose that $T$ is bilaterally power bounded, and let $C=\sup_{n\in\mathbb Z}\|T^n\|<\infty$.
For every $n\in\mathbb Z$, we write $T^n=u^{mn}F_n(v)$.
It is clear that $\|F_n\|_\infty=\|T^n\|\leqslant C$ and $F_{n+1}(z)=f(z)F_n(\alpha^mz)$.
Let
\begin{equation*}
    g(z)=\sup_{n\in\mathbb Z}|F_n(z)|.
\end{equation*}
Then $1\leqslant g(z)\leqslant C$.
In particular, $g(v)$ is a positive invertible operator in $W^*(v)$.
It is clear that
\begin{equation*}
    |f(z)|g(\alpha^mz)
    =\sup_{n\in\mathbb Z}|f(z)F_n(\alpha^mz)|
    =\sup_{n\in\mathbb Z}|F_{n+1}(z)|
    =g(z).
\end{equation*}
Consequently, $g(v)Tg(v)^{-1}=u^m(f/|f|)(v)$ is a unitary operator in $\mathcal R$.
\end{proof}

Recall from \cite{FK-52} that the \emph{Fuglede-Kadison determinant} $\Delta\colon\mathcal R\to[0,\infty)$ is defined by
\begin{equation*}
	\Delta(T)=\exp\tau(\log|T|),\qquad T\in\mathcal R,
\end{equation*}
with the convention $\exp(-\infty)=0$.
For a homogeneous operator $T=u^mf(v)$ in $\mathcal R$, one has
\begin{equation*}
	\Delta(T)=\exp\int_{\mathbb T}\log|f(z)|d\mu(z).
\end{equation*}
If $T$ is similar to a unitary operator in $\mathcal R$, then $\Delta(T)=1$ and $T$ is power bounded.
The converse also holds, as follows from the corollary below.
In contrast with \Cref{prop power-bounded} (2), we do not assume bilateral power boundedness.

\begin{corollary}\label{cor power-bounded}
Let $T$ be a homogeneous operator in $\mathcal R$.
If $\Delta(T)=1$ and $T$ is power bounded, then $T$ is similar to a unitary operator in $\mathcal R$.
\end{corollary}

\begin{proof}
By \Cref{prop power-bounded}, $T$ is similar to a contraction $S$ in $\mathcal R$.
Since
\begin{equation*}
	0=\tau(\log|T|)=\log\Delta(T)=\log\Delta(S)=\tau(\log|S|)\leqslant 0,
\end{equation*}
we have $|S|=1$.
Thus, $S$ is a unitary operator.
\end{proof}

\begin{remark}\label{rem power-bounded}
As an application of \Cref{cor power-bounded}, if $T$ is homogeneous and noninvertible in $\mathcal R$ and $\Delta(T)=1$, then $T$ cannot be power bounded.
In particular, $T=u^m(1+v)$ is not power bounded for every $m\in\mathbb Z$.
\end{remark}

Next we study unitary equivalence and similarity of homogeneous operators in $L^0(\mathcal R,\tau)$.
The following lemma gives a sufficient condition for similarity implemented by an invertible homogeneous operator in $\mathcal R$.

\begin{lemma}\label{lem 2-intertwining}
Let $T_j=u^{m_j}f_j(v)\in L^0(\mathcal R,\tau)$, where $m_j\in\mathbb Z$ and $f_j\ne0$ almost everywhere, for $j=1,2$.
Suppose that $m_1\ne0$.
If there exist $x,y\in\mathcal R$ such that
\begin{equation*}
    T_1x=xT_2,\qquad T_2y=yT_1,\qquad \tau(xy)\ne0,
\end{equation*}
then $m_1=m_2$ and $T_1w=wT_2$ for some invertible homogeneous operator $w\in\mathcal R$.
\end{lemma}

\begin{proof}
Let $d=m_1-m_2$, $x_n=\Pi_n(x)$, and $y_n=\Pi_n(y)$ for $n\in\mathbb Z$.
It follows from \Cref{lem frequently} that
\begin{equation}\label{equ 2-intertwining}
	T_1x_n=x_{n+d}T_2,\qquad
	T_2y_{-n}=y_{-n-d}T_1.
\end{equation}
Let $g_n(v)=x_ny_{-n}\in W^*(v)$ for every $n\in\mathbb Z$.
By \eqref{equ 2-intertwining}, we have
\begin{equation*}
    T_1g_n(v)=T_1x_ny_{-n}
    =x_{n+d}T_2y_{-n}
    =x_{n+d}y_{-n-d}T_1
    =g_{n+d}(v)T_1.
\end{equation*}
Since $T_1=u^{m_1}f_1(v)$ and $f_1\ne0$ almost everywhere, it follows that
\begin{equation}\label{equ 2-intertwining-again}
	g_n(z)=g_{n+d}(\alpha^{m_1}z)
	\qquad\text{almost everywhere}.
\end{equation}
In particular, $\|g_n\|_1=\|g_{n+d}\|_1$.
The Cauchy--Schwarz inequality implies that
\begin{equation*}
	\sum_{n\in\mathbb Z}\|g_n\|_1
	\leqslant\sum_{n\in\mathbb Z}\|x_n\|_2\|y_{-n}\|_2
	\leqslant\|x\|_2\|y\|_2<\infty.
\end{equation*}
Moreover, we have
\begin{equation*}
	\sum_{n\in\mathbb Z}\tau(g_n(v))=\lim_{N\to\infty}\sum_{|n|\leqslant N}\tau(x_ny_{-n})
	=\lim_{N\to\infty}\sum_{|n|\leqslant N}\tau(x_ny)=\tau(xy)\ne0.
\end{equation*}
Thus, $g_n$ cannot be all zero, and it follows from $\|g_n\|_1=\|g_{n+d}\|_1$ for all $n\in\mathbb Z$ that $d=0$, that is, $m_1=m_2=m$.
It follows that $g_n(z)=g_n(\alpha^mz)$ by \eqref{equ 2-intertwining-again}.
Since $m\ne 0$, each function $g_n$ is essentially constant.
Consequently, $x_ny_{-n}=c1$ for some $n\in\mathbb Z$ and $c\ne 0$.
Therefore, $x_n$ is invertible and homogeneous in $\mathcal R$ with inverse $c^{-1}y_{-n}$.
Since $d=0$, \eqref{equ 2-intertwining} gives $T_1x_n=x_nT_2$.
This completes the proof.
\end{proof}

The following proposition characterizes unitary equivalence and similarity of homogeneous operators in terms of their degrees and weight functions.
In both cases, the implementing operator in $\mathcal R$ can be chosen to be homogeneous.

\begin{proposition}\label{prop equivalence}
Let $T_j=u^{m_j}f_j(v)\in L^0(\mathcal R,\tau)$, where $m_j\in\mathbb Z$ and $f_j\ne0$ almost everywhere, for $j=1,2$.
Suppose that $m_1\ne0$.
\begin{enumerate}
\item[$(1)$] $T_1w=wT_2$ for some unitary operator $w$ in $\mathcal{R}$ if and only if $m_1=m_2=m$ and there exist $n\in\mathbb Z$ and $g\colon\mathbb T\to\mathbb T$ measurable such that
\begin{equation*}
	f_2(z)=f_1(\alpha^nz)\frac{g(z)}{g(\alpha^mz)}.
\end{equation*}

\item[$(2)$] $T_1x=xT_2$ for some invertible operator $x$ in $\mathcal{R}$ if and only if $m_1=m_2=m$ and there exist $n\in\mathbb Z$ and $g,g^{-1}\in L^\infty(\mathbb T,\mu)$ such that
\begin{equation*}
	f_2(z)=f_1(\alpha^nz)\frac{g(z)}{g(\alpha^mz)}.
\end{equation*}
\end{enumerate}
\end{proposition}

\begin{proof}
\textbf{Necessity.}
For $j=1,2$, let
\begin{equation*}
	\varphi_j=\frac{f_j}{|f_j|},\qquad h_j=\frac{|f_j|}{1+|f_j|}, \qquad U_j=u^{m_j}\varphi_j(v),\qquad H_j=h_j(v).
\end{equation*}
Then $T_j=U_j\cdot |f_j|(v)$ is the polar decomposition of $T_j$.

(1) If $T_1w=wT_2$, then we have
\begin{equation*}
	U_1w\cdot w^*|f_1|(v)w=wU_2\cdot|f_2|(v).
\end{equation*}
By the uniqueness of the polar decomposition, we obtain that $U_1w=wU_2$ and $w^*|f_1|(v)w=|f_2|(v)$.
Hence $H_1w=wH_2$.
Let $d=m_1-m_2$ and $u^ng_n(v)=\Pi_n(w)$ for $n\in\mathbb Z$.
By \Cref{lem frequently} and $U_1w=wU_2$, we have
\begin{equation}\label{equ UW=WU}
	\varphi_1(\alpha^n z)g_n(z)=g_{n+d}(\alpha^{m_2}z)\varphi_2(z)
	\qquad\text{almost everywhere}.
\end{equation}
It follows that $\|g_n\|_2=\|g_{n+d}\|_2$.
Since $\sum_{n\in\mathbb Z}\|g_n\|_2^2=\|w\|_2^2=1<\infty$, it is clear that $d=0$, that is, $m_1=m_2=m$.
Consequently, $|g_n(z)|=|g_n(\alpha^mz)|$ by \eqref{equ UW=WU}.
Thus, each $|g_n|$ is essentially constant.
Choose a nonzero $g_n$ such that $|g_n|=c>0$ almost everywhere, and let $g(z)=c^{-1}g_n(z)\colon\mathbb T\to\mathbb T$.
By \Cref{lem frequently} and $H_1w=wH_2$, we have
\begin{equation*}
	\bigl(h_1(\alpha^n z)-h_2(z)\bigr)g(z)=0
    \qquad\text{almost everywhere}.
\end{equation*}
It follows that $h_1(\alpha^n z)=h_2(z)$, and hence $|f_1(\alpha^n z)|=|f_2(z)|$.
Combining this with \eqref{equ UW=WU}, we have $f_1(\alpha^n z)g(z)=f_2(z)g(\alpha^mz)$.
This proves the necessity in (1).

(2) Suppose that $T_1x=xT_2$ for some invertible $x\in\mathcal R$.
Then $T_2x^{-1}=x^{-1}T_1$ and $\tau(xx^{-1})=1$.
By \Cref{lem 2-intertwining}, $m_1=m_2=m$, and there exists an invertible homogeneous operator $u^ng(v)\in\mathcal R$ such that $T_1u^ng(v)=u^ng(v)T_2$.
It follows that $g,g^{-1}\in L^\infty(\mathbb T,\mu)$ and
\begin{equation*}
	f_1(\alpha^nz)g(z)=g(\alpha^mz)f_2(z)\qquad\text{almost everywhere}.
\end{equation*}
This proves the necessity in (2).

\textbf{Sufficiency.}
Suppose that $m_1=m_2=m$, $n\in\mathbb Z$, and
\begin{equation*}
	f_2(z)=f_1(\alpha^nz)\frac{g(z)}{g(\alpha^mz)}
	\qquad\text{almost everywhere}.
\end{equation*}
Let $x=u^ng(v)$.
Then
\begin{equation*}
	T_1x
	=u^{m+n}f_1(\alpha^nv)g(v)
	=u^{m+n}g(\alpha^mv)f_2(v)
	=xT_2.
\end{equation*}
In (1), $|g|=1$ and hence $x$ is unitary in $\mathcal R$.
In (2), $g,g^{-1}\in L^\infty(\mathbb T,\mu)$.
Thus, $x$ is invertible in $\mathcal R$.
This completes the proof.
\end{proof}

As a direct application of \Cref{prop equivalence}, we have the following result.

\begin{corollary}\label{cor scalar-similarity}
Let $T=u^mf(v)\in\mathcal R$, where $m\ne 0$ and $f\ne0$ almost everywhere on $\mathbb T$.
Suppose that $T$ is not invertible in $\mathcal R$.
Then for $\lambda\in\mathbb C$, the following are equivalent:
\begin{enumerate}
\item[$(1)$] $\lambda\in\{\alpha^{km}\colon k\in\mathbb Z\}$;
\item[$(2)$] $T$ and $\lambda T$ are unitarily equivalent in $\mathcal R$;
\item[$(3)$] $T$ and $\lambda T$ are similar in $\mathcal R$.
\end{enumerate}
\end{corollary}

\begin{proof}
$(1)\Rightarrow(2)$.
If $\lambda=\alpha^{mk}$ for some $k\in\mathbb Z$, then $v^kT=\lambda Tv^k$.
Hence $T$ and $\lambda T$ are unitarily equivalent in $\mathcal R$.

$(2)\Rightarrow(3)$.
It is immediate.

$(3)\Rightarrow(1)$.
Since the nonzero operator $T$ is similar to $\lambda T$ in $\mathcal R$, it is clear that $\lambda\ne0$.
By \Cref{prop equivalence}, there exists an invertible homogeneous operator $x=u^ng(v)\in\mathcal R$ such that $\lambda Tx=xT$.
Since $T^nx^{-m}$ is a homogeneous operator in $L^0(\mathcal R,\tau)$ of degree zero, we can write $h(v)=T^nx^{-m}$.
In particular, $h(v)\in W^*(v)$ if $n\geqslant 0$.
Then $h\ne 0$ almost everywhere and
\begin{equation*}
	Th(v)=T(T^nx^{-m})=\lambda^m(T^nx^{-m})T=\lambda^mh(v)T.
\end{equation*}
Consequently, we have
\begin{equation}\label{equ hT=Th}
    h(z)=\lambda^mh(\alpha^mz)
    \qquad\text{almost everywhere}.
\end{equation}
It follows that $|\lambda|=1$ and $|h|$ is essentially constant.
In particular, $h\in L^\infty(\mathbb T,\mu)$.
Comparing Fourier coefficients of $h$ in \eqref{equ hT=Th}, we have
\begin{equation*}
	\lambda^m=\alpha^{-jm},\qquad h(v)=cv^j
\end{equation*}
for some $j\in\mathbb Z$ and $c\ne 0$.
It follows that $T^n=cv^jx^m$ is invertible in $\mathcal R$.
However, $T$ is not invertible in $\mathcal R$.
Thus, $n=0$ and $x=g(v)$.
By $\lambda Tx=xT$, we have
\begin{equation*}
	\lambda g(z)=g(\alpha^mz).
\end{equation*}
Comparing Fourier coefficients of $g$, we have $\lambda=\alpha^{km}$ for some $k\in\mathbb Z$.
\end{proof}

\begin{remark}\label{rem scalar-similarity}
The noninvertibility assumption in \Cref{cor scalar-similarity} is essential.
Indeed, there exists an invertible operator $T=u^2f(v)\in\mathcal R$ such that $T$ and $\alpha T$ are similar, but not unitarily equivalent.
Let
\begin{equation*}
	f(z)=z\frac{2+z}{2+\alpha^2z},\qquad T=u^2f(v),\qquad x=2+v.
\end{equation*}
Then $T=x^{-1}(u^2v)x$.
Since $u^{-1}(u^2v)u=\alpha u^2v$, $T$ and $\alpha T$ are similar in $\mathcal R$.
Suppose, to the contrary, that $T$ and $\alpha T$ are unitarily equivalent.
By \Cref{prop equivalence}, we have
\begin{equation*}
    \alpha f(z)
    =f(\alpha^nz)\frac{g(z)}{g(\alpha^2z)}
    \qquad\text{almost everywhere}
\end{equation*}
for some $n\in\mathbb Z$ and some measurable function $g\colon\mathbb T\to\mathbb T$.
It follows that
\begin{equation*}
    |f(z)|=|f(\alpha^nz)|
    \qquad\text{almost everywhere}.
\end{equation*}
If $n\ne 0$, then $|f|$ is essentially constant.
However, $|f|$ is continuous and
\begin{equation*}
    |f(1)|=\frac{3}{|2+\alpha^2|}>1,
    \qquad
    |f(\alpha^{-2})|
    =\frac{|2+\alpha^{-2}|}{3}<1.
\end{equation*}
Therefore, $n=0$, and hence
\begin{equation*}
    \alpha g(\alpha^2z)=g(z)
    \qquad\text{almost everywhere}.
\end{equation*}
Comparing Fourier coefficients of $g$, we obtain that
$\alpha^{2k+1}=1$ for some $k\in\mathbb Z$.
That is a contradiction.
\end{remark}

Using \Cref{lem 2-intertwining}, we now show that similarity of finite direct sums of homogeneous operators is determined by the similarity classes of their summands, counted with multiplicity.

\begin{corollary}\label{cor similarity-direct-sum}
Suppose that $m_j\in\mathbb Z\setminus\{0\}$, $n_j\in\mathbb Z$, and $f_j,g_j\in L^\infty(\mathbb T,\mu)$ are nonzero almost everywhere, for $1\leqslant j\leqslant N$.
Then the following are equivalent:
\begin{enumerate}
\item[$(1)$] $\bigoplus_{j=1}^Nu^{m_j}f_j(v)$ and $\bigoplus_{j=1}^Nu^{n_j}g_j(v)$
are similar in $M_N(\mathcal R)$;
\item[$(2)$]
there exists a permutation $\pi$ of $\{1,\ldots,N\}$
such that $u^{m_j}f_j(v)$ and $u^{n_{\pi(j)}}g_{\pi(j)}(v)$ are similar for every $1\leqslant j\leqslant N$.
\end{enumerate}
\end{corollary}

\begin{proof}
Since $(2)\Rightarrow(1)$ is obvious, we will prove $(1)\Rightarrow(2)$.
Let $T_j=u^{m_j}f_j(v)$ and $S_j=u^{n_j}g_j(v)$ for $1\leqslant j\leqslant N$.
Choose an invertible operator $X=(x_{ij})_{1\leqslant i,j\leqslant N}$ in $M_N(\mathcal R)$ such that
\begin{equation*}
    \operatorname{diag}(T_1,\ldots,T_N)X
    =X\operatorname{diag}(S_1,\ldots,S_N),
\end{equation*}
and write $Y=X^{-1}=(y_{ij})_{1\leqslant i,j\leqslant N}$.
Then
\begin{equation*}
    T_ix_{ij}=x_{ij}S_j,\qquad
    S_jy_{ji}=y_{ji}T_i.
\end{equation*}
By \Cref{lem 2-intertwining}, $\tau(x_{ij}y_{ji})=0$ whenever $T_i$ and $S_j$ are not similar in $\mathcal R$.
For each similarity class $\mathcal T$ represented among the $T_i$, let
\begin{equation*}
    I=\{i\colon T_i\in\mathcal T\},\qquad
    J=\{j\colon S_j\in\mathcal T\}.
\end{equation*}
Taking traces in the diagonal entries of $XY=1$ and $YX=1$ implies that
\begin{align*}
	|I| &=\sum_{i\in I}\sum_{j=1}^N\tau(x_{ij}y_{ji})
	=\sum_{i\in I}\sum_{j\in J}\tau(x_{ij}y_{ji}),\\
	|J| &=\sum_{j\in J}\sum_{i=1}^N\tau(y_{ji}x_{ij})
	=\sum_{j\in J}\sum_{i\in I}\tau(y_{ji}x_{ij}).
\end{align*}
It follows that $|I|=|J|$.
Thus, every such class occurs with the same multiplicity
in both lists.
This completes the proof.
\end{proof}

\begin{remark}\label{rem Kaplansky}
Suppose that $T_1$ and $T_2$ are injective homogeneous operators in $\mathcal R$ with nonzero degree.
If $T_1\oplus T_1$ and $T_2\oplus T_2$ are similar in $M_2(\mathcal R)$, then $T_1$ and $T_2$ are similar in $\mathcal R$.
This gives an affirmative answer to the similarity version of the Ulam--Kaplansky problem for injective homogeneous operators of nonzero degree; see \cite{JMSZ-26} for the formulation of the problem and related results.
\end{remark}

\section{Irreducible operators and reducible approximations}\label{sec IR}

An operator $T$ in $\mathcal R$ is called \emph{irreducible} relative to $\mathcal R$ if the relative commutant $\{T\}'\cap\mathcal R$ contains no projections other than 0 and 1; equivalently, $\{T,T^*\}'\cap\mathcal R=\mathbb C1$.
Otherwise, $T$ is called \emph{reducible} relative to $\mathcal R$.
In \Cref{subsec IR}, we will determine the irreducibility of homogeneous operators in $\mathcal R$.
Moreover, we show that every operator in $C^*(u,W^*(v))$ is a norm limit of reducible ones in \Cref{subsec Red}.

\subsection{Irreducible operators}\label{subsec IR}

For a measurable function $f$ on $\mathbb T$, its zero set is denoted by
\begin{equation}\label{equ zero-set}
	Z(f)=\{z\in\mathbb T\colon f(z)=0\}.
\end{equation}
Clearly, $f\ne 0$ almost everywhere if and only if $\mu(Z(f))=0$.

\begin{lemma}\label{lem diffuse-commutant}
Let $T=u^mf(v)\in L^0(\mathcal R,\tau)$, where $m\in\mathbb Z$ and $f$ is measurable on $\mathbb T$.
If $\mu(Z(f))>0$, then $\{T,T^*\}'\cap\mathcal R$ is diffuse.
\end{lemma}

\begin{proof}
Let $\mathcal M=\{T,T^*\}'\cap\mathcal R$.
It suffices to prove that $\mathcal M$ contains a unital diffuse abelian von Neumann subalgebra $\mathcal A$.
If $m=0$, then we can take $\mathcal A=W^*(v)$.

Suppose that $m\ne0$.
Since the rotation $z\mapsto\alpha^mz$ is ergodic and $\mu(Z(f))>0$, the measurable function
\begin{equation*}
    \gamma\colon\mathbb T\to\mathbb Z_{\geqslant 0},\qquad \gamma(z)=\min\{n\geqslant0\colon\alpha^{mn}z\in Z(f)\}
\end{equation*}
is well defined almost everywhere.
Moreover, $\gamma(z)=0$ if and only if $z\in Z(f)$.
Let
\begin{equation*}
	h(z)=\alpha^{m\gamma(z)}z,\qquad \mathcal A=W^*(h(v)).
\end{equation*}
If $z\notin Z(f)$ and $\gamma(z)<\infty$, then
\begin{equation*}
    \gamma(\alpha^mz)=\gamma(z)-1,
    \qquad
    h(\alpha^mz)=h(z).
\end{equation*}
Consequently, $f(z)\bigl(h(\alpha^mz)-h(z)\bigr)=0$ almost everywhere.
It follows that
\begin{equation*}
    h(v)T
    =u^mh(\alpha^mv)f(v)
    =u^mf(v)h(v)
    =Th(v).
\end{equation*}
Since $h(v)$ is a unitary operator in $\mathcal R$, it also commutes with $T^*$, and hence $\mathcal A\subseteq\mathcal M$.
For every $\lambda\in\mathbb T$, the set $\{z\in\mathbb T\colon h(z)=\lambda\}$ is contained, modulo a null set, in the countable set $\{\alpha^{-mn}\lambda\colon n\geqslant0\}$.
Therefore, the spectral distribution of $h(v)$ has no atoms.
As a consequence, $\mathcal A$ is a diffuse von Neumann algebra.
Since $\mathcal M$ is finite and $\mathcal A$ is unital, $\mathcal M$ is also diffuse.
\end{proof}

The following theorem characterizes the triviality of the relative commutant of a homogeneous measurable operator and its adjoint.

\begin{theorem}\label{thm IR-criterion}
Let $T=u^mf(v)\in L^0(\mathcal R,\tau)$ and $\mathcal{M}=\{T,T^*\}'\cap\mathcal{R}$, where $m\in\mathbb Z$ and $f$ is a measurable function on $\mathbb T$.
Then $\mathcal M=\mathbb C1$ if and only if
\begin{equation*}
	(1)~m\ne 0,\quad (2)~f\ne 0~\text{almost everywhere},\quad(3)~|f|~\text{is not essentially constant}.
\end{equation*}
Otherwise, $\mathcal M$ is diffuse.
\end{theorem}

\begin{proof}
Let
\begin{equation*}
    S=T(1+T^*T)^{-1/2}=u^mg(v),
    \qquad
    g(z)=\frac{f(z)}{\sqrt{1+|f(z)|^2}}.
\end{equation*}
Since $T=S(1-S^*S)^{-1/2}$, we have
$\mathcal M=\{S,S^*\}'\cap\mathcal R$.
Moreover, $Z(g)=Z(f)$, and $|g|$ is essentially constant if and only if $|f|$ is.
Therefore, there is no loss of generality in assuming that $f\in L^\infty(\mathbb T,\mu)$.

If $m=0$, then $W^*(v)\subseteq\mathcal M$ and hence $\mathcal M$ is diffuse.
If $\mu(Z(f))>0$, then $\mathcal M$ is diffuse by \Cref{lem diffuse-commutant}.
If $m\ne0$ and $|f|=c>0$ almost everywhere, then $U=c^{-1}T$ is a Haar unitary operator in $\mathcal M$.
Thus, $\mathcal M$ contains the unital diffuse von Neumann subalgebra $W^*(U)$ and hence is diffuse.

Now we assume that the conditions (1)--(3) hold.
Suppose that $p=\sum_{n\in\mathbb Z}u^ng_n(v)$ is a projection in $\mathcal M$ and let
\begin{equation*}
    \varphi=\frac{f}{|f|},
    \qquad
    U=u^m\varphi(v).
\end{equation*}
Then $T=U\cdot|f|(v)$ is the polar decomposition of $T$.
Since $T$ commutes with $2p-1$, we have
\begin{equation*}
	T=(2p-1)U(2p-1)\cdot(2p-1)|f|(v)(2p-1).
\end{equation*}
By the uniqueness of the polar decomposition, $p$ commutes with both $U$
and $|f|(v)$.
Together with \Cref{lem frequently}, we have almost everywhere
\begin{equation}\label{equ Tp=pT}
    \varphi(\alpha^n z)g_n(z)
    =g_n(\alpha^m z)\varphi(z),\qquad \bigl(|f(\alpha^n z)|-|f(z)|\bigr)g_n(z)=0.
\end{equation}
The first equality of \eqref{equ Tp=pT} gives $|g_n(\alpha^m z)|=|g_n(z)|$ for every $n\in\mathbb Z$.
Hence $|g_n|$ is essentially constant.
If $g_n$ is not identically zero for some $n\ne 0$, then the second equality of \eqref{equ Tp=pT} gives $|f(\alpha^nz)|=|f(z)|$ almost everywhere.
It follows that $|f|$ is essentially constant, a contradiction.
Hence $g_n=0$ for every $n\ne0$.
By the first equality of \eqref{equ Tp=pT} again, we have $g_0(\alpha^m z)=g_0(z)$.
Thus, $g_0$ is essentially constant, and $p=g_0(v)$ is a scalar projection.
Therefore, $p=0$ or $p=1$, and hence $\mathcal M=\mathbb C1$.
\end{proof}

As an immediate application of \Cref{thm IR-criterion}, we obtain a complete criterion for the irreducibility of homogeneous operators in $\mathcal R$.
See \cite[Corollary 4.3]{ZFS-17} for a degree-one irreducibility result.

\begin{corollary}\label{cor IR-criterion}
Let $T=u^mf(v)\in\mathcal R$, where $m\in\mathbb Z$ and $f\in L^\infty(\mathbb T,\mu)$.
Then $T$ is irreducible relative to $\mathcal R$ if and only if
\begin{equation*}
	(1)~m\ne 0,\quad (2)~f\ne 0~\text{almost everywhere},\quad(3)~|f|~\text{is not essentially constant}.
\end{equation*}
\end{corollary}

\subsection{Density of reducible operators}\label{subsec Red}

By \Cref{lem pi-z-continuous}, \Cref{prop Fejer}, and \Cref{prop Abel-Poisson}, for every $x\in\mathcal R$, the following statements are equivalent:
\begin{enumerate}
\item[$(1)$] $x\in C^*(u,W^*(v))$;
\item[$(2)$] the map $z\mapsto\pi_z(x)$ is norm continuous;
\item[$(3)$] $\|\sigma_N(x)-x\|\to 0$ as $N\to\infty$;
\item[$(4)$] $\|\Phi_t(x)-x\|\to 0$ as $t\to 1$.
\end{enumerate}
For simplicity, let $\mathfrak A=C^*(u,W^*(v))$.
In this subsection, we show that every operator in $\mathfrak A$ is a limit of reducible ones in norm topology.

Recall that a function $f$ in $L^\infty(\mathbb T,\mu)$ is called a \emph{simple function} if there are scalars $\{c_k\}_{k=1}^n$ and disjoint measurable sets $\{\Omega_k\}_{k=1}^n$ such that $f=\sum_{k=1}^nc_k\chi_{\Omega_k}$; equivalently, its spectrum $\sigma(f)$ is a finite set.
We first establish two auxiliary lemmas.

\begin{lemma}\label{lem simple-function}
Let $f_1,\ldots,f_r$ be simple functions on $\mathbb T$.
Then for every integer $N\geqslant1$, there exist an integer $q>N+1$ and a measurable set $\Omega\subseteq\mathbb T$ such that
\begin{equation*}
    0<\mu(\Omega)<q^{-1},
    \qquad
    \Omega,\alpha\Omega,\ldots,\alpha^{q+N+1}\Omega~
    \text{are pairwise disjoint},
\end{equation*}
and for almost every $z\in\Omega$, we have
\begin{equation}\label{equ simple-function}
    f_j(\alpha^kz)=f_j(\alpha^{q+k}z)
    \qquad
    1\leqslant j\leqslant r,~0\leqslant k\leqslant N.
\end{equation}
\end{lemma}

\begin{proof}
The finitely many simple functions $f_j(\alpha^k\cdot)$, $1\leqslant j\leqslant r$ and $0\leqslant k\leqslant N$, are simultaneously constant on some measurable set $\Omega_1$ with $\mu(\Omega_1)>0$.
By recurrence, choose an integer $q>N+1$ such that $\mu(\Omega_1\cap\alpha^{-q}\Omega_1)>0$ and let $\Omega_2=\Omega_1\cap\alpha^{-q}\Omega_1$.
Since $z,\alpha^qz\in\Omega_1$ for every $z\in\Omega_2$, the equality \eqref{equ simple-function} holds almost everywhere on $\Omega_2$.
Since $1,\alpha,\ldots,\alpha^{q+N+1}$ are distinct, there exists $\Omega\subseteq\Omega_2$ with $\mu(\Omega)>0$ such that $\Omega,\alpha\Omega,\ldots,\alpha^{q+N+1}\Omega$ are pairwise disjoint.
Disjointness gives $(q+N+2)\mu(\Omega)\leqslant1$, and hence $0<\mu(\Omega)<1/q$.
The set $\Omega$ has the required properties.
\end{proof}

\begin{lemma}\label{lem almost-commuting-projection}
Let $f_1,\ldots,f_r$ be simple functions on $\mathbb T$.
Then for every integer $N\geqslant1$, there exists a nontrivial projection $P\in\mathfrak A$ such that
\begin{equation*}
    f_j(v)P=Pf_j(v)\qquad(1\leqslant j\leqslant r),
    \qquad
    \|uP-Pu\|=\sin\frac{\pi}{2N}.
\end{equation*}
\end{lemma}

\begin{proof}
Choose $q$ and $\Omega$ as in \Cref{lem simple-function}.
Let $E=\chi_\Omega(v)$ and
\begin{equation*}
    E_{ij}=u^iEu^{-j},
    \qquad 0\leqslant i,j\leqslant q+N+1.
\end{equation*}
Since $\Omega,\alpha\Omega,\ldots,\alpha^{q+N+1}\Omega$ are pairwise disjoint, $\{E_{ij}\}_{0\leqslant i,j\leqslant q+N+1}$ forms a system of matrix units in $\mathfrak A$ such that
\begin{equation*}
	E_{jj}\in W^*(v)\quad(0\leqslant j\leqslant q+N+1),\qquad
	uE_{ij}u^*=E_{i+1,j+1}\quad
	(0\leqslant i,j\leqslant q+N).
\end{equation*}
For every $0\leqslant k\leqslant N$, let $\theta_k=\frac{k\pi}{2N}$ and
\begin{equation*}
	P_k=\sin^2\theta_k\,E_{kk}+\sin\theta_k\cos\theta_k(E_{k,q+k}+E_{q+k,k})+\cos^2\theta_k\,E_{q+k,q+k}.
\end{equation*}
We define a projection $P$ by
\begin{equation*}
    P=\sum_{k=0}^N P_k
      +\sum_{k=N+1}^{q-1}E_{kk}.
\end{equation*}
Its summands are mutually orthogonal projections, and $\tau(P)=q\mu(\Omega)\in(0,1)$.
Thus, $P$ is a nontrivial projection in $\mathfrak A$.

For $1\leqslant j\leqslant r$, let $D_j=f_j(v)$ and
\begin{equation*}
    D_{j,k}=Eu^{-k}D_ju^kE=f_j(\alpha^kv)E,\qquad 0\leqslant k\leqslant q+N+1.
\end{equation*}
It follows from \eqref{equ simple-function} that $D_{j,k}=D_{j,q+k}$, and hence $D_j$ commutes with $P_k$ for every
$0\leqslant k\leqslant N$.
Clearly, $D_j$ also commutes with $E_{kk}$ for each $N+1\leqslant k\leqslant q-1$.
Thus, we obtain that $D_jP=PD_j$, that is, $f_j(v)P=Pf_j(v)$, for every $1\leqslant j\leqslant r$.

Finally, since $P_0=E_{qq}$ and $P_N=E_{NN}$, we have
\begin{equation*}
    uPu^*-P=\sum_{k=1}^N(uP_{k-1}u^*-P_k).
\end{equation*}
In the reduced subalgebra $(E_{kk}+E_{q+k,q+k})\mathcal R(E_{kk}+E_{q+k,q+k})$, we can write
\begin{equation*}
	uP_{k-1}u^*=
	\begin{pmatrix}
		\sin^2\theta_{k-1} & \sin\theta_{k-1}\cos\theta_{k-1}\\
		\sin\theta_{k-1}\cos\theta_{k-1} & \cos^2\theta_{k-1}
	\end{pmatrix},
\end{equation*}
and
\begin{equation*}
	P_k=
	\begin{pmatrix}
		\sin^2\theta_k & \sin\theta_k\cos\theta_k\\
		\sin\theta_k\cos\theta_k & \cos^2\theta_k
	\end{pmatrix}.
\end{equation*}
Therefore, $\|uP_{k-1}u^*-P_k\|=\sin(\theta_k-\theta_{k-1})=\sin\frac{\pi}{2N}$.
Since $\{E_{kk}+E_{q+k,q+k}\}_{k=1}^N$ are mutually orthogonal, we obtain that
\begin{equation*}
    \|uP-Pu\|=\|uPu^*-P\|=\sin\frac{\pi}{2N}.
\end{equation*}
This completes the proof.
\end{proof}

The following proposition shows that reducible operators belonging to $\mathfrak A$ are norm dense in $\mathfrak A$.

\begin{proposition}\label{prop reducible-approximation}
For every $x\in\mathfrak A$ and $\varepsilon>0$, there exist an operator $y\in\mathfrak A$ and a nontrivial projection $P\in\mathfrak A$ such that $\|x-y\|<\varepsilon$ and $yP=Py$.
\end{proposition}

\begin{proof}
Let $x\in\mathfrak A$ and $\varepsilon>0$.
By \Cref{prop Fejer} and the norm density of simple functions in $L^\infty(\mathbb T,\mu)$, we can choose
\begin{equation*}
    y_0=\sum_{|n|\leqslant M}u^nf_n(v),
    \qquad \|x-y_0\|<\frac{\varepsilon}{2},
\end{equation*}
where each $f_n$ is a simple function in $L^\infty(\mathbb T,\mu)$.
Choose an integer $N\geqslant1$ such that
\begin{equation*}
	\sin\frac{\pi}{2N}\sum_{|n|\leqslant M}|n|\|f_n\|_\infty<\frac{\varepsilon}{2}.
\end{equation*}
By \Cref{lem almost-commuting-projection}, there exists a nontrivial projection $P\in\mathfrak A$ commuting with every $f_n(v)$ such that $\|uP-Pu\|=\sin\frac{\pi}{2N}$.
Since $\|u^nP-Pu^n\|\leqslant |n|\|uP-Pu\|$ for every $n\in\mathbb Z$, we have
\begin{equation*}
    \|y_0P-Py_0\|
    \leqslant
    \sin\frac{\pi}{2N}
    \sum_{|n|\leqslant M}|n|\|f_n\|_\infty<\frac{\varepsilon}{2}.
\end{equation*}
Let $y=Py_0P+(1-P)y_0(1-P)\in\mathfrak A$.
Then $y$ commutes with the projection $P\in\mathfrak A$.
Since $2P-1$ is a unitary operator, we have
\begin{equation*}
	\|y_0-y\|=\|(y_0-y)(2P-1)\|=\|y_0P-Py_0\|<\frac{\varepsilon}{2}.
\end{equation*}
It follows that $\|x-y\|<\varepsilon$.
This completes the proof.
\end{proof}

\begin{remark}\label{rem reducible-approximation}
Although homogeneous operators in $\mathcal R$ satisfying the conditions in \Cref{thm IR-criterion} are irreducible relative to $\mathcal R$, \Cref{prop reducible-approximation} shows that every homogeneous operator in $\mathcal R$ is a norm limit of reducible operators.
\end{remark}

\section{Relative commutants}\label{sec commutant}

In this section, we determine the relative commutant $\{T\}'\cap\mathcal R$ of injective homogeneous operators $T\in L^0(\mathcal R,\tau)$ with nonzero degree.
For this purpose, we introduce homogeneously power-free operators in \Cref{subsec power-free} and use them to describe the relative commutant in \Cref{subsec commutant}.

\subsection{Homogeneously power-free operators}\label{subsec power-free}

Throughout this subsection, let $T=u^mf(v)\in L^0(\mathcal R,\tau)$, where $m\ne 0$ and $f\ne 0$ almost everywhere.
For every $n\in\mathbb Z$, let $\mathcal C_n(T)$ denote the vector space of all operators of the form $u^ng(v)\in L^0(\mathcal R,\tau)$ commuting with $T$, that is,
\begin{equation*}
	\mathcal{C}_n(T)=\{u^ng(v)\in L^0(\mathcal R,\tau)\colon Tu^ng(v)=u^ng(v)T\}.
\end{equation*}
Note that every nonzero operator in $\mathcal{C}_n(T)$ is $n$-homogeneous.
We show that $\mathcal C_n(T)$ is at most one dimensional for every $n\in\mathbb Z$.

\begin{lemma}\label{lem commutant-degree}
For every $n\in\mathbb Z$, we have $\dim\mathcal C_n(T)\leqslant 1$.
In particular, $\mathcal C_0(T)=\mathbb C1$.
If $\mathcal C_n(T)=\mathbb Cu^ng(v)\ne\{0\}$, then $u^ng(v)$ is invertible in $L^0(\mathcal R,\tau)$.
Moreover, the set
\begin{equation*}
	\Gamma(T)=\{n\in\mathbb Z\colon\mathcal C_n(T)\ne\{0\}\}
\end{equation*}
is an additive group such that $m\mathbb Z\subseteq\Gamma(T)\subseteq\mathbb Z$.
\end{lemma}

\begin{proof}
Suppose that $n\in\Gamma(T)$ and $u^ng(v)$ is a nonzero operator in $\mathcal C_n(T)$.
Then
\begin{equation}\label{equ commutant-degree}
    g(\alpha^mz)f(z)=f(\alpha^nz)g(z)
    \qquad\text{almost everywhere}.
\end{equation}
If $n=0$, then $g(\alpha^mz)=g(z)$ almost everywhere by \eqref{equ commutant-degree} and hence $g$ is essentially constant.
Thus, $\mathcal C_0(T)=\mathbb C1$ and $0\in\Gamma(T)$.

Since $f\ne0$ almost everywhere, it follows from \eqref{equ commutant-degree} that $Z(g)=\alpha^mZ(g)$ modulo a null set, where $Z(g)$ is the zero set of $g$ as in \eqref{equ zero-set}.
Since the rotation $z\mapsto\alpha^mz$ is ergodic, $\mu(Z(g))$ is either $0$ or $1$.
Since $u^ng(v)$ is nonzero, we have $\mu(Z(g))=0$.
It follows that $g\ne 0$ almost everywhere and $u^ng(v)$ is invertible in $L^0(\mathcal R,\tau)$.
For every operator $u^nh(v)\in\mathcal C_n(T)$, we have
\begin{equation*}
	(u^ng(v))^{-1}u^nh(v)=g(v)^{-1}h(v)\in\mathcal C_0(T)=\mathbb C1.
\end{equation*}
It follows that $u^nh(v)\in\mathbb Cu^ng(v)$ and hence $\mathcal C_n(T)=\mathbb Cu^ng(v)$.
Therefore, every nonzero $\mathcal C_n(T)$ has dimension one.

If $n_1\in\Gamma(T)$ and $u^{n_1}g_1(v)$ is a nonzero operator in $\mathcal C_{n_1}(T)$, then
\begin{equation*}
    (u^ng(v))^{-1}u^{n_1}g_1(v)
    \in\mathcal C_{n_1-n}(T)\setminus\{0\}.
\end{equation*}
Therefore, $n_1-n\in\Gamma(T)$.
Combining this with $0\in\Gamma(T)$, it is clear that $\Gamma(T)$ is a subgroup of $\mathbb Z$.
Since $m\in\Gamma(T)$ by $T\in\mathcal C_m(T)$, we have $m\mathbb Z\subseteq\Gamma(T)$.
\end{proof}

It is clear that $\mathcal C_n(XTX^{-1})=X\mathcal C_n(T)X^{-1}$ and $\Gamma(XTX^{-1})=\Gamma(T)$ for every invertible homogeneous operator $X$ in $L^0(\mathcal R,\tau)$.
Next we introduce the following definition based on \Cref{lem commutant-degree}.

\begin{definition}\label{def power-free}
Let $T=u^mf(v)\in L^0(\mathcal R,\tau)$, where $m\ne 0$ and $f\ne 0$ almost everywhere.
Its \emph{commutant degree group} $\Gamma(T)$ is the set of all integers $n\in\mathbb Z$ such that $\mathcal C_n(T)\ne\{0\}$.
We say that $T$ is \emph{homogeneously power-free} if $\Gamma(T)=m\mathbb Z$.
\end{definition}

By definition, if $T=u^{\pm1}f(v)\in L^0(\mathcal R,\tau)$ and $f\ne 0$ almost everywhere, then $T$ is homogeneously power-free.
The following lemma verifies the terminology ``homogeneously power-free''.

\begin{lemma}\label{lem power-free}
Let $T=u^mf(v)\in L^0(\mathcal R,\tau)$, where $m\ne 0$ and $f\ne 0$ almost everywhere.
Then $T$ is homogeneously power-free if and only if $T\ne S^d$ for any integer $d\geqslant 2$ and any homogeneous operator $S$ in $L^0(\mathcal R,\tau)$.
\end{lemma}

\begin{proof}
If $T=S^d$ for some integer $d\geqslant2$ and some $k$-homogeneous operator $S$, then $m=kd$ and $S$ commutes with $T$.
Thus, $k\in\Gamma(T)\setminus m\mathbb Z$, and hence $T$ is not homogeneously power-free.
Conversely, if $T$ is not homogeneously power-free, then by \Cref{lem commutant-degree}, there exists $k\in\mathbb Z$ such that $\Gamma(T)=k\mathbb Z$ and $km>0$.
Since $m\mathbb Z\subseteq\Gamma(T)$, we can write $m=kd$ for some $d\geqslant 2$.
Let $X$ be a nonzero operator in $\mathcal C_k(T)$.
By \Cref{lem commutant-degree}, there is a scalar $c\ne 0$ such that $X^d=cT$.
Choose $\lambda\in\mathbb C$ with $\lambda^d=c^{-1}$ and let $S=\lambda X$.
Then $S$ is $k$-homogeneous and $S^d=T$.
This completes the proof.
\end{proof}

The following proposition shows that every injective homogeneous operator with nonzero degree is a power of a homogeneously power-free operator.

\begin{proposition}\label{prop power-free-factorization}
Let $T=u^mf(v)\in L^0(\mathcal R,\tau)$, where $m\ne0$ and $f\ne0$ almost everywhere.
Then
\begin{equation*}
	T=S^d
\end{equation*}
for some positive integer $d\geqslant1$ and homogeneously power-free operator $S\in L^0(\mathcal R,\tau)$.
Moreover, $d$ is unique, $S$ is unique up to multiplication by a $d$-th root of unity, and $\Gamma(T)=\Gamma(S)$.
\end{proposition}

\begin{proof}
By \Cref{lem commutant-degree}, there is a unique $k\in\mathbb Z$ such that $\Gamma(T)=k\mathbb Z$ and $km>0$.
Let $d=m/k$, which is a positive integer.
As in the proof of \Cref{lem power-free}, there exists a $k$-homogeneous operator $S$ such that $S^d=T$.
Since every homogeneous operator commuting with $S$ commutes with $T$, we have
\begin{equation*}
    k\mathbb Z\subseteq\Gamma(S)\subseteq\Gamma(T)=k\mathbb Z.
\end{equation*}
Thus, $\Gamma(S)=k\mathbb Z$, and $S$ is homogeneously power-free.

To prove uniqueness, suppose that $T=S_1^{d_1}$ for some integer $d_1\geqslant1$ and homogeneously power-free operator $S_1$ of degree $k_1$.
Then $k_1d_1=m$ and $S_1\in\mathcal C_{k_1}(T)$.
It follows that $k_1\in\Gamma(T)=k\mathbb Z$.
Write $k_1=qk$, where $q$ is a positive integer.
Since both $S_1$ and $S^q$ are nonzero operators in $\mathcal C_{k_1}(T)$, by \Cref{lem commutant-degree}, we have $S_1=\lambda S^q$ for some scalar $\lambda\ne 0$.
It follows from \Cref{lem power-free} that $q=1$ as $S_1$ is homogeneously power-free.
Therefore, $k_1=k$, $d_1=d$, and $S_1=\lambda S$.
Moreover, since $S_1^d=S^d=T\ne0$ and $S_1=\lambda S$, it is clear that $\lambda^d=1$.
\end{proof}

We have the following immediate application of \Cref{prop power-free-factorization}.

\begin{corollary}\label{cor T^j-Gamma}
Let $T=u^mf(v)\in L^0(\mathcal R,\tau)$, where $m\ne0$ and $f\ne0$ almost everywhere.
Then for every $j\in\mathbb Z\setminus\{0\}$ and $n\in\mathbb Z$, we have
\begin{equation*}
	\mathcal C_n(T^j)=\mathcal C_n(T),\qquad \Gamma(T^j)=\Gamma(T).
\end{equation*}
In particular, $\{T^j\}'\cap L^2(\mathcal R,\tau)=\{T\}'\cap L^2(\mathcal R,\tau)$ and $\{T^j\}'\cap\mathcal R=\{T\}'\cap\mathcal R$.
\end{corollary}

\begin{proof}
By \Cref{prop power-free-factorization}, we write $T=S^d$, where $d\geqslant 1$, $S$ is homogeneously power-free, and $\Gamma(T)=\Gamma(S)$.
Clearly, $\mathcal C_n(S)\subseteq\mathcal C_n(T)$ for every $n\in\mathbb Z$.
Combining this with \Cref{lem commutant-degree}, we have $\mathcal C_n(S)=\mathcal C_n(T)$.
By applying \Cref{prop power-free-factorization} to $T^{|j|}=S^{|j|d}$, we similarly have $\mathcal C_n(S)=\mathcal C_n(T^{|j|})=\mathcal C_n(T^j)$ for all $n\in\mathbb Z$.
Therefore, $\mathcal C_n(T^j)=\mathcal C_n(T)$ for all $n\in\mathbb Z$ and hence $\Gamma(T^j)=\Gamma(T)$.
The last two equalities follow from \Cref{lem frequently}.
\end{proof}

\subsection{Criteria for homogeneous power-freeness}

The definition of homogeneous power-freeness involves measurable operators and can be difficult to verify directly.
The following lemma gives a sufficient condition in terms of the weight function.

\begin{lemma}\label{lem f-f-bar}
Let $T=u^mf(v)\in L^0(\mathcal R,\tau)$, where $m\ne0$ and $f\ne0$ almost everywhere.
If
\begin{equation*}
    \frac{f(z)}{\overline{f(z)}}=cz^k
    \qquad\text{almost everywhere}
\end{equation*}
for some constant $c\in\mathbb T$ and $k\in\mathbb Z$ with $\gcd(|m|,|k|)=1$, then $T$ is homogeneously power-free.
\end{lemma}

\begin{proof}
Let $n\in\Gamma(T)$ and choose a nonzero operator
$u^ng(v)\in\mathcal C_n(T)$.
By \Cref{lem commutant-degree}, $g\ne0$ almost everywhere.
Since $Tu^ng(v)=u^ng(v)T$, we have
\begin{equation*}
    g(\alpha^mz)f(z)=f(\alpha^nz)g(z)
    \qquad\text{almost everywhere}.
\end{equation*}
Let $h=g/\overline g$.
Dividing the preceding equality by its complex conjugate, we obtain that $h(\alpha^mz)=\alpha^{nk}h(z)$.
Since $|h|=1$ almost everywhere, $h$ is a unit vector in $L^2(\mathbb T,\mu)$.
Comparing Fourier coefficients gives $\alpha^{mj}=\alpha^{nk}$ for some $j\in\mathbb Z$.
It follows that $mj=nk$.
Since $\gcd(|m|,|k|)=1$, we have $m\mid n$.
It follows that $\Gamma(T)=m\mathbb Z$.
This completes the proof.
\end{proof}

\begin{remark}\label{rem f-f-bar}
If $g$ is a real-valued measurable function on $\mathbb T$ with $g\ne 0$ almost everywhere, then $u^m(1+v)g(v)$ is homogeneously power-free for every $m\ne 0$.
In particular, $u^m(1+v)$ is homogeneously power-free for $m\ne 0$.
\end{remark}

It is surprising that the homogeneous power-freeness of $T=u^mf(v)$ can depend on the arithmetic properties of the parameter $\theta$, even when $m$ and $f$ are fixed.
We present two technical lemmas concerning Diophantine approximation and use them to illustrate this phenomenon in \Cref{cor Diophantine}.

\begin{lemma}\label{lem Diophantine-1}
Let $T=u^m(1+\rho v)\in\mathcal R$, where $m\geqslant2$ and $0<\rho<1$.
If
\begin{equation}\label{equ Diophantine-1}
    \sum_{r=1}^{\infty}
    \frac{\rho^r}{r|\alpha^{mr}-1|}<\infty,
\end{equation}
then $\Gamma(T)=\mathbb Z$.
In particular, $T$ is not homogeneously power-free.
\end{lemma}

\begin{proof}
By \eqref{equ Diophantine-1}, we define a continuous function $H$ on $\mathbb T$ by
\begin{equation*}
    H(z)=\sum_{r=1}^{\infty}
    \frac{(-1)^{r+1}\rho^rz^r}
    {r(\alpha^{mr}-1)}.
\end{equation*}
Let $h=e^H$ and $x=h(v)$.
Then $x$ is an invertible operator in $W^*(v)$.
Since
\begin{equation*}
    \frac{h(\alpha^mz)}{h(z)}
    =e^{H(\alpha^mz)-H(z)}=
    \exp\left(
        \sum_{r=1}^{\infty}
        \frac{(-1)^{r+1}\rho^rz^r}{r}
    \right)
    =1+\rho z,
\end{equation*}
we have
\begin{equation*}
	T=u^m(1+\rho v)=u^mh(\alpha^mv)h(v)^{-1}=xu^mx^{-1}=(xux^{-1})^m.
\end{equation*}
It is clear that $xux^{-1}=uh(\alpha v)h(v)^{-1}$.
Since $xux^{-1}$ is $1$-homogeneous and commutes with $T$, we have $1\in\Gamma(T)$.
Therefore, $\Gamma(T)=\mathbb Z$.
\end{proof}

\begin{lemma}\label{lem Diophantine-2}
Let $T=u^m(1+\rho v)\in\mathcal R$, where $m\geqslant2$ and $0<\rho<1$.
If there are reduced rational numbers $\{p_j/q_j\}_{j=1}^{\infty}$ such that $m\mid q_j$, $q_j\to\infty$, and 
\begin{equation}\label{equ Diophantine-2}
    \left|\theta-\frac{p_j}{q_j}\right|<e^{-q_j^2},
\end{equation}
then $\Gamma(T)=m\mathbb Z$, that is, $T$ is homogeneously power-free.
\end{lemma}

\begin{proof}
Suppose to the contrary that $n\in\Gamma(T)\setminus m\mathbb Z$.
Choose a nonzero operator $u^ng(v)\in\mathcal C_n(T)$.
Then $g\ne0$ almost everywhere by \Cref{lem commutant-degree}.
Let
\begin{equation*}
    G(z)=\log|g(z)|,
    \qquad F(z)=\log|1+\rho z|.
\end{equation*}
Since $Tu^ng(v)=u^ng(v)T$, we have $G(\alpha^mz)-G(z)=F(\alpha^nz)-F(z)$, and hence
\begin{equation}\label{equ G-F}
    G(z)-G(\alpha^{mN}z)
    =\sum_{r=0}^{N-1}
    \bigl(F(\alpha^{mr}z)-F(\alpha^{mr+n}z)\bigr)\qquad\text{for}~N\geqslant 1.
\end{equation}
We define
\begin{equation*}
    \beta_j=e^{2\pi i p_j/q_j},
    \qquad \ell_j=q_j/m,
    \qquad k_j=\lfloor\rho^{-\ell_j}\rfloor,
    \qquad N_j=\ell_j k_j.
\end{equation*}
Since $m\nmid n$ and $\gcd(p_j,m)=1$, it is clear that $e^{2\pi i p_jn/m}$ is a nontrivial $m$-th root of unity for every $j\geqslant 1$.
Passing to a subsequence, we may assume that
\begin{equation*}
	\beta_j^{n\ell_j}=e^{2\pi i p_jn/m}=\zeta\ne 1\qquad\text{for all}~j\geqslant 1.
\end{equation*}
By the assumption \eqref{equ Diophantine-2}, we have
\begin{equation*}
    N_j^2|\alpha-\beta_j|
    \leqslant
    \ell_j^2\rho^{-2\ell_j}2\pi\left|\theta-\frac{p_j}{q_j}\right|
    \leqslant
    2\pi\left(\frac{q_j}{m}\right)^2
    \exp\left(-q_j^2+\frac{2q_j}{m}\log\frac1\rho\right)
    \to0.
\end{equation*}
Since $\beta_j^{mN_j}=\beta_j^{k_jq_j}=1$, it follows that $|\alpha^{mN_j}-1|\leqslant mN_j|\alpha-\beta_j|\to 0$.
Since $F$ is a smooth function on $\mathbb T$, by \eqref{equ G-F} we have the following uniform estimates
\begin{equation}\label{equ G-F-again}
	G(z)-G(\alpha^{mN_j}z)=
    \sum_{r=0}^{N_j-1}
    \bigl(F(\beta_j^{mr}z)-F(\beta_j^{mr+n}z)\bigr)+o(1).
\end{equation}
Since $\beta_j^m$ is a primitive $\ell_j$-th root of unity, it is clear that
\begin{equation*}
    \prod_{r=0}^{\ell_j-1}(1+\rho\beta_j^{mr}z)
    =1-(-\rho z)^{\ell_j}.
\end{equation*}
Using $N_j=k_j\ell_j$ and $\beta_j^{n\ell_j}=\zeta$,
we obtain from \eqref{equ G-F-again} that
\begin{align*}
    G(z)-G(\alpha^{mN_j}z)
    &=k_j\log\left|\frac{1-(-\rho z)^{\ell_j}}{1-\zeta(-\rho z)^{\ell_j}}\right|+o(1)\\
    &=k_j\rho^{\ell_j}\operatorname{Re}\bigl((\zeta-1)(-z)^{\ell_j}\bigr)+O(k_j\rho^{2\ell_j})+o(1)\\
    &=\operatorname{Re}\bigl((\zeta-1)(-z)^{\ell_j}\bigr)+o(1).
\end{align*}
Since $G$ is finite almost everywhere and $\alpha^{mN_j}\to1$, we have $G(z)-G(\alpha^{mN_j}z)\to0$ in measure.
However, the functions $\operatorname{Re}\bigl((\zeta-1)(-z)^{\ell_j}\bigr)$ have the same distribution for every $j\geqslant 1$.
Since $\zeta\ne1$, they cannot converge to zero in measure.
This contradicts the preceding uniform asymptotic formula.
Therefore, $\Gamma(T)=m\mathbb Z$.
\end{proof}

In the following corollary, we show that the homogeneous power-freeness of the operator $T=u^m(1+\rho v)$ depends on the choice of $\theta$.

\begin{corollary}\label{cor Diophantine}
There are irrational numbers $\theta_1$ and $\theta_2$ such that for every integer $m\geqslant 2$ and $0<\rho<1$, we have
\begin{equation*}
	\Gamma(u^m(1+\rho v))=
	\begin{cases}
		\mathbb Z, & \theta=\theta_1,\\
		m\mathbb Z, & \theta=\theta_2.
	\end{cases}
\end{equation*}
\end{corollary}

\begin{proof}
Let $\theta_1=\sqrt2$.
For every integer $q\geqslant1$, there exists an integer $p$ such that $\delta:=|q\sqrt2-p|\leqslant1/2$.
Since $|2q^2-p^2|\geqslant1$ and $q\sqrt2+p\leqslant4q$,
we have
\begin{equation*}
    \delta
    =\frac{|2q^2-p^2|}{q\sqrt2+p}
    \geqslant\frac1{4q}.
\end{equation*}
Consequently, $|e^{2\pi iq\sqrt2}-1|=2\sin(\pi\delta)\geqslant4\delta\geqslant q^{-1}$.
Taking $q=mr$, we obtain
\begin{equation*}
    \sum_{r=1}^{\infty}
    \frac{\rho^r}{r|e^{2\pi i mr\theta_1}-1|}
    \leqslant m\sum_{r=1}^{\infty}\rho^r
    =\frac{m\rho}{1-\rho}<\infty.
\end{equation*}
Thus, \Cref{lem Diophantine-1} gives
$\Gamma(u^m(1+\rho v))=\mathbb Z$ when $\theta=\theta_1$.

Define a sequence of positive integers $\{q_j\}_{j=1}^{\infty}$ by $q_1=3$ and $q_{j+1}=(j+1)q_j^{q_j^2}$ for every $j\geqslant 1$.
Let
\begin{equation*}
    \theta_2=\sum_{j=1}^{\infty}\frac1{q_j},
    \qquad
    p_j=q_j\sum_{r=1}^{j}\frac1{q_r}.
\end{equation*}
Since $q_j\mid q_{j+1}$, each $p_j$ is an integer, and
\begin{equation*}
    p_{j+1}=(j+1)q_j^{q_j^2-1}p_j+1.
\end{equation*}
Every prime divisor of $q_{j+1}$ divides $(j+1)q_j$ and therefore does not divide $p_{j+1}$.
Together with $p_1=1$, we have $\gcd(p_j,q_j)=1$ for every $j\geqslant1$.
Since $q_{j+1}\geqslant2q_j$ and $q_j\geqslant3$, we have
\begin{equation*}
    0<\theta_2-\frac{p_j}{q_j}
    \leqslant\frac2{q_{j+1}}
    =\frac2{(j+1)q_j^{q_j^2}}
    \leqslant q_j^{-q_j^2}
    <e^{-q_j^2}.
\end{equation*}
Thus, $\theta_2$ is irrational.
Since $j!\mid q_j$, we have $m\mid q_j$ whenever $j\geqslant m$.
By \Cref{lem Diophantine-2}, we obtain that $\Gamma(u^m(1+\rho v))=m\mathbb Z$ when $\theta=\theta_2$.
\end{proof}

\subsection{Relative commutants}\label{subsec commutant}

In analogy with \Cref{def power-free}, we introduce the following definition, which is useful for determining the relative commutant $\{T\}'\cap\mathcal R$ of injective homogeneous operators $T$ in $L^0(\mathcal R,\tau)$ with nonzero degree.

\begin{definition}\label{def Gamma-b}
Let $T=u^mf(v)\in L^0(\mathcal R,\tau)$, where $m\ne 0$ and $f\ne 0$ almost everywhere.
The \emph{bounded commutant degree semigroup} $\Gamma_b(T)$ is the set of all integers $n\in\mathbb Z$ such that $\mathcal C_n(T)\cap\mathcal R\ne\{0\}$.
\end{definition}

If $n_j\in\Gamma_b(T)$ and $u^{n_j}g_j(v)\in\mathcal C_{n_j}(T)\cap\mathcal R$ are nonzero operators for $j=1,2$, then the product $u^{n_1}g_1(v)u^{n_2}g_2(v)$ is a nonzero operator in $\mathcal C_{n_1+n_2}(T)\cap\mathcal R$ by \Cref{lem commutant-degree}.
Thus, $\Gamma_b(T)$ is an additive semigroup of $\mathbb Z$ containing $0$.

\begin{lemma}\label{lem semigroup}
Let $T=u^mf(v)\in L^0(\mathcal R,\tau)$, where $m\ne 0$ and $f\ne 0$ almost everywhere.
Then at least one of the following statements holds:
\begin{enumerate}
\item[$(1)$] $\Gamma_b(T)$ is a subsemigroup of $\mathbb Z_{\geqslant 0}$;
\item[$(2)$] $\Gamma_b(T)$ is a subsemigroup of $\mathbb Z_{\leqslant 0}$;
\item[$(3)$] $\Gamma_b(T)$ is a subgroup of $\mathbb Z$.
\end{enumerate}
\end{lemma}

\begin{proof}
If neither (1) nor (2) holds, then $\Gamma_b(T)$ contains both a positive integer $p$ and a negative integer $-q$.
Since $\Gamma_b(T)$ is an additive semigroup containing $0$, it follows that $pq,-pq\in\Gamma_b(T)$ and hence $pq\mathbb Z\subseteq\Gamma_b(T)$.
For every $n\in\Gamma_b(T)$, we have
\begin{equation*}
    -n=(pq-1)n-npq\in\Gamma_b(T).
\end{equation*}
Thus, $\Gamma_b(T)$ is a subgroup of $\mathbb Z$.
This completes the proof.
\end{proof}

By \Cref{prop power-free-factorization}, every injective homogeneous operator $T$ of nonzero degree is a power of a homogeneously power-free operator $S$.
It follows from \Cref{cor T^j-Gamma} that $\{T\}'\cap\mathcal R=\{S\}'\cap\mathcal R$.
Thus, it suffices to consider the relative commutant of homogeneously power-free operators.

\begin{theorem}\label{thm commutant}
Let $T=u^mf(v)\in L^0(\mathcal R,\tau)$, where $m\ne0$ and $f\ne0$ almost everywhere.
Then there exist an integer $d\geqslant 1$ and a homogeneously power-free operator $S$ of degree $k$ such that
\begin{equation*}
	T=S^d,\qquad m=kd.
\end{equation*}
Moreover, every operator $x$ in $\{T\}'\cap\mathcal R$ has a unique $L^2$-Fourier expansion
\begin{equation*}
    x=\sum_{\substack{n\in\mathbb Z\\kn\in\Gamma_b(S)}}a_nS^n,
    \qquad a_n\in\mathbb C.
\end{equation*}
\end{theorem}

\begin{proof}
The existence of $S$ is guaranteed by \Cref{prop power-free-factorization}.
According to \Cref{cor T^j-Gamma}, we have $\{T\}'\cap\mathcal R=\{S\}'\cap\mathcal R$.
Thus, it suffices to consider the relative commutant of $S$.
Since $S$ is homogeneously power-free, it follows from \Cref{lem commutant-degree} that $\Gamma_b(S)=\{kn\colon n\in\mathbb Z,~S^n\in\mathcal R\}$ and
\begin{equation*}
	\mathcal C_j(S)\cap\mathcal R=
	\begin{cases}
		\mathbb CS^n, & j=kn\in\Gamma_b(S),\\
		\{0\}, & \mbox{otherwise}.
	\end{cases}
\end{equation*}
Now we choose $x\in\{S\}'\cap\mathcal R$.
Then $\Pi_j(x)\in\mathcal C_j(S)\cap\mathcal R$ for each $j\in\mathbb Z$ by \Cref{lem frequently}.
The asserted expansion and its uniqueness now follow from the unique $L^2$-Fourier expansion of $x$.
\end{proof}

\begin{remark}\label{rem commutant}
The commutant $\{u^m(1+\rho v)\}'\cap\mathcal R$ for $m\geqslant 2$ and $0<\rho<1$ depends on the choice of $\theta$ by \Cref{cor Diophantine}.
\end{remark}

As a consequence of \Cref{prop Fejer} and \Cref{thm commutant}, the relative commutant $\{T\}'\cap\mathcal R$ is always abelian.

\begin{corollary}\label{cor commutant-maximal-abelian}
Let $T=u^mf(v)\in L^0(\mathcal R,\tau)$, where $m\ne0$ and $f\ne0$ almost everywhere.
If $\{T\}'\cap\mathcal R\ne\mathbb C1$, then $\{T\}'\cap\mathcal R$ is a maximal commutative subalgebra of $\mathcal R$.
\end{corollary}

\begin{proof}
By \Cref{thm commutant}, we may assume that $T$ is homogeneously power-free.
Let $\mathcal A=\{T\}'\cap\mathcal R$.
For any $x,y\in\mathcal A$ and $M,N\geqslant 1$, both Fej\'er sums $\sigma_M(x)$ and $\sigma_N(y)$ are finite linear combinations of integer powers of $T$.
Hence
\begin{equation*}
	\sigma_M(x)\sigma_N(y)=\sigma_N(y)\sigma_M(x).
\end{equation*}
Letting $M\to\infty$ and then $N\to\infty$, it follows that $xy=yx$ by \Cref{prop Fejer}.
Therefore, $\mathcal A$ is abelian.
Suppose that $\mathcal A\ne\mathbb C1$.
Then $T^j\in\mathcal R$ for some nonzero $j\in\mathbb Z$ by \Cref{thm commutant}.
Combining this with \Cref{cor T^j-Gamma}, we obtain that
\begin{equation*}
    \mathcal A'\cap\mathcal R
    \subseteq \{T^j\}'\cap\mathcal R
    =\{T\}'\cap\mathcal R
    =\mathcal A.
\end{equation*}
This completes the proof.
\end{proof}

\begin{remark}\label{rem power-free-commutant}
The relative commutant $\{T\}'\cap\mathcal R$ is maximal commutative whenever $T\in\mathcal R$.
In contrast, if $f(z)=(z+1)/(z-1)$ and $T=uf(v)$, then $T^n$ is unbounded for every $n\ne 0$.
Thus, $\Gamma_b(T)=\{0\}$ and $\{T\}'\cap\mathcal R=\mathbb C1$.
\end{remark}

Let $T$ be an injective homogeneous operator in $\mathcal R$ of nonzero degree.
If $T$ has no $d$-th roots in $\mathcal R$ for $d\geqslant 2$, one might expect that every operator in $\{T\}'\cap\mathcal R$ can be written as a formal power series of $T$.
However, as the following example illustrates, this expectation is not justified.

\begin{example}
Let
\begin{equation*}
	g(z)=\frac{(z-\alpha)(z-\alpha^{-1})}{z-1},\qquad S=ug(v),\qquad T=S^2.
\end{equation*}
Since $S^2=u^2\cdot\alpha(v-\alpha)(v-\alpha^{-2})$, $S^3=u^3\cdot\alpha^3(v-\alpha)(v-\alpha^{-1})(v-\alpha^{-3})$, and $S$ is homogeneously power-free, we have $\Gamma_b(T)=\{0,2,3,4,\ldots\}$.
By \Cref{thm commutant}, every operator $x$ in $\{T\}'\cap\mathcal R$ has a unique $L^2$-Fourier expansion
\begin{equation}\label{equ sum-S}
	x=a_0+\sum_{n=2}^{\infty}a_nS^n.
\end{equation}
In particular, $S^3$ cannot be represented by an $L^2$-convergent series of powers of $T$.
Moreover, $T$ has no $d$-th roots in $\mathcal R$ for $d\geqslant 2$.
Suppose to the contrary that $T=x^d$ for some $x\in\mathcal R$.
Since $x\in\{T\}'\cap\mathcal R$, it is of the form \eqref{equ sum-S}.
Comparing the Fourier components of $T=x^d$ gives $a_0^d=0$, and hence $a_0=0$.
Thus, the lowest nonzero degree of $T=x^d$ is at least $2d\geqslant 4$.
That is a contradiction.
\end{example}

\section{Strongly irreducible operators}\label{sec SI}

An operator $T$ in $\mathcal R$ is said to be \emph{strongly irreducible} relative to $\mathcal R$ if the \emph{relative commutant} $\{T\}'\cap\mathcal R$ contains no idempotents other than 0 and 1.
In this section, we give sufficient conditions for the strong irreducibility of homogeneous operators in $\mathcal R$.

\subsection{Strongly irreducible operators}

Based on \Cref{thm IR-criterion} and \Cref{cor IR-criterion}, every irreducible homogeneous operator is injective.
Accordingly, unless otherwise stated, we shall assume in the sequel that $T=u^mf(v)$ with $f\ne0$ almost everywhere.
In the following lemma, we give a sufficient condition for the absence of nontrivial idempotents in $\{T\}'\cap\mathcal R$ in terms of the bounded commutant degree semigroup $\Gamma_b(T)$ defined in \Cref{def Gamma-b}.

\begin{lemma}\label{lem SI}
Let $T=u^mf(v)\in L^0(\mathcal R,\tau)$, where $m\ne0$ and $f\ne0$ almost everywhere.
If $\Gamma_b(T)$ is not a subgroup of $\mathbb Z$, then the only idempotents in $\{T\}'\cap\mathcal R$ are $0$ and $1$.
\end{lemma}

\begin{proof}
By \Cref{lem semigroup}, we have either $\Gamma_b(T)\subseteq\mathbb Z_{\geqslant0}$ or $\Gamma_b(T)\subseteq\mathbb Z_{\leqslant0}$.
Without loss of generality, we assume that $\Gamma_b(T)\subseteq\mathbb Z_{\geqslant0}$.
Let $p$ be an idempotent in $\{T\}'\cap\mathcal R$.
By \Cref{thm commutant}, $p$ has a unique $L^2$-Fourier expansion
\begin{equation*}
	p=\tau(p)+\sum_{n=1}^{\infty}\Pi_n(p).
\end{equation*}
Since $p^2=p$, we have $\tau(p)^2=\tau(p)$, that is, $\tau(p)=0$ or $\tau(p)=1$.
Since $p$ is similar to a projection in $\mathcal R$, it follows that $p$ is either $0$ or $1$.
\end{proof}

If the weight function satisfies $f\in L^\infty(\mathbb T,\mu)$ and $f^{-1}\notin L^\infty(\mathbb T,\mu)$, then we have a very simple criterion for the strong irreducibility of $u^mf(v)$.

\begin{theorem}\label{thm SI-noninvertible}
Let $T=u^mf(v)\in\mathcal R$, where $m\ne0$ and $f\ne0$ almost everywhere.
If $f$ is not invertible in $L^{\infty}(\mathbb T,\mu)$, then $T$ is strongly irreducible relative to $\mathcal R$.
\end{theorem}

\begin{proof}
Since $T\in\mathcal C_m(T)\cap\mathcal R$ is nonzero, we have $m\in\Gamma_b(T)$.
On the other hand, \Cref{lem commutant-degree} gives
\begin{equation*}
    \mathcal C_{-m}(T)=\mathbb CT^{-1}.
\end{equation*}
Since $T$ is not invertible in $\mathcal R$, its inverse $T^{-1}\in L^0(\mathcal R,\tau)$ is unbounded.
It follows that $-m\notin\Gamma_b(T)$.
Therefore, $\Gamma_b(T)$ is not a subgroup of $\mathbb Z$, and \Cref{lem SI} implies that $T$ is strongly irreducible relative to $\mathcal R$.
\end{proof}

Having established strong irreducibility in the noninvertible case in \Cref{thm SI-noninvertible}, we now turn to invertible homogeneous operators.
In the following, we assume that $T=u^mf(v)\in\mathcal R$, where $m\ne0$ and $f$ is invertible in $L^\infty(\mathbb T,\mu)$.
Since both $T$ and $T^{-1}$ are bounded, we have $m,-m\in\Gamma_b(T)$.
Thus, $\Gamma_b(T)$ is a subgroup of $\mathbb Z$ by \Cref{lem semigroup}.
Since $m\in\Gamma_b(T)$, there are integers $k,d\in\mathbb Z$ such that
\begin{equation*}
	\Gamma_b(T)=k\mathbb Z,\qquad m=kd,\qquad d\geqslant 1.
\end{equation*}
Let $S=u^kg(v)$ be a nonzero operator in $\mathcal C_k(T)\cap\mathcal R$.
Then $\mathcal C_{-k}(T)=\mathbb CS^{-1}$ by \Cref{lem commutant-degree}.
Since $-k\in\Gamma_b(T)$, the inverse $S^{-1}$ is bounded, and hence $S$ is invertible in $\mathcal R$.
Up to a scalar multiple, we may assume that $\Delta(S)=1$, where $\Delta$ is the Fuglede--Kadison determinant.
Since $S^d$ and $T$ are nonzero operators in $\mathcal C_m(T)$, \Cref{lem commutant-degree} implies that $S^d=cT$ for some scalar $c\ne0$.
It follows from \Cref{cor T^j-Gamma} that
\begin{equation*}
    \{S\}'\cap\mathcal R=\{T\}'\cap\mathcal R,
    \qquad
    \Gamma_b(S)=\Gamma_b(T)=k\mathbb Z.
\end{equation*}
Consequently, $T$ is strongly irreducible relative to $\mathcal R$ if and only if $S$ is.
It therefore suffices to consider invertible homogeneous operators whose bounded commutant degree semigroup is the subgroup of $\mathbb Z$ generated by their degree and whose Fuglede--Kadison determinant is $1$.
The following lemma gives a criterion for strong irreducibility in this setting.

Let $E\subseteq\mathbb T$ be a measurable set and $\chi_E$ its characteristic function.
The \emph{Fourier coefficients} of $\chi_E$ are defined by
\begin{equation*}
	\widehat{\chi_E}(n)=\int_Ez^{-n}d\mu(z),\qquad n\in\mathbb Z,
\end{equation*}
and the \emph{Fourier series} $\sum_{n\in\mathbb Z}\widehat{\chi_E}(n)z^n$ converges to $\chi_E$ in $L^2(\mathbb T,\mu)$.
As in \eqref{equ Fejer}, for the formal series $p_E=\sum_{n\in\mathbb Z}\widehat{\chi_E}(n)T^n$, we define its Fej\'er sums
\begin{equation*}
    \sigma_N(p_E)=
    \sum_{\substack{n\in\mathbb Z\\|mn|\leqslant N}}
    \left(1-\frac{|mn|}{N+1}\right)
    \widehat{\chi_E}(n)T^n,\qquad N\geqslant 0.
\end{equation*}

\begin{lemma}\label{lem SI-invertible}
Let $T=u^mf(v)\in\mathcal R$, where $m\ne0$ and $f$ is an invertible function in $L^\infty(\mathbb T,\mu)$ such that $\Gamma_b(T)=m\mathbb Z$ and $\Delta(T)=1$.
Then the following statements are equivalent:
\begin{enumerate}
\item[$(1)$] $T$ is strongly irreducible relative to $\mathcal R$;
\item[$(2)$] for each measurable subset $E\subseteq\mathbb T$ with $0<\mu(E)<1$, the Fej\'er sums of the formal series $p_E=\sum_{n\in\mathbb Z}\widehat{\chi_E}(n)T^n$ satisfy that $\sup_{N\geqslant0}\|\sigma_N(p_E)\|=\infty$.
\end{enumerate}
\end{lemma}

\begin{proof}
For every $n\in\mathbb Z$, we have $\Delta(T^n)=\Delta(T)^n=1$.
It follows that $\|T^n\|_2\geqslant 1$ by Jensen's inequality.
Suppose first that $p$ is a nontrivial idempotent in $\{T\}'\cap\mathcal R$.
By \Cref{thm commutant} and $\Gamma_b(T)=m\mathbb Z$, the idempotent $p$ has a unique $L^2$-Fourier expansion
\begin{equation*}
    p=\sum_{n\in\mathbb Z}a_nT^n.
\end{equation*}
By orthogonality, we have
\begin{equation*}
    \sum_{n\in\mathbb Z}|a_n|^2
    \leqslant
    \sum_{n\in\mathbb Z}|a_n|^2\|T^n\|_2^2
    =\|p\|_2^2.
\end{equation*}
Hence there exists $\varphi\in L^2(\mathbb T,\mu)$ with Fourier series $\varphi(z)=\sum_{n\in\mathbb Z}a_nz^n$.
Comparing the homogeneous components of $p^2=p$ gives
\begin{equation*}
    \sum_{j\in\mathbb Z}a_ja_{n-j}=a_n
    \qquad\text{for every }n\in\mathbb Z,
\end{equation*}
where the series converge absolutely.
It follows that $\varphi^2=\varphi$.
Thus, $\varphi=\chi_E$ for some measurable set $E\subseteq\mathbb T$, $a_n=\widehat{\chi_E}(n)$ for every $n\in\mathbb Z$, and $p=p_E$.
Since $p$ is nontrivial, we have $0<\mu(E)<1$.
It follows from \Cref{prop Fejer} that
\begin{equation*}
    \sup_{N\geqslant0}\|\sigma_N(p_E)\|\leqslant\|p_E\|=\|p\|<\infty.
\end{equation*}
Conversely, suppose that $0<\mu(E)<1$ and $C=\sup_{N\geqslant0}\|\sigma_N(p_E)\|<\infty$.
Then orthogonality and Fatou's lemma imply that
\begin{equation*}
    \sum_{n\in\mathbb Z}|\widehat{\chi_E}(n)|^2\|T^n\|_2^2
    \leqslant C^2.
\end{equation*}
Therefore, $\sigma_N(p_E)$ converges to $p_E=\sum_{n\in\mathbb Z}\widehat{\chi_E}(n)T^n$ in $L^2(\mathcal R,\tau)$.
Since the Fej\'er sums $\{\sigma_N(p_E)\}_{N=0}^{\infty}$ form a bounded sequence in $\mathcal R$, we have $p_E\in\mathcal R$ with $\|p_E\|\leqslant C$. 
Moreover, the identity $\chi_E^2=\chi_E$ gives
\begin{equation*}
    \sum_{j\in\mathbb Z}\widehat{\chi_E}(j)\widehat{\chi_E}(n-j)=\widehat{\chi_E}(n),
    \qquad\text{for every}~n\in\mathbb Z.
\end{equation*}
Hence $p_E^2=p_E$.
Since $\tau(p_E)=\widehat{\chi_E}(0)=\mu(E)\in(0,1)$, $p_E$ is a nontrivial idempotent in
$\{T\}'\cap\mathcal R$.
This proves the equivalence.
\end{proof}

\begin{remark}\label{rem SI-noninvertible}
By \Cref{cor T^j-Gamma}, for every integer $j\geqslant 1$ and injective homogeneous operator $T\in\mathcal R$ of nonzero degree, $T^j$ is strongly irreducible relative to $\mathcal R$ if and only if $T$ is.
\end{remark}

The preceding lemma reduces strong irreducibility to the unboundedness of Fej\'er sums associated with nontrivial characteristic functions.
The following corollary gives a sufficient condition in terms of the operator norms and $L^2$-norms of the integer powers of $T$.

\begin{corollary}\label{cor SI-summability}
Let $T=u^mf(v)\in\mathcal R$, where $m\ne0$ and
$f$ is invertible in $L^\infty(\mathbb T,\mu)$ such that $\Delta(T)=1$.
If
\begin{equation*}
    \sum_{n=1}^{\infty}
    \frac1{\|T^n\|+\|T^{-n}\|+\|T^n\|_2^2+\|T^{-n}\|_2^2}
    <\infty,
\end{equation*}
then $T$ is strongly irreducible relative to $\mathcal R$.
\end{corollary}

\begin{proof}
By the argument preceding \Cref{lem SI-invertible}, there are integers $k,d\in\mathbb Z$ such that
\begin{equation*}
	\Gamma_b(T)=k\mathbb Z,\qquad m=kd,\qquad d\geqslant 1.
\end{equation*}
Let $S=u^kg(v)$ be a nonzero operator in $\mathcal C_k(T)\cap\mathcal R$.
Then $S$ is invertible in $\mathcal R$.
After multiplying $S$ by a suitable scalar, we may assume that $T=S^d$.
Together with \Cref{cor T^j-Gamma}, we have
\begin{equation*}
	\Gamma_b(S)=k\mathbb Z,\qquad\Delta(S)=1,\qquad \{S\}'\cap\mathcal R=\{T\}'\cap\mathcal R.
\end{equation*}
It therefore suffices to prove that $S$ is strongly
irreducible relative to $\mathcal R$.
Next we claim that
\begin{equation*}
    \sum_{n=1}^{\infty}
    \frac1{\|S^n\|+\|S^{-n}\|+\|S^n\|_2^2+\|S^{-n}\|_2^2}
    <\infty.
\end{equation*}
For an invertible operator $X\in\mathcal R$, let
\begin{equation*}
    D_n(X)=\|X^n\|+\|X^{-n}\|+\|X^n\|_2^2+\|X^{-n}\|_2^2.
\end{equation*}
Let $C=\max_{0\leqslant r<d}\{\|S^r\|,\|S^{-r}\|\}\geqslant 1$.
For every $j\geqslant1$ and $0\leqslant r<d$, we have
\begin{equation*}
    S^{jd}=S^{jd+r}S^{-r},
    \qquad
    S^{-jd}=S^{-(jd+r)}S^r.
\end{equation*}
It follows that $D_j(T)=D_j(S^d)=D_{jd}(S)\leqslant C^2D_{jd+r}(S)$.
Consequently, we have
\begin{equation*}
    \sum_{n=d}^{\infty}\frac1{D_n(S)}
    =\sum_{j=1}^{\infty}\sum_{r=0}^{d-1}
      \frac1{D_{jd+r}(S)}
    \leqslant dC^2\sum_{j=1}^{\infty}\frac1{D_j(T)}
    <\infty.
\end{equation*}
This proves the claim.
Let
\begin{equation*}
    A_n=\|S^n\|+\|S^{-n}\|,
    \qquad
    B_n=\|S^n\|_2^2+\|S^{-n}\|_2^2,
\end{equation*}
and
\begin{equation*}
    J_1=\{n\geqslant1\colon A_n\geqslant B_n\},
    \qquad
    J_2=\{n\geqslant1\colon A_n<B_n\}.
\end{equation*}
Then $\sum_{n\in J_1}A_n^{-1}<\infty$, $\sum_{n\in J_2}B_n^{-1}<\infty$, and we can define
\begin{equation*}
    K=2\sum_{n\in J_1}\frac1{A_n}
      +\left(\sum_{n\in J_2}\frac1{B_n}\right)^{1/2}<\infty.
\end{equation*}
Let $E\subseteq\mathbb T$ be a measurable set with $0<\mu(E)<1$.
Suppose that the characteristic function $\chi_E$ has Fourier series $\sum_{n\in\mathbb Z}a_nz^n$, where $a_n=\widehat{\chi_E}(n)$.
Then $a_{-n}=\overline{a_n}$ and $\sum_{n=1}^{\infty}|a_n|=\infty$.
Otherwise, $\chi_E$ would be a continuous function on $\mathbb T$, a contradiction.
As in \Cref{lem SI-invertible}, we define a formal series $p_E=\sum_{n\in\mathbb Z}a_nS^n$.
For simplicity, let
\begin{equation*}
    w_{n,N}=\max\left\{1-\frac{|kn|}{N+1},0\right\}.
\end{equation*}
Then the Fej\'er sums of $p_E$ are given by the finite sums
\begin{equation*}
    \sigma_{N}(p_E)=a_0+\sum_{n=1}^{\infty}w_{n,N}\bigl(a_nS^n+\overline{a_n}S^{-n}\bigr).
\end{equation*}
It follows that $w_{n,N}|a_n|\|S^{\pm n}\|\leqslant\|\sigma_N(p_E)\|$ and hence $w_{n,N}|a_n|A_n\leqslant2\|\sigma_N(p_E)\|$.
Thus, we have
\begin{equation}\label{equ J-1}
    \sum_{n\in J_1}w_{n,N}|a_n|
    \leqslant2\|\sigma_{N}(p_E)\|\sum_{n\in J_1}\frac1{A_n}.
\end{equation}
On the other hand, orthogonality gives
\begin{equation*}
    \sum_{n=1}^{\infty}w_{n,N}^2|a_n|^2B_n
    =\|\sigma_{N}(p_E)\|_2^2-|a_0|^2
    \leqslant\|\sigma_{N}(p_E)\|^2.
\end{equation*}
By the Cauchy--Schwarz inequality, we have
\begin{equation}\label{equ J-2}
    \sum_{n\in J_2}w_{n,N}|a_n|
    \leqslant
    \|\sigma_{N}(p_E)\|
    \left(\sum_{n\in J_2}\frac1{B_n}\right)^{1/2}.
\end{equation}
Combining \eqref{equ J-1} and \eqref{equ J-2}, we have
\begin{equation*}
    \sum_{n=1}^{\infty}w_{n,N}|a_n|
    \leqslant K\|\sigma_{N}(p_E)\|.
\end{equation*}
Since $w_{n,N}\to 1$ (as $N\to\infty$) for every $n\geqslant1$ and $\sum_{n=1}^{\infty}|a_n|=\infty$, it follows that $\|\sigma_{N}(p_E)\|\to\infty$.
We complete the proof by \Cref{lem SI-invertible}.
\end{proof}

\begin{remark}
Let $2<p<\infty$ and $1/p+1/q=1$.
Then for every $X\in\mathcal R$, we have
\begin{equation*}
    \|X\|_p^q
    \leqslant
    \|X\|^{(p-2)/(p-1)}
    \bigl(\|X\|_2^2\bigr)^{1/(p-1)}
    \leqslant
    \frac{p-2}{p-1}\|X\|
    +\frac1{p-1}\|X\|_2^2.
\end{equation*}
Therefore, under the remaining hypotheses of that corollary, the condition
\begin{equation*}
    \sum_{n=1}^{\infty}
    \frac1{\|T^n\|_p^q+\|T^{-n}\|_p^q}<\infty
\end{equation*}
also implies that $T$ is strongly irreducible relative
to $\mathcal R$.
\end{remark}

The preceding summability criterion yields a simple sufficient condition for strong irreducibility in terms of the spectral radii of $T$ and $T^{-1}$.
Note that for every invertible operator $T$, we have
\begin{equation*}
	r(T)r(T^{-1})=
	\frac{\max_{\lambda\in\sigma(T)}|\lambda|}{\min_{\lambda\in\sigma(T)}|\lambda|}.
\end{equation*}

\begin{corollary}\label{cor SI-spectral-radius}
Let $T=u^mf(v)\in\mathcal R$, where $m\ne0$ and $f$ is invertible in $L^\infty(\mathbb T,\mu)$.
If $r(T)r(T^{-1})>1$, then $T$ is strongly irreducible relative to $\mathcal R$.
\end{corollary}

\begin{proof}
Let $\rho=\Delta(T)>0$ and $S=\rho^{-1}T$.
Then $\Delta(S)=1$ and
\begin{equation*}
    r(S)r(S^{-1})=r(T)r(T^{-1})>1.
\end{equation*}
It follows that $R=r(S^\varepsilon)>1$ for some $\varepsilon\in\{1,-1\}$.
For every $n\geqslant1$, we have
\begin{equation*}
    \|S^{\varepsilon n}\|
    \geqslant r(S^{\varepsilon n})=R^n.
\end{equation*}
Consequently, we have
\begin{equation*}
    \sum_{n=1}^{\infty}
    \frac1{\|S^n\|+\|S^{-n}\|+\|S^n\|_2^2+\|S^{-n}\|_2^2}
    \leqslant
    \sum_{n=1}^{\infty}R^{-n}
    <\infty.
\end{equation*}
We complete the proof by \Cref{cor SI-summability}.
\end{proof}

\subsection{Strong irreducibility and similarity of \texorpdfstring{$u^m+v^n$}{u\string^m+v\string^n}} \label{subsec um-vn}

In this subsection, we show that $u^m+v^n$ is strongly irreducible whenever $mn\ne0$, and classify these operators up to similarity in $\mathcal R$.

First, we consider $u^m+v$, where $m\ne0$.
Let $U=u$ and $V=u^{-m}v$.
Then
\begin{equation*}
    \mathcal R=W^*(U,V),\qquad VU=\alpha UV.
\end{equation*}
Since $u^m+v=U^m(1+V)$, it follows from \Cref{thm SI-noninvertible} that $u^m+v$ is strongly irreducible relative to $\mathcal R$, which generalizes the result in \cite[Proposition 8.5]{FJLX-13}.
The same idea extends to $u^m+v^n$ by choosing suitable unitary generators of $\mathcal R$.

\begin{corollary}\label{cor SI-sum-power}
Let $m,n\in\mathbb Z$ with $mn\ne0$.
Then $u^m+v^n$ is strongly irreducible relative to $\mathcal R$.
\end{corollary}

\begin{proof}
Write $m=m_1d$ and $n=n_1d$, where $d=\gcd(|m|,|n|)\geqslant 1$.
Choose $r,s\in\mathbb Z$ such that $m_1r+n_1s=1$.
Let
\begin{equation*}
    U=u^sv^r,\qquad
    V=u^{-m_1}v^{n_1},\qquad
    N=m_1n_1d\ne0.
\end{equation*}
Since $(u^pv^q)^k=\alpha^{pq\frac{k(k-1)}2}u^{kp}v^{kq}$ for all $p,q,k\in\mathbb Z$, we have
\begin{equation*}
	u=\alpha^{n_1r(s+m_1-1)/2}U^{n_1}V^{-r},\qquad v=\alpha^{m_1s(r-n_1+1)/2}U^{m_1}V^s.
\end{equation*}
Thus, $\mathcal R=W^*(U,V)$.
Moreover, we have $VU=\alpha^{m_1r+n_1s}UV=\alpha UV$.
A direct computation shows that
\begin{equation}\label{equ U-N-F-V}
	u^m+v^n=U^NF(V),
\end{equation}
where $F(z)=\alpha^{Ns(r+n_1(d-1))/2}z^{-mr}\bigl(\alpha^{-N(d-1)/2}+z^d\bigr)$.
Since $N\ne 0$ and $F$ is not invertible in $L^\infty(\mathbb T,\mu)$, \Cref{thm SI-noninvertible} shows that $U^NF(V)$ in \eqref{equ U-N-F-V} is strongly irreducible relative to $\mathcal R$.
This completes the proof.
\end{proof}

Before studying similarity, we record a lemma on intertwining between homogeneous operators of the same degree.

\begin{lemma}\label{lem mutual-homogeneous-intertwining}
Let $T_j=u^mf_j(v)\in\mathcal R$, where $m\ne0$
and $f_j\ne0$ almost everywhere, for $j=1,2$.
Suppose that there are nonzero operators $X,Y\in\mathcal R$ such that
\begin{equation*}
    T_1X=XT_2,\qquad T_2Y=YT_1.
\end{equation*}
Then there are nonzero homogeneous operators $x,y\in\mathcal R$
and an integer $a\geqslant1$ such that
\begin{equation*}
    T_1x=xT_2,\qquad T_2y=yT_1,
    \qquad
    (xy)^{|m|}=T_1^a,\qquad (yx)^{|m|}=T_2^a.
\end{equation*}
\end{lemma}

\begin{proof}
We assume that $m>0$, since the case $m<0$ is analogous.
By \Cref{lem frequently}, we can choose nonzero homogeneous intertwiners
\begin{equation*}
    x=u^ng(v),\qquad y=u^kh(v).
\end{equation*}
For integers $r,s\geqslant0$, the operators $T_1^rx$ and $T_2^sy$ remain nonzero bounded intertwiners and have degrees $n+rm$ and $k+sm$, respectively.
Choosing $r,s$ sufficiently large and relabeling these operators, we may therefore assume that $n\geqslant1$ and $k\geqslant1$.

The identity $T_1x=xT_2$ gives
\begin{equation*}
    f_1(\alpha^nz)g(z)
    =
    g(\alpha^mz)f_2(z)
    \qquad\text{almost everywhere}.
\end{equation*}
Since $f_1\ne 0$ and $f_2\ne0$ almost everywhere, we obtain that $g\ne0$ almost everywhere.
The same argument applies to $h$.
Thus, both $x$ and $y$ are invertible in $L^0(\mathcal R,\tau)$.
Let
\begin{equation*}
    a=n+k\geqslant2,\qquad A=xy,\qquad B=yx.
\end{equation*}
Then $A$ commutes with $T_1$, $B$ commutes with $T_2$, and both have degree $a$.
Thus, $A^mT_1^{-a}$ is a nonzero element of $\mathcal C_0(T_1)=\mathbb C1$ by \Cref{lem commutant-degree}.
Up to a scalar multiple, we may assume that $A^m=T_1^a$.
Moreover, $A=xy$, $B=yx$, and $T_1x=xT_2$ imply that
\begin{equation*}
    x^{-1}Ax=B,\qquad x^{-1}T_1x=T_2.
\end{equation*}
It follows that $B^m=x^{-1}A^mx=x^{-1}T_1^ax=T_2^a$.
This completes the proof.
\end{proof}

We now apply the preceding lemma to scalar changes of the two summands in $u^m+v^n$.
For this family, the existence of nonzero bounded intertwiners in both directions is equivalent to similarity.

\begin{lemma}\label{lem mutual-phase-changes-sum-powers}
Let $m,n\in\mathbb Z\setminus\{0\}$, $\lambda_1,\lambda_2\in\mathbb C$, and
\begin{equation*}
    T=u^m+v^n,\qquad
    S=\lambda_1u^m+\lambda_2v^n.
\end{equation*}
Then the following statements are equivalent:
\begin{enumerate}
\item[$(1)$] $\lambda_1\in\alpha^{m\mathbb Z}$ and
$\lambda_2\in\alpha^{n\mathbb Z}$;
\item[$(2)$] $T$ and $S$ are similar in $\mathcal R$;
\item[$(3)$] there exist nonzero $x,y\in\mathcal R$ such that
$Tx=xS$ and $Sy=yT$.
\end{enumerate}
\end{lemma}

\begin{proof}
If $\lambda_1=\alpha^{k_1m}$ and
$\lambda_2=\alpha^{k_2n}$ for some integers $k_1,k_2\in\mathbb Z$, then
\begin{equation*}
    T u^{k_2}v^{-k_1}=u^{k_2}v^{-k_1}S.
\end{equation*}
Thus, $(1)\Rightarrow(2)$.
Clearly, $(2)\Rightarrow(3)$.
It remains to prove $(3)\Rightarrow(1)$.

Since $T$ is injective and $x\ne0$, we have $S\ne0$.
We adopt the notation in the proof of \Cref{cor SI-sum-power}.
As in \eqref{equ U-N-F-V}, write
\begin{equation*}
    T=U^NF(V),\qquad S=U^NG(V),
\end{equation*}
where
\begin{equation*}
    F(z)=a_1z^p(1+a_2z^d),\qquad
    G(z)=a_1z^p(\lambda_1+\lambda_2a_2z^d),
\end{equation*}
with $d=\gcd(|m|,|n|)\geqslant 1$ and $a_1,a_2\in\mathbb T$.
We assume that $N>0$, since the case $N<0$ is analogous.
By \Cref{lem mutual-homogeneous-intertwining}, we can choose nonzero homogeneous intertwiners $x=U^rg(V)$, $y=U^sh(V)$ in $\mathcal R$, and an integer $a\geqslant1$ such that
\begin{equation*}
    Tx=xS,\qquad Sy=yT,\qquad
    (xy)^N=T^a,\qquad (yx)^N=S^a.
\end{equation*}
Comparing degrees gives $r+s=a$.
Moreover, $x,y$ are invertible in $L^0(\mathcal R,\tau)$.

Since $T$ is not invertible in $\mathcal R$, $xy$ is not invertible.
Since $\mathcal R$ is finite, $yx$ is not invertible.
Thus, $S$ is not invertible, and hence $|\lambda_1|=|\lambda_2|>0$.
Moreover,
\begin{equation*}
    \Delta(T)^a
    =\Delta(xy)^N
    =\Delta(yx)^N
    =\Delta(S)^a.
\end{equation*}
Since $\Delta(T)=1$ and $\Delta(S)=\max\{|\lambda_1|,|\lambda_2|\}$, we obtain that $\lambda_1,\lambda_2\in\mathbb T$.

The identity $(xy)^N=T^a$ gives
\begin{equation*}
    T^a=xZ,\qquad Z=y(xy)^{N-1}\in\mathcal R.
\end{equation*}
Write $Z=U^{Na-r}q(V)$ with $q\in L^\infty(\mathbb T,\mu)$.
It follows from $T^a=xZ$ that
\begin{equation*}
    \prod_{j=0}^{a-1}F(\alpha^{jN}z)
    =
    g(\alpha^{Na-r}z)q(z)
    \qquad\text{almost everywhere}.
\end{equation*}
Replacing $z$ by $\alpha^{r-Na}z$ and using boundedness of $q$,
we conclude that
\begin{equation*}
    |g(z)|\geqslant\frac{|H(z)|}{\|q\|_\infty}
    \qquad\text{almost everywhere},
    \qquad
    H(z)=\prod_{j=0}^{a-1}F(\alpha^{r-Na+jN}z).
\end{equation*}
Moreover, the identity $Tx=xS$ gives
\begin{equation*}
    F(\alpha^rz)g(z)=g(\alpha^Nz)G(z)
    \qquad\text{almost everywhere}.
\end{equation*}
Combining this with the preceding lower bound for $|g(z)|$, we obtain
\begin{equation*}
    |H(z)F(\alpha^rz)|
    \leqslant\|q\|_\infty\|g\|_\infty|G(z)|
    \qquad\text{almost everywhere}.
\end{equation*}
By continuity, this inequality holds everywhere on $\mathbb T$.
Choose $z_0\in\mathbb T$ such that $G(z_0)=0$.
Then
\begin{equation*}
    \prod_{j=0}^{a}|F(\alpha^{r-Na+jN}z_0)|
    =|H(z_0)F(\alpha^rz_0)|=0.
\end{equation*}
Thus, $F(\alpha^{r-Nk}z_0)=0$ for some integer $0\leqslant k\leqslant a$.
Let $\ell=r-Nk$.
Then
\begin{equation*}
    \lambda_2=\lambda_1\alpha^{\ell d},\qquad
    G(z)=\eta F(\alpha^\ell z),\qquad
    \eta=\lambda_1\alpha^{-\ell p}\in\mathbb T.
\end{equation*}
Consequently, $S=\eta U^{-\ell}TU^\ell$.
Since $r-\ell=Nk$, the operator
\begin{equation*}
    \varphi(V)=T^{-k}xU^{-\ell}\in L^0(\mathcal R,\tau)
\end{equation*}
has degree zero.
The identity $Tx=xS$ gives
$T\varphi(V)=\eta\varphi(V)T$, and hence
\begin{equation*}
    \varphi(\alpha^Nz)=\eta^{-1}\varphi(z)
    \qquad\text{almost everywhere}.
\end{equation*}
Since $|\eta|=1$, ergodicity implies that $|\varphi|$ is essentially constant.
Thus, $\varphi\in L^\infty(\mathbb T,\mu)$.
Comparing Fourier coefficients gives $\varphi(z)=c_0z^j$ for some $j\in\mathbb Z$ and $c_0\ne0$.
Let $W=\varphi(V)U^\ell$.
Then $W$ is invertible in $\mathcal R$ and
\begin{equation*}
    TW=\eta\varphi(V)TU^\ell=WS.
\end{equation*}
Since $U,V$ are monomials in $u,v$, we have $W=c_1u^{k_2}v^{-k_1}$ for some integers $k_1,k_2\in\mathbb Z$
and $c_1\ne0$.
Comparing the two monomials in $TW=WS$, we obtain
\begin{equation*}
    \lambda_1=\alpha^{k_1m},\qquad
    \lambda_2=\alpha^{k_2n}.
\end{equation*}
This proves $(1)$.
\end{proof}

Using the preceding lemma, we now show that sums corresponding to distinct exponent pairs are not similar in $\mathcal R$.

\begin{proposition}\label{prop similarity-sum-powers}
Let $m_1,n_1,m_2,n_2\in\mathbb Z$.
Then $u^{m_1}+v^{n_1}$ and $u^{m_2}+v^{n_2}$ are similar in $\mathcal R$ if and only if $(m_1,n_1)=(m_2,n_2)$.
\end{proposition}

\begin{proof}
Let $T_j=u^{m_j}+v^{n_j}$ for $j=1,2$.
Suppose that $T_1x=xT_2$ for some invertible operator $x$ in $\mathcal R$.
The proof is divided into three cases.

\textbf{Case I.}
$(m_1,n_1)=(0,0)$.
Since $\tau(T_2)=\tau(T_1)=2$, we have $(m_2,n_2)=(0,0)$.

\textbf{Case II.}
Exactly one of $m_1$ and $n_1$ is zero.
Since $\tau(T_2)=\tau(T_1)=1$, we see that exactly one of $m_2$ and $n_2$ is zero.
In this case, we have
\begin{equation*}
    T_j=1+u^{m_j}v^{n_j},\qquad j=1,2.
\end{equation*}
Since $T_1x=xT_2$, it follows that $u^{m_1}v^{n_1}x=xu^{m_2}v^{n_2}$.
Let
\begin{equation*}
	x=\sum_{m,n\in\mathbb Z}c_{m,n}u^mv^n\in\mathcal R\subseteq L^2(\mathcal R,\tau).
\end{equation*}
Then $|c_{m,n}|=|c_{m+m_1-m_2,n+n_1-n_2}|$.
If $(m_1,n_1)\ne (m_2,n_2)$, then $c_{m,n}=0$ for all $m,n\in\mathbb Z$, a contradiction.
Thus, $(m_1,n_1)=(m_2,n_2)$.

\textbf{Case III.}
$m_1n_1\ne 0$.
Since $\tau(T_2)=\tau(T_1)=0$, we have $m_2n_2\ne0$.
For $z\in\mathbb T$, let $\rho_z$ be the automorphism of $\mathcal R$ determined by
\begin{equation*}
    \rho_z(u)=z^{n_1}u,\qquad
    \rho_z(v)=z^{m_1}v.
\end{equation*}
Then $\rho_z(T_1)=z^{m_1n_1}T_1$, $z^{-m_1n_1}\rho_z(T_2)=z^{n_1(m_2-m_1)}u^{m_2}+z^{m_1(n_2-n_1)}v^{n_2}$, and
\begin{equation*}
	z^{m_1n_1}T_1\rho_z(x)=\rho_z(T_1)\rho_z(x)=\rho_z(T_1x)=\rho_z(xT_2)=\rho_z(x)\rho_z(T_2).
\end{equation*}
It follows that
\begin{equation*}
	T_2\cdot x^{-1}\rho_z(x)=x^{-1}\rho_z(x)\cdot z^{-m_1n_1}\rho_z(T_2).
\end{equation*}
Therefore, for every $z\in\mathbb T$, the operators
\begin{equation*}
	u^{m_2}+v^{n_2}\qquad\text{and}\qquad z^{n_1(m_2-m_1)}u^{m_2}+z^{m_1(n_2-n_1)}v^{n_2}
\end{equation*}
are similar in $\mathcal R$.
By \Cref{lem mutual-phase-changes-sum-powers}, we have
\begin{equation*}
    n_1(m_2-m_1)=0,\qquad
    m_1(n_2-n_1)=0.
\end{equation*}
Since $m_1n_1\ne0$, we obtain that $(m_1,n_1)=(m_2,n_2)$.
\end{proof}

\subsection{Roots}

In this subsection, we study the $d$-th roots of strongly irreducible homogeneous operators in $\mathcal R$ for $d\geqslant 2$.
We first record an elementary algebraic fact that will be used in the proof of \Cref{prop SI-root}.

\begin{lemma}\label{lem finite-order-scalar}
Let $\mathcal A$ be a unital complex algebra whose only idempotents are $0$ and $1$.
If $V\in\mathcal A$ satisfies $V^d=1$ for some integer $d\geqslant2$, then $V=\lambda1$ for some scalar $\lambda\in\mathbb C$ with $\lambda^d=1$.
\end{lemma}

\begin{proof}
Let $\Lambda_d=\{\lambda\in\mathbb C\colon\lambda^d=1\}$.
For every $\lambda\in\Lambda_d$, let
\begin{equation*}
    p_\lambda(t)=
    \prod_{\eta\in\Lambda_d\setminus\{\lambda\}}
    \frac{t-\eta}{\lambda-\eta},
    \qquad P_\lambda=p_\lambda(V).
\end{equation*}
By the polynomial identities $1=\sum_{\lambda\in\Lambda_d}p_\lambda(t)$ and $t=\sum_{\lambda\in\Lambda_d}\lambda p_\lambda(t)$, we obtain that
\begin{equation*}
	1=\sum_{\lambda\in\Lambda_d}P_\lambda,\qquad V=\sum_{\lambda\in\Lambda_d}\lambda P_\lambda.
\end{equation*}
Since $p_\lambda^2-p_\lambda$ vanishes on $\Lambda_d$, there exists a polynomial $h_\lambda$ such that
\begin{equation*}
    p_\lambda(t)^2-p_\lambda(t)=(t^d-1)h_\lambda(t).
\end{equation*}
Combining this with $V^d=1$, we have $P_\lambda^2=P_\lambda$.
Thus, every $P_\lambda$ equals $0$ or $1$.
Since their sum equals $1$, exactly one equals $1$.
Consequently, $V=\lambda1$ for some $\lambda\in\Lambda_d$.
\end{proof}

We now show that every root in $\mathcal R$ of a strongly irreducible homogeneous operator is itself homogeneous.
This allows us to determine the existence of roots in terms of the bounded commutant degree semigroup defined in \Cref{def Gamma-b} and shows that they are unique up to multiplication by roots of unity.

\begin{proposition}\label{prop SI-root}
Let $T=u^mf(v)\in\mathcal R$, where $m\ne0$ and $f\ne0$ almost everywhere.
Suppose that $T$ is strongly irreducible relative to $\mathcal R$.
Then, for every integer $d\geqslant2$, the following
statements are equivalent:
\begin{enumerate}
\item[$(1)$] $T$ has a $d$-th root in $\mathcal R$;
\item[$(2)$] $T$ has a homogeneous $d$-th root in $\mathcal R$;
\item[$(3)$] $m/d\in\Gamma_b(T)$.
\end{enumerate}
When these equivalent conditions hold, every $d$-th root of $T$ belongs to $\mathcal C_{m/d}(T)\cap\mathcal R$.
If $S$ is one such root, then all the $d$-th roots of $T$ are precisely $\{\lambda S\colon \lambda^d=1\}$.
\end{proposition}

\begin{proof}
We assume that $m>0$, since the case $m<0$ is analogous.
Let $\mathcal A=\{T\}'\cap\mathcal R$, which is abelian in $\mathcal R$ by \Cref{cor commutant-maximal-abelian}.

We first prove that every $d$-th root of $T$ in
$\mathcal R$ is $m/d$-homogeneous.
Suppose that $T=S^d$ for some $S\in\mathcal R$.

\noindent
\textbf{Case 1.}
$T$ is not invertible.

Since $\mathcal C_{-m}(T)=\mathbb CT^{-1}$ and $T^{-1}\notin\mathcal R$, we have $-m\notin\Gamma_b(T)$.
As $m\in\Gamma_b(T)$, \Cref{lem semigroup} gives $\Gamma_b(T)\subseteq\mathbb Z_{\geqslant0}$.
By \Cref{thm commutant}, we may write
\begin{equation*}
    S=a_0+\sum_{n=1}^{\infty}S_n,
\end{equation*}
where $a_0\in\mathbb C$ and $S_n\in\mathcal C_n(T)\cap\mathcal R$.
Comparing the coefficients of degree zero in $S^d=T$, we have $a_0^d=0$, and hence $a_0=0$.
Let $k\geqslant1$ be the smallest integer such that $S_k\ne0$.
Since $S_k$ is invertible in $L^0(\mathcal R,\tau)$ by \Cref{lem commutant-degree}, we have $S_k^d\ne0$.
The lowest nonzero homogeneous component of $S^d$ equals $S_k^d$ and has degree $kd$.
Consequently, we have $m=kd$ and $T=S_k^d$.
If $S\ne S_k$, let $\ell>k$ be the smallest integer such that $S_\ell\ne0$.
Then the homogeneous component of $S^d$ of degree $(d-1)k+\ell$ equals
\begin{equation*}
    dS_k^{d-1}S_\ell\ne0,
\end{equation*}
which contradicts $T=S^d$.
Thus, $S=S_k$.
In particular, $S$ is $m/d$-homogeneous.

\noindent
\textbf{Case 2.}
$T$ is invertible.

Since $T$ is invertible, $S$ is invertible with inverse $S^{-1}=S^{d-1}T^{-1}$.
For every $z\in\mathbb T$, let $\pi_z$ be the automorphism of $\mathcal R$ such that $\pi_z(u)=zu$ and $\pi_z(v)=v$ as defined in \eqref{equ pi-z}.
Let
\begin{equation*}
    V(z)=z^{-m}\pi_{z^d}(S)S^{-1}\qquad (z\in\mathbb T),\qquad \Lambda_d=\{\lambda\in\mathbb C\colon\lambda^d=1\}.
\end{equation*}
Since $\pi_{z^d}(T)=z^{dm}T$, the operator $\pi_{z^d}(S)$ commutes with $T$.
Thus, $V(z)\in\mathcal A$ and
\begin{equation*}
    V(z)^d=z^{-md}\pi_{z^d}(T)T^{-1}=1.
\end{equation*}
By \Cref{lem finite-order-scalar}, we have $V(z)=\lambda(z)1$ for some scalar $\lambda(z)\in\Lambda_d$.
Since the map $z\mapsto V(z)$ is continuous in $L^2(\mathcal R,\tau)$, $\mathbb T$ is connected, and $V(1)=1$, we obtain that $V(z)=1$ for every $z\in\mathbb T$, that is,
\begin{equation*}
    \pi_{z^d}(S)=z^mS
    \qquad\text{for every}~z\in\mathbb T.
\end{equation*}
In particular, $S=z^mS$ for every $z\in\Lambda_d$.
Thus, $d\mid m$ and $\pi_w(S)=w^{m/d}S$ for every $w\in\mathbb T$.
Hence $S$ is $m/d$-homogeneous.

We have proved $(1)\Rightarrow(2)$ and the homogeneity assertion.
If $S$ is a homogeneous $d$-th root of $T$, then $S$ commutes with $T$ and has degree $m/d$.
Hence $(2)\Rightarrow(3)$.
Suppose that $m/d\in\Gamma_b(T)$.
Choose a nonzero $S\in\mathcal C_{m/d}(T)\cap\mathcal R$.
By \Cref{lem commutant-degree}, $S^d$ is a nonzero element of $\mathcal C_m(T)=\mathbb CT$.
Up to a scalar multiple, we have $T=S^d$.
This proves $(3)\Rightarrow(1)$.
Moreover, since $\mathcal C_{m/d}(T)$ is one dimensional, any two $d$-th roots of $T$ differ by a scalar $d$-th root of unity.
\end{proof}

For continuous weights, the preceding proposition yields an obstruction to the existence of roots in terms of the zeros of the weight.

\begin{corollary}\label{cor root-continuous}
Let $T=u^mf(v)\in\mathcal R$, where $m\ne0$ and $f\ne0$ almost everywhere.
Suppose that $f\in C(\mathbb T)$ and there exists $z_0\in\mathbb T$ such that
\begin{equation*}
	f(z_0)=0,\qquad f(\alpha^{-n}z_0)f(\alpha^nz_0)\ne0\qquad\text{whenever}~n\mid|m|.
\end{equation*}
Then $T$ has no $d$-th roots in $\mathcal R$ for any integer $d\geqslant 2$.
\end{corollary}

\begin{proof}
We assume that $m\geqslant 1$ since the case $m\leqslant -1$ is analogous.
Since $f$ is continuous and $f(z_0)=0$, it is not invertible in $L^\infty(\mathbb T,\mu)$.
Then $T$ is strongly irreducible relative to $\mathcal R$ by \Cref{thm SI-noninvertible}.
By \Cref{prop SI-root}, it suffices to prove that $T$ has no homogeneous $d$-th roots in $\mathcal R$ for any $d\geqslant 2$.
Suppose, towards a contradiction, that $T=S^d$ for some $d\geqslant2$ and $n$-homogeneous operator $S=u^ng(v)\in\mathcal R$.
Then $n\mid m$ and
\begin{equation*}
    f(z)=\prod_{j=0}^{d-1}g(\alpha^{jn}z)
    \qquad\text{almost everywhere}.
\end{equation*}
Consequently, we have
\begin{align*}
    |f(\alpha^{-n}z)f(\alpha^nz)|
    &=
    |f(z)|\,
    |g(\alpha^{-n}z)g(\alpha^{dn}z)|
    \prod_{j=1}^{d-2}|g(\alpha^{jn}z)|\\
    &\leqslant \|g\|_\infty^d|f(z)|
    \qquad\text{almost everywhere}.
\end{align*}
Since both sides of the preceding inequality are continuous, it holds everywhere on $\mathbb T$.
Evaluating at $z_0$ gives $f(\alpha^{-n}z_0)f(\alpha^nz_0)=0$.
That is a contradiction.
\end{proof}

\subsection{\texorpdfstring{$K_0$}{K0}-groups of relative commutants in \texorpdfstring{$C^*$}{C*}-algebras}

Let $\mathfrak A=C^*(u,W^*(v))$.
In \Cref{subsec Red}, we showed that every operator in $\mathfrak A$ is a norm limit of reducible ones relative to $\mathcal R$.
In this subsection, we show that the $K_0$-group of the relative commutant $\{T\}'\cap\mathfrak A$ is isomorphic to $\mathbb Z$ for every injective homogeneous operator $T\in L^0(\mathcal R,\tau)$ of nonzero degree.
Recall that $\Gamma_b(T)$ was introduced in \Cref{def Gamma-b}.
We first consider the case in which $\Gamma_b(T)$ is contained in either $\mathbb Z_{\geqslant0}$ or $\mathbb Z_{\leqslant0}$.

\begin{lemma}\label{lem one-sided-commutant-V}
Let $T=u^mf(v)\in L^0(\mathcal R,\tau)$, where $m\ne0$ and $f\ne0$ almost everywhere, and let $\mathcal A=\{T\}'\cap\mathfrak A$.
If  $\Gamma_b(T)$ is contained in either $\mathbb Z_{\geqslant0}$ or $\mathbb Z_{\leqslant0}$, then
\begin{equation*}
    V(\mathcal A)\cong\mathbb Z_{\geqslant 0},
    \qquad
    K_0(\mathcal A)\cong\mathbb Z,
\end{equation*}
with $[1]$ corresponding to $1$.
\end{lemma}

\begin{proof}
As in \eqref{equ Abel-Poisson}, let $\Phi_r$ be the Abel--Poisson maps 
\begin{equation*}
	\Phi_r(x)=\sum_{n\in\mathbb Z}r^{|n|}\Pi_n(x),\qquad 0\leqslant r<1.
\end{equation*}
Let $\Phi_1(x)=x$.
Since $\Pi_n(\mathcal A)\subseteq\mathcal A$ by \Cref{lem frequently}, it follows that $\Phi_r(\mathcal A)\subseteq\mathcal A$.
Since $\Gamma_b(T)$ is contained in either $\mathbb Z_{\geqslant0}$ or $\mathbb Z_{\leqslant0}$, each $\Phi_r|_{\mathcal A}$ is an algebra homomorphism by the formula \eqref{equ Pi-n-xy} at the end of \Cref{subsec 2-approximations}.
Moreover, $r\mapsto\Phi_r(x)$ is norm-continuous on $[0,1]$ for every $x\in\mathcal A$ by \Cref{prop Abel-Poisson}.
By \Cref{lem commutant-degree}, we have
$\mathcal C_0(T)=\mathbb C1$, and hence
\begin{equation*}
    \Phi_0(x)=\tau(x)1,\qquad x\in\mathcal A.
\end{equation*}
Let $P\in M_N(\mathcal A)$ be an idempotent.
Applying $\Phi_r$ entrywise gives a norm-continuous path of idempotents
\begin{equation*}
    P_r=\Phi_r^{(N)}(P),
    \qquad 0\leqslant r\leqslant1,
\end{equation*}
with $P_1=P$ and $P_0=(\tau(P_{ij}))_{1\leqslant i,j\leqslant N}\in M_N(\mathbb C)$.
Hence $P$ and $P_0$ are similar in $M_N(\mathcal A)$.
The matrix idempotent $P_0$ is similar to $\operatorname{diag}(I_k,0_{N-k})$, where $k=\operatorname{rank}(P_0)$.
Thus, every class in $V(\mathcal A)$ is represented by an identity matrix.
Applying $\Phi_0$ entrywise shows that identity matrices of different ranks cannot represent the same class.
Consequently, we obtain that $V(\mathcal A)\cong\mathbb Z_{\geqslant 0}$, with $[1]$ corresponding to $1$.
Taking the Grothendieck group gives $K_0(\mathcal A)\cong\mathbb Z$.
\end{proof}

We now consider the case in which $\Gamma_b(T)$ contains both positive and negative integers.
By \Cref{lem semigroup}, it is a nonzero subgroup of $\mathbb Z$.
The following lemma identifies the relative commutant as the norm closure of a Laurent polynomial algebra.

\begin{lemma}\label{lem group-commutant-Laurent}
Let $T=u^mf(v)\in L^0(\mathcal R,\tau)$, where $m\ne0$ and $f\ne0$ almost everywhere.
Suppose that $\Gamma_b(T)=k\mathbb Z$ for some integer $k\geqslant1$, and let $\mathcal A=\{T\}'\cap\mathfrak A$.
Then every nonzero $S\in\mathcal C_k(T)\cap\mathcal R$ is invertible in $\mathfrak A$, and $\mathcal A=\overline{\mathbb C[S,S^{-1}]}^{\|\cdot\|}$.
\end{lemma}

\begin{proof}
Choose a nonzero $S\in\mathcal C_k(T)\cap\mathcal R$.
By \Cref{lem commutant-degree}, $S$ is invertible
in $L^0(\mathcal R,\tau)$ and $\mathcal C_{-k}(T)=\mathbb CS^{-1}$.
Since $-k\in\Gamma_b(T)$, it follows that $S^{-1}$ is a bounded homogeneous operator.
Thus, $S$ is invertible in $\mathfrak A$.
By \Cref{lem commutant-degree} again, we have
\begin{equation*}
    \mathcal C_j(T)\cap\mathcal R
    =
    \begin{cases}
        \mathbb CS^n, & j=kn,\quad n\in\mathbb Z,\\
        \{0\}, & j\notin k\mathbb Z.
    \end{cases}
\end{equation*}
It is clear that $\overline{\mathbb C[S,S^{-1}]}^{\|\cdot\|}\subseteq\mathcal A$.
Conversely, let $x\in\mathcal A\subseteq\mathfrak A$.
Then the Fej\'er sums $\sigma_N(x)$ converge to $x$ in operator norm by \Cref{prop Fejer}.
Since $\sigma_N(x)\in\mathbb C[S,S^{-1}]$, it follows that $x\in\overline{\mathbb C[S,S^{-1}]}^{\|\cdot\|}$.
\end{proof}

Next we identify the maximal ideal space of the relative commutant $\{T\}'\cap\mathfrak A$.

\begin{lemma}\label{lem Laurent-commutant-spectrum}
Under the hypotheses of \Cref{lem group-commutant-Laurent}, we have
\begin{equation*}
    \sigma(S)=\left\{\lambda\in\mathbb C\colon
    r(S^{-1})^{-1}\leqslant|\lambda|\leqslant r(S)\right\}.
\end{equation*}
Moreover, the map
\begin{equation*}
    \mathfrak M(\mathcal A)\to\sigma(S),
    \qquad
    \varphi\mapsto\varphi(S),
\end{equation*}
is a homeomorphism from the maximal ideal space $\mathfrak M(\mathcal A)$ onto $\sigma(S)$.
\end{lemma}

\begin{proof}
By \Cref{lem group-commutant-Laurent}, we have $\mathcal A=\overline{\mathbb C[S,S^{-1}]}^{\|\cdot\|}$.
Thus, $\mathcal A$ is a commutative unital Banach algebra.
We first observe that $\sigma_{\mathcal A}(S)=\sigma(S)$, where the latter spectrum is computed in $\mathcal R$.
Indeed, if $\lambda-S$ is invertible in $\mathcal R$, then its inverse belongs to the unital $C^*$-subalgebra $\mathfrak A$ of $\mathcal R$.
Since $\lambda-S$ commutes with $T$, its inverse also commutes with $T$ in $L^0(\mathcal R,\tau)$.
Hence $(\lambda1-S)^{-1}\in\mathcal A$, which proves
the asserted equality of spectra.
As in \eqref{equ pi-z}, let $\pi_z$ be the automorphism of $\mathcal R$ determined by
\begin{equation*}
	\pi_z(u)=zu,\qquad\pi_z(v)=v.
\end{equation*}
Since $\pi_z(S)=z^kS$ for every $z\in\mathbb T$, $\sigma(S)$ is rotation-invariant.
Let $a=r(S^{-1})^{-1}$ and $b=r(S)$.
Then $a$ and $b$ are the minimum and maximum moduli of points in $\sigma(S)$, respectively.
Suppose that some $\rho\in(a,b)$ does not occur as the modulus of a point in $\sigma(S)$.
Then the circle $|\lambda|=\rho$ lies in the resolvent set, and the Riesz idempotent
\begin{equation*}
    E=\frac1{2\pi i}\int_{|\lambda|=\rho}(\lambda-S)^{-1}\,d\lambda
\end{equation*}
belongs to $\mathcal A$.
Since $\sigma(S)$ has points both inside and outside this circle, $E$ is neither $0$ nor $1$.
On the other hand, for every $z\in\mathbb T$, we have
\begin{align*}
    \pi_z(E)
    &=
    \frac1{2\pi i}
    \int_{|\lambda|=\rho}
    (\lambda-z^kS)^{-1}\,d\lambda
    =\frac1{2\pi i}
    \int_{|\lambda|=\rho}
    (z^k\lambda-z^kS)^{-1}\,d(z^k\lambda)
    \\
    &=
    \frac1{2\pi i}
    \int_{|\mu|=\rho}
    (\mu1-S)^{-1}\,d\mu
    =E.
\end{align*}
Thus, $E$ is $0$-homogeneous.
Since $E\in\mathcal A$, \Cref{lem commutant-degree} implies that $E\in\mathcal C_0(T)=\mathbb C1$.
Thus, $E$ is either $0$ or $1$.
That is a contradiction.
Therefore, every modulus in $[a,b]$ occurs in $\sigma(S)$.
Together with rotation invariance, we obtain that
\begin{equation*}
    \sigma(S)=\{\lambda\in\mathbb C\colon
    a\leqslant|\lambda|\leqslant b\}.
\end{equation*}
Finally, since $\mathcal A$ is a commutative unital Banach algebra, we have
\begin{equation*}
    \sigma(S)=\sigma_{\mathcal A}(S)
    =
    \{\varphi(S)\colon\varphi\in\mathfrak M(\mathcal A)\}.
\end{equation*}
Thus, the map $\varphi\mapsto\varphi(S)$ is surjective.
Since $S$ is invertible in $\mathcal A$, every character satisfies $\varphi(S^{-1})=\varphi(S)^{-1}$.
Consequently, any two characters agreeing on $S$ agree on every Laurent polynomial in $S$.
By norm density and continuity, they agree on $\mathcal A$.
Hence the map is also injective.
It is continuous by the definition of the topology on $\mathfrak M(\mathcal A)$.
Since $\mathfrak M(\mathcal A)$ is compact and $\sigma(S)$ is Hausdorff, this map is a homeomorphism.
\end{proof}

Combining the preceding lemmas with the Novodvorskii--Taylor theorem gives the following result.

\begin{proposition}\label{prop norm-commutant-K0}
Let $T=u^mf(v)\in L^0(\mathcal R,\tau)$, where
$m\ne0$ and $f\ne0$ almost everywhere.
Let $\mathfrak A=C^*(u,W^*(v))$ and $\mathcal A=\{T\}'\cap\mathfrak A$.
Then
\begin{equation*}
    V(\mathcal A)\cong\mathbb Z_{\geqslant 0},
    \qquad
    K_0(\mathcal A)\cong\mathbb Z,
\end{equation*}
with $[1]$ corresponding to $1$.
\end{proposition}

\begin{proof}
If $\Gamma_b(T)$ is contained in either $\mathbb Z_{\geqslant0}$ or $\mathbb Z_{\leqslant0}$, then the conclusion follows from \Cref{lem one-sided-commutant-V}.
Otherwise, \Cref{lem semigroup} implies that $\Gamma_b(T)=k\mathbb Z$ for some integer $k\geqslant1$.

By \Cref{lem group-commutant-Laurent} and \Cref{lem Laurent-commutant-spectrum}, $\mathcal A$ is a commutative unital Banach algebra whose maximal ideal space is homeomorphic to a closed annulus, possibly a circle.
Thus, $\mathfrak M(\mathcal A)$ is connected and has the homotopy type of $\mathbb T$.
By the Novodvorskii--Taylor theorem \cite[Section~7.5]{Taylor-76} and homotopy invariance, we have
\begin{equation*}
    V(\mathcal A)\cong
    V\bigl(C(\mathfrak M(\mathcal A))\bigr)
    \cong\mathbb Z_{\geqslant 0},
\end{equation*}
with $[1]$ corresponding to $1$.
Passing to the Grothendieck group yields $K_0(\mathcal A)\cong\mathbb Z$.
\end{proof}

For injective homogeneous operators of degree one in $C^*(u,v)$, their relative commutants in $C^*(u,v)$ and $C^*(u,W^*(v))$ coincide.
Thus, \Cref{prop norm-commutant-K0} yields the following corollary.

\begin{corollary}\label{cor rotation-commutant-K0}
Let $T=uf(v)\in C^*(u,v)$, where $f\in C(\mathbb T)$ is nonzero almost everywhere, and let $\mathcal B=\{T\}'\cap C^*(u,v)$.
Then
\begin{equation*}
    V(\mathcal B)\cong\mathbb Z_{\geqslant 0},
    \qquad
    K_0(\mathcal B)\cong\mathbb Z,
\end{equation*}
with $[1]$ corresponding to $1$.
\end{corollary}

\begin{proof}
Let $\mathfrak A=C^*(u,W^*(v))$.
By \Cref{prop norm-commutant-K0}, it suffices to show that
\begin{equation*}
    \{T\}'\cap\mathfrak A
    =
    \{T\}'\cap C^*(u,v).
\end{equation*}
Indeed, by \Cref{prop Fejer} and \Cref{thm commutant}, we have
\begin{equation*}
	\{T\}'\cap\mathfrak A=
	\begin{cases}
		\overline{\mathbb C[T]}^{\|\cdot\|}, & T~\mbox{is not invertible in}~\mathcal R,\\
		\overline{\mathbb C[T,T^{-1}]}^{\|\cdot\|}, & \mbox{otherwise}.
	\end{cases}
\end{equation*}
It follows that $\{T\}'\cap\mathfrak A\subseteq C^*(u,v)$.
This completes the proof.
\end{proof}

\section{Bilateral weighted shift operators}\label{sec bilateral}

In \Cref{prop norm-commutant-K0}, we computed the $K_0$-group of $\{T\}'\cap C^*(u,W^*(v))$.
To study the $K_0$-group of $\{T\}'\cap\mathcal R$, we develop weighted shift models in this section and the next.
For the sake of convenience, we will focus on $T=uf(v)\in\mathcal R$.

\subsection{Bilateral weighted shift model}\label{subsec bilateral-model}

A \emph{bilateral weighted shift operator} $B$ on $\ell^2(\mathbb Z)$ is determined by
\begin{equation*}
	Be_n=w_ne_{n+1},\qquad n\in\mathbb Z,
\end{equation*}
where $\{e_n\}_{n\in\mathbb Z}$ is the standard orthonormal basis of $\ell^2(\mathbb Z)$ and the weight sequence $\{w_n\}_{n\in\mathbb Z}$ is bounded.
In this subsection, we show that for $T=uf(v)\in\mathcal R$ with $f\ne 0$ almost everywhere, the relative commutant $\{T\}'\cap\mathcal R$ is isomorphic to $\{B\}'$ for some bilateral weighted shift operator $B$.

For $f\in L^\infty(\mathbb T,\mu)$, choose a measurable representative satisfying $|f(z)|\leqslant\|f\|_\infty$ for every $z\in\mathbb T$.
For $z\in\mathbb T$, let $B_{f,z}$ be the bilateral weighted shift operator defined by
\begin{equation}\label{equ B-f-z}
    B_{f,z}e_n=f(\alpha^nz)e_{n+1},
    \qquad n\in\mathbb Z.
\end{equation}
Let $T=uf(v)$ and $\mathcal V=W^*(v)$.
Then
\begin{equation*}
    L^2(\mathcal R,\tau)
    =
    \bigoplus_{n\in\mathbb Z}u^nL^2(\mathcal V,\tau),
    \qquad
    L^2(\mathcal V,\tau)\cong L^2(\mathbb T,\mu),
\end{equation*}
where the second identification sends $g(v)$ to $g$.
Consequently, there exists a unitary operator
\begin{equation*}
	W\colon L^2(\mathcal R,\tau)\to L^2(\mathbb T,\ell^2(\mathbb Z))
\end{equation*}
defined by $\bigl(Wu^ng_n(v)\bigr)(z)=g_n(z)e_n$ for $n\in\mathbb Z$.
Let $L_T$ denote the left multiplication operator by $T$ on $L^2(\mathcal R,\tau)$.
For the vector-valued function $\eta(z)=g_n(z)e_n$, it is straightforward to verify that
\begin{equation*}
    (WL_TW^*\eta)(z)
    =
    f(\alpha^nz)g_n(z)e_{n+1}
    =
    B_{f,z}\eta(z)\qquad\text{almost everywhere}.
\end{equation*}
By linearity and density, this identity holds for every $\eta\in L^2(\mathbb T,\ell^2(\mathbb Z))$.
Therefore, we have
\begin{equation}\label{equ T-direct-integral}
	L^2(\mathbb T,\ell^2(\mathbb Z))
    =
    \int_{\mathbb T}^{\oplus}\ell^2(\mathbb Z)\,d\mu(z),\qquad
    WL_TW^*
    =
    \int_{\mathbb T}^{\oplus}B_{f,z}\,d\mu(z).
\end{equation}
Hence left multiplication by $T=uf(v)$ acts fiberwise as $B_{f,z}$ defined in \eqref{equ B-f-z}.

\begin{lemma}\label{lem fiber-Laurent-norm}
With the notation of \eqref{equ B-f-z} and $T=uf(v)\in\mathcal R$, if $f$ is invertible in $L^\infty(\mathbb T,\mu)$, then there exists a null set $E_f\subseteq\mathbb T$ such that
\begin{equation*}
    \|p(T)\|=\|p(B_{f,z})\|,
    \qquad z\in\mathbb T\setminus E_f,
\end{equation*}
for every Laurent polynomial $p$.
\end{lemma}

\begin{proof}
Since $f$ is invertible in $L^\infty(\mathbb T,\mu)$, we can choose a measurable representative satisfying
\begin{equation*}
	0<\|f^{-1}\|_\infty^{-1}\leqslant|f(z)|\leqslant\|f\|_\infty\qquad\text{for every}~z\in\mathbb T.
\end{equation*}
Thus, $B_{f,z}$ is invertible and \eqref{equ T-direct-integral} implies that $WL_{T^{-1}}W^*=\int_{\mathbb T}^{\oplus}B_{f,z}^{-1}\,d\mu(z)$.
It follows that
\begin{equation*}
    WL_{p(T)}W^*
    =
    \int_{\mathbb T}^{\oplus}p(B_{f,z})\,d\mu(z)
\end{equation*}
for every Laurent polynomial $p$.
In particular, we obtain that
\begin{equation}\label{equ fiber-Laurent-essential-norm}
    \|p(T)\|
    =
    \operatorname*{ess\,sup}_{z\in\mathbb T}
    \|p(B_{f,z})\|.
\end{equation}
Let $U$ be the unweighted bilateral shift operator defined by $Ue_n=e_{n+1}$.
A direct computation gives $B_{f,\alpha z}=U^*B_{f,z}U$.
Hence $\|p(B_{f,z})\|=\|p(B_{f,\alpha z})\|$ for every Laurent polynomial $p$.
By ergodicity of the irrational rotation, the measurable function $z\mapsto\|p(B_{f,z})\|$ is essentially constant.
Together with \eqref{equ fiber-Laurent-essential-norm}, we have
\begin{equation*}
    \|p(T)\|=\|p(B_{f,z})\|\qquad\text{almost everywhere}.
\end{equation*}
Apply the preceding equality to the countable set of Laurent polynomials with coefficients in $\mathbb Q+i\mathbb Q$, and let $E_f$ denote the union of their exceptional null sets.
For an arbitrary Laurent polynomial $p$, approximate its coefficients by elements of $\mathbb Q+i\mathbb Q$.
Passing to the limit then establishes the asserted equality for every Laurent polynomial $p$ and every $z\in\mathbb T\setminus E_f$.
\end{proof}

\begin{remark}\label{rem fiber-Laurent-norm}
By the same proof, if $f$ is not invertible in $L^\infty(\mathbb T,\mu)$, then there exists a null set $E_f\subseteq\mathbb T$ such that
\begin{equation*}
    \|p(uf(v))\|=\|p(B_{f,z})\|,
    \qquad z\in\mathbb T\setminus E_f,
\end{equation*}
for every polynomial $p$.
\end{remark}

Recall that a bounded linear map between weak-operator closed operator spaces is called \emph{normal} if it is ultraweakly continuous; equivalently, it is weak-operator continuous on bounded sets.

\begin{lemma}\label{lem Laurent-comparison-commutant}
Let $T=uf(v)\in\mathcal R$, where $f$ is invertible in $L^\infty(\mathbb T,\mu)$, and let $B$ be an invertible bilateral weighted shift operator.
Suppose that there exists $C\geqslant1$ such that
\begin{equation}\label{equ Laurent-comparison-commutant}
    \|p(B)\|\leqslant\|p(T)\|
    \leqslant C\|p(B)\|
\end{equation}
for every Laurent polynomial $p$.
Then there exists a unique unital normal algebra isomorphism
\begin{equation*}
    \Phi\colon\{T\}'\cap\mathcal R\to\{B\}',
    \qquad \Phi(T)=B.
\end{equation*}
Moreover, $\|\Phi\|=1$, $\|\Phi^{-1}\|\leqslant C$, and $\Phi^{-1}$ is also normal.
\end{lemma}

\begin{proof}
By \Cref{prop Fejer,thm commutant}, every $X\in\{T\}'\cap\mathcal R$ has a unique $L^2$-Fourier expansion $X=\sum_{n\in\mathbb Z}a_nT^n$.
Its Fej\'er sums
\begin{equation*}
    \sigma_N(X)
    =
    \sum_{|n|\leqslant N}
    \left(1-\frac{|n|}{N+1}\right)a_nT^n
\end{equation*}
converge to $X$ in the weak-operator topology and satisfy $\|\sigma_N(X)\|\leqslant\|X\|$.
Similarly, by Theorems 3 and 12 of \cite[pp.~62,~90]{Shields-74}, every $Y\in\{B\}'$ has a unique formal Fourier expansion $Y=\sum_{n\in\mathbb Z}b_nB^n$.
Its Fej\'er sums
\begin{equation*}
    \sigma_N(Y)
    =
    \sum_{|n|\leqslant N}
    \left(1-\frac{|n|}{N+1}\right)b_nB^n
\end{equation*}
converge to $Y$ in the weak-operator topology and satisfy $\|\sigma_N(Y)\|\leqslant\|Y\|$.
By \eqref{equ Laurent-comparison-commutant}, for $X=\sum_{n\in\mathbb Z}a_nT^n\in\{T\}'\cap\mathcal R$, the formal series $Y=\sum_{n\in\mathbb Z}a_nB^n$ satisfies
\begin{equation*}
	\|\sigma_N(Y)\|\leqslant\|\sigma_N(X)\|\leqslant\|X\|.
\end{equation*}
Clearly, every weak-operator cluster point of $\{\sigma_N(Y)\}_{N=1}^\infty$ has Fourier expansion $Y=\sum_{n\in\mathbb Z}a_nB^n$.
By weak-operator compactness of bounded closed balls in $\{B\}'$ and uniqueness of Fourier expansion, the sequence converges in the weak-operator topology.
Thus, $Y$ represents an operator in $\{B\}'$ with $\|Y\|\leqslant\|X\|$.
Define $\Phi(X)=Y$.
Conversely, for $Y=\sum_{n\in\mathbb Z}b_nB^n\in\{B\}'$, the formal series $X=\sum_{n\in\mathbb Z}b_nT^n$ satisfies
\begin{equation*}
	\|\sigma_N(X)\|\leqslant C\|\sigma_N(Y)\|\leqslant C\|Y\|.
\end{equation*}
It follows that
\begin{equation*}
	\sum_{|n|\leqslant N}
	\left(1-\frac{|n|}{N+1}\right)^2
	|b_n|^2\|T^n\|_2^2
	=\|\sigma_N(X)\|_2^2
	\leqslant\|\sigma_N(X)\|^2
	\leqslant C^2\|Y\|^2.
\end{equation*}
By Fatou's lemma, $\sum_{n\in\mathbb Z}|b_n|^2\|T^n\|_2^2\leqslant C^2\|Y\|^2$.
Hence $X\in L^2(\mathcal R,\tau)$ and its Fej\'er sums converge to $X$ in $L^2$.
Since the Fej\'er sums are uniformly bounded, it follows that $X\in\mathcal R$.
Since the Fej\'er sums converge ultraweakly to $X$, we have $X\in\{T\}'\cap\mathcal R$ and $\|X\|\leqslant C\|Y\|$.
Define $\Psi(Y)=X$.
Since both maps preserve the Fourier coefficients, they are linear and $\Psi=\Phi^{-1}$.
Moreover, by $\Phi(1)=I$, we have $\|\Phi\|=1$ and $\|\Psi\|\leqslant C$.

If a bounded net $X_\lambda$ converges to $X$ in the weak-operator topology, then the Fourier coefficients of $X_\lambda$ converge to those of $X$.
Hence every weak-operator cluster point of the bounded net $\Phi(X_\lambda)$ has the same Fourier expansion as $\Phi(X)$.
It follows that $\Phi(X_\lambda)\to\Phi(X)$ in the weak-operator topology.
Thus, $\Phi$ is normal.
The same argument proves that $\Psi$ is normal.

For any $X_1,X_2\in\{T\}'\cap\mathcal R$, their Fej\'er sums are Laurent polynomials of $T$ and we have
\begin{equation*}
	\Phi\bigl(\sigma_{N_1}(X_1)\sigma_{N_2}(X_2)\bigr)
	=\Phi(\sigma_{N_1}(X_1))\Phi(\sigma_{N_2}(X_2)).
\end{equation*}
Then $\Phi(X_1X_2)=\Phi(X_1)\Phi(X_2)$ by letting $N_1\to\infty$ and then $N_2\to\infty$.
Clearly, $\Phi(1)=I$ and $\Phi(T)=B$.
Thus, $\Phi$ is an algebra isomorphism.

Finally, any unital algebra isomorphism sending $T$ to $B$ sends $T^{-1}$ to $B^{-1}$, and hence agrees with $\Phi$ on Laurent polynomials in $T$.
If it is normal, Fej\'er approximation implies that it agrees with $\Phi$ on $\{T\}'\cap\mathcal R$.
This proves uniqueness.
\end{proof}

\begin{remark}\label{rem polynomial-comparison-commutant}
Let $T=uf(v)$, where $f\in L^\infty(\mathbb T,\mu)$ is nonzero almost everywhere and not invertible.
Let $B$ be an injective noninvertible bilateral weighted shift operator.
Suppose that there exists $C\geqslant1$ such that
\begin{equation*}
    \|p(B)\|\leqslant\|p(T)\|
    \leqslant C\|p(B)\|
\end{equation*}
for every polynomial $p$.
Then there exists a unique unital normal algebra isomorphism
\begin{equation*}
    \Phi\colon\{T\}'\cap\mathcal R\to\{B\}',
    \qquad \Phi(T)=B.
\end{equation*}
Moreover, $\|\Phi\|=1$, $\|\Phi^{-1}\|\leqslant C$, and $\Phi^{-1}$ is also normal.
\end{remark}

The preceding norm identities and comparison results yield the desired bilateral weighted shift model.

\begin{theorem}\label{thm bilateral-model}
Let $T=uf(v)\in\mathcal R$, where $f\in L^\infty(\mathbb T,\mu)$ and $f\ne 0$ almost everywhere.
Then there exists a null set $E_f\subseteq\mathbb T$ such that, for every $z\in\mathbb T\setminus E_f$, the bilateral weighted shift operator
\begin{equation*}
	B_{f,z}e_n=f(\alpha^nz)e_{n+1},\qquad n\in\mathbb Z
\end{equation*}
is bounded and injective, and there exists a unique unital normal isometric algebra isomorphism
\begin{equation*}
	\Phi\colon\{T\}'\cap\mathcal R\to\{B_{f,z}\}',\qquad \Phi(T)=B_{f,z}.
\end{equation*}
Moreover, $\Phi^{-1}$ is also normal.
\end{theorem}

\begin{proof}
If $f$ is invertible in $L^\infty(\mathbb T,\mu)$, then
\begin{equation*}
	0<\|f^{-1}\|_\infty^{-1}\leqslant|f(z)|\leqslant\|f\|_\infty\qquad\text{almost everywhere}.
\end{equation*}
Thus, $B_{f,z}$ is invertible almost everywhere, and the conclusion follows from \Cref{lem fiber-Laurent-norm,lem Laurent-comparison-commutant} with $C=1$.
If $f$ is not invertible, then
\begin{equation*}
	0<|f(z)|\leqslant\|f\|_\infty,
	\qquad
	\inf_{n\in\mathbb Z}|f(\alpha^nz)|=0\qquad\text{almost everywhere}.
\end{equation*}
Thus, $B_{f,z}$ is bounded, injective, and noninvertible almost everywhere.
The conclusion follows from \Cref{rem fiber-Laurent-norm,rem polynomial-comparison-commutant} with $C=1$.
\end{proof}

\begin{remark}\label{rem bilateral-model}
Let $T=u^mf(v)\in L^0(\mathcal R,\tau)$, where $m\ne 0$ and $f\ne0$ almost everywhere.
Then there exist a densely defined closed injective bilateral weighted shift operator $B$ on $\mathcal H=\ell^2(\mathbb Z)$ and a unital normal isometric algebra isomorphism
\begin{equation*}
	\Phi\colon\{T\}'\cap\mathcal R
	\to\{B\}'\cap\mathcal B(\mathcal H).
\end{equation*}
Moreover, $\Phi^{-1}$ is normal.
Here the bounded commutant of $B$ is defined by
\begin{equation*}
    \{B\}'\cap\mathcal B(\mathcal H)
    =
    \left\{
    X\in\mathcal B(\mathcal H)\colon
    \begin{array}{l}
        X\mathcal D(B)\subseteq\mathcal D(B),\\
        BX\xi=XB\xi
        ~\text{for every}~\xi\in\mathcal D(B)
    \end{array}
    \right\},
\end{equation*}
where $\mathcal D(B)$ is the domain of $B$.
If $\Gamma_b(T)$ is contained in either $\mathbb Z_{\geqslant0}$ or $\mathbb Z_{\leqslant0}$, then $B$ may instead be chosen to be a densely defined closed injective unilateral weighted shift operator on $\ell^2(\mathbb N)$ as in \Cref{thm unilateral-model}.
These extensions follow from the same Fej\'er approximation argument for a homogeneously power-free root given by \Cref{prop power-free-factorization}.
We omit the details of these extensions.
\end{remark}

For continuous weights, the bilateral weighted shift model is valid along every orbit on which the weights are nonzero.

\begin{corollary}\label{cor bilateral-model-continuous}
Let $T=uf(v)\in\mathcal R$, where $f\in C(\mathbb T)$ and $f\ne 0$ almost everywhere.
Then for every $z\in\mathbb T$ with $f(\alpha^nz)\ne 0$ for all $n\in\mathbb Z$, there exists a unique unital normal isometric algebra isomorphism
\begin{equation*}
	\Phi\colon\{T\}'\cap\mathcal R\to\{B_{f,z}\}',\qquad \Phi(T)=B_{f,z}.
\end{equation*}
Moreover, $\Phi^{-1}$ is normal.
\end{corollary}

\begin{proof}
Since $f$ is uniformly continuous and
\begin{equation*}
    \|B_{f,z}-B_{f,w}\|
    =
    \sup_{n\in\mathbb Z}
    |f(\alpha^nz)-f(\alpha^nw)|,
\end{equation*}
the map $z\mapsto B_{f,z}$ is operator-norm continuous.
Consequently, $z\mapsto\|p(B_{f,z})\|$ is continuous for every polynomial $p$.
If $f$ is invertible, then the same holds for every Laurent polynomial $p$.
By \Cref{lem fiber-Laurent-norm,rem fiber-Laurent-norm}, we obtain
\begin{equation*}
    \|p(T)\|=\|p(B_{f,z})\|,
    \qquad z\in\mathbb T,
\end{equation*}
for every polynomial $p$, and for every Laurent polynomial $p$ when $f$ is invertible.

Let $z\in\mathbb T$ with $f(\alpha^nz)\ne0$ for every $n\in\mathbb Z$.
Then $B_{f,z}$ is injective.
By continuity of $f$, we have
\begin{equation*}
    \inf_{n\in\mathbb Z}|f(\alpha^nz)|
    =
    \min_{w\in\mathbb T}|f(w)|.
\end{equation*}
Thus, $B_{f,z}$ is invertible if and only if $f$ is invertible in $L^\infty(\mathbb T,\mu)$.
We complete the proof by applying \Cref{lem Laurent-comparison-commutant,rem polynomial-comparison-commutant} with $C=1$.
\end{proof}

\subsection{Bilateral weighted shift realization}

We now consider the converse realization problem of \Cref{thm bilateral-model}: whether the commutant of every injective bilateral weighted shift operator can be realized as the relative commutant of a homogeneous operator in $\mathcal R$ of degree one.
We begin with an elementary lemma.

\begin{lemma}\label{lem finite-shift-blocks}
Let $B$ be an injective bilateral weighted shift operator.
For every integer $m\geqslant 1$, let $B_m$ be its compression to the $m$-dimensional subspace
\begin{equation*}
    H_m=\operatorname{span}
    \{e_{-\lfloor m/2\rfloor},\ldots,
    e_{m-1-\lfloor m/2\rfloor}\}.
\end{equation*}
Then for every polynomial $p$, we have
\begin{equation*}
    \|p(B_m)\|\leqslant\|p(B)\|,
    \qquad
    \lim_{m\to\infty}\|p(B_m)\|=\|p(B)\|.
\end{equation*}
Moreover, for any contractions $C$ and $D$ with $C^m=D^m=0$, we have
\begin{equation*}
	\|p(C)-p(D)\|
    \leqslant
    \kappa_m\|C-D\|\|p(B)\|,
    \qquad
    \kappa_m=\sum_{j=1}^{m-1}\frac{j}{\|B^j\|}.
\end{equation*}
\end{lemma}

\begin{proof}
Let $P_m$ be the projection from $\ell^2(\mathbb Z)$ onto $H_m$.
It is straightforward to verify that
\begin{equation*}
    B_m^j=P_mB^j|_{H_m},
    \qquad j\geqslant0.
\end{equation*}
Thus, $p(B_m)=P_mp(B)|_{H_m}$ and $\|p(B_m)\|\leqslant\|p(B)\|$ for every polynomial $p$.
Since $P_m\to I$ in the strong-operator topology and $\|p(B_m)\|=\|P_mp(B)P_m\|$, we have
\begin{equation*}
    \|p(B)\|
    \leqslant\liminf_{m\to\infty}\|p(B_m)\|
    \leqslant\limsup_{m\to\infty}\|p(B_m)\|
    \leqslant\|p(B)\|.
\end{equation*}
This proves the first assertion.

For $z\in\mathbb T$, let $U_z$ be a unitary operator such that $U_ze_n=z^ne_n$.
It is clear that $U_zBU_z^*=zB$.
Then for every polynomial $p(t)=\sum_{j=0}^Na_jt^j$, we have
\begin{equation*}
    a_jB^j
    =
    \int_{\mathbb T}
    z^{-j}U_zp(B)U_z^*\,d\mu(z).
\end{equation*}
It follows that $|a_j|\|B^j\|\leqslant\|p(B)\|$.
Since $\|C\|,\|D\|\leqslant 1$ and
\begin{equation*}
    C^j-D^j
    =
    \sum_{k=0}^{j-1}
    C^{j-1-k}(C-D)D^k,
\end{equation*}
we have $\|C^j-D^j\|\leqslant j\|C-D\|$ for each $j\geqslant 0$.
Thus, we obtain that
\begin{equation*}
	\|p(C)-p(D)\|\leqslant\sum_{j=1}^{m-1}|a_j|\|C^j-D^j\|
	\leqslant\sum_{j=1}^{m-1}j|a_j|\|C-D\|
	\leqslant\kappa_m\|C-D\|\|p(B)\|.
\end{equation*}
Note that both sides are $0$ when $m=1$.
This proves the second assertion.
\end{proof}

In the following lemma, we partition $\mathbb T$ according to the first return times of the irrational rotation to a measurable set.

\begin{lemma}\label{lem return-time-partition}
Let $D$ be a measurable subset of $\mathbb T$ with $\mu(D)>0$ and
\begin{equation*}
    D_m=\{z\in D\colon
    \alpha^jz\notin D\ (1\leqslant j<m),
    \ \alpha^mz\in D\},\qquad m\geqslant 1.
\end{equation*}
Then the sets $\{\alpha^jD_m\colon 0\leqslant j<m,~m\geqslant 1\}$ are pairwise disjoint and
\begin{equation*}
	\mathbb T=\bigsqcup_{m=1}^\infty\bigsqcup_{j=0}^{m-1}\alpha^jD_m\qquad\text{modulo a null set}.
\end{equation*}
\end{lemma}

\begin{proof}
Suppose that $\alpha^ix=\alpha^jy$, where $x\in D_m$, $y\in D_n$, $0\leqslant i<m$, and $0\leqslant j<n$.
If $i>j$, then
\begin{equation*}
    y=\alpha^{i-j}x\in D\qquad\text{and}
    \qquad 1\leqslant i-j<m,
\end{equation*}
contrary to $x\in D_m$.
Similarly, we cannot have $i<j$.
Thus $i=j$, $x=y$, and $m=n$.
This proves pairwise disjointness.

For every $N\geqslant0$, let
\begin{equation*}
    A_N=\bigcup_{n\geqslant N}\alpha^{-n}D,\qquad B_N=\bigcup_{n\geqslant N}\alpha^nD.
\end{equation*}
Since $A_N\subseteq\alpha A_N$ and rotation preserves measure, $A_N=\alpha A_N$ modulo a null set.
Combining this with $\mu(A_N)\geqslant\mu(D)>0$, it follows that $\mu(A_N)=1$.
Similarly, we have $\mu(B_N)=1$.
Thus, for almost every $z\in\mathbb T$, the sets
\begin{equation*}
	\{n\geqslant0\colon\alpha^nz\in D\}\qquad\text{and}\qquad\{n\geqslant0\colon\alpha^{-n}z\in D\}
\end{equation*}
are infinite.
For such $z$, let
\begin{equation*}
    j=\min\{k\geqslant0\colon\alpha^{-k}z\in D\},
    \qquad x=\alpha^{-j}z,\qquad
    m=\min\{k\geqslant1:\alpha^kx\in D\}.
\end{equation*}
Then $x\in D_m$, $0\leqslant j<m$, and $z=\alpha^jx\in\alpha^jD_m$.
This completes the proof.
\end{proof}

The following lemma allows us to remove the zeros of a weight function while keeping the polynomial norms of the corresponding homogeneous operator close to their original values.

\begin{lemma}\label{lem filling-zero-weights}
Let $B$ be an injective bilateral weighted shift operator.
Suppose that $g\in L^\infty(\mathbb T,\mu)$ satisfies that $|g|\leqslant1$ and its zero set $Z(g)$ has positive measure.
Then for every $\varepsilon>0$, there exists
$f\in L^\infty(\mathbb T,\mu)$ such that
\begin{enumerate}
\item[$(1)$] $f=g$ on $\mathbb T\setminus Z(g)$;
\item[$(2)$] $0<|f|\leqslant 1$ almost everywhere and $f$ is not invertible in $L^\infty(\mathbb T,\mu)$;
\item[$(3)$] $\|p(uf(v))-p(ug(v))\|\leqslant\varepsilon\|p(B)\|$ for every polynomial $p$.
\end{enumerate}
\end{lemma}

\begin{proof}
Let $\{Z_k\}_{k=0}^{\infty}$ be a sequence of measurable subsets of $\mathbb T$ such that
\begin{equation*}
	Z(g)=Z_0\supseteq Z_1\supseteq Z_2\supseteq\cdots,
	\qquad
	0<\mu(Z_k)\to 0.
\end{equation*}
For $k,m\geqslant 1$, let $\varepsilon_{k,m}=\min\left\{2^{-k},\varepsilon\min_{1\leqslant j\leqslant m}2^{-j}\|B^j\|\right\}$ and
\begin{equation*}
	Z_{k,m}=\{z\in Z_k\colon\alpha^jz\notin Z_k\,(1\leqslant j<m),\,\alpha^mz\in Z_k\}.
\end{equation*}
For each $k\geqslant1$, it follows from \Cref{lem return-time-partition} that
\begin{equation}\label{equ zero-set-partition}
    \mathbb T
    =
    \bigsqcup_{m=1}^\infty\bigsqcup_{j=0}^{m-1}
    \alpha^jZ_{k,m}
    \qquad\text{modulo a null set}.
\end{equation}
Next we define $f$ almost everywhere on $\mathbb T$ by
\begin{equation*}
    f(z)=
    \begin{cases}
      g(z),
        &z\in\mathbb T\setminus Z_0,\\
      \varepsilon_{k,m},
        &z\in(Z_{k-1}\setminus Z_k)
          \cap\alpha^jZ_{k,m},
          \qquad 1\leqslant j<m,~m\geqslant1,~k\geqslant1.
    \end{cases}
\end{equation*}
Clearly, $f=g$ outside $Z_0$ and $0<|f|\leqslant 1$ almost everywhere.
For every $k\geqslant1$, the definition gives $0<f\leqslant2^{-k-1}$ almost everywhere on $Z_k$.
Since $\mu(Z_k)>0$, we have $\operatorname*{ess\,inf}|f|=0$.
Thus $f$ is not invertible in $L^\infty(\mathbb T,\mu)$.
This proves $(1)$ and $(2)$.

Let $T=uf(v)$ and $S=ug(v)$.
We claim that
\begin{equation}\label{equ filling-power-bound}
    \|T^n-S^n\|
    \leqslant\varepsilon2^{-n}\|B^n\|,
    \qquad n\geqslant1.
\end{equation}
Fix $n\geqslant1$.
Let
\begin{equation*}
	\gamma(z)=\min\{k\geqslant 0\colon
      \alpha^jz\notin Z_k~
      \text{for every}~0\leqslant j<n\}.
\end{equation*}
Then $\gamma<\infty$ almost everywhere since $\bigcap_{k=1}^\infty Z_k$ is a null set.
If $\gamma(z)=0$, then $f(\alpha^jz)=g(\alpha^jz)$
for every $0\leqslant j<n$, and hence
\begin{equation*}
    \prod_{j=0}^{n-1}f(\alpha^jz)
    =
    \prod_{j=0}^{n-1}g(\alpha^jz).
\end{equation*}
Let $k=\gamma(z)\geqslant 1$.
Then $\alpha^jz\notin Z_k$ for $0\leqslant j<n$ and $\alpha^{j_0}z\in Z_{k-1}\setminus Z_k$ for some $0\leqslant j_0<n$.
By \eqref{equ zero-set-partition}, we can write
\begin{equation*}
    z=\alpha^iw,
    \qquad w\in Z_{k,m},\qquad 1\leqslant i<m.
\end{equation*}
Then $i<m-n+1$.
Otherwise, $1\leqslant m-i<n$ and
\begin{equation*}
    \alpha^{m-i}z=\alpha^mw\in Z_k.
\end{equation*}
That contradicts the definition of $\gamma(z)$.
It follows that
\begin{equation*}
    \alpha^jz\in\alpha^{i+j}Z_{k,m},
    \qquad 1\leqslant i+j<m
    \qquad(0\leqslant j<n).
\end{equation*}
In particular, $f(\alpha^{j_0}z)=\varepsilon_{k,m}$.
Since $n<m$ and $|f|\leqslant1$, we have
\begin{equation*}
    \left|\prod_{j=0}^{n-1}f(\alpha^jz)\right|
    \leqslant\varepsilon_{k,m}
    \leqslant\varepsilon2^{-n}\|B^n\|.
\end{equation*}
Since $Z_{k-1}\subseteq Z_0=Z(g)$, it follows that $g(\alpha^{j_0}z)=0$.
Thus, in both cases, we have
\begin{equation*}
    \left|
      \prod_{j=0}^{n-1}f(\alpha^jz)
      -
      \prod_{j=0}^{n-1}g(\alpha^jz)
    \right|
    \leqslant\varepsilon2^{-n}\|B^n\|\qquad\text{almost everywhere}.
\end{equation*}
It is clear that
\begin{equation*}
    T^n-S^n
    =
    u^n\left(
    \prod_{j=0}^{n-1}f(\alpha^jv)
    -
    \prod_{j=0}^{n-1}g(\alpha^jv)
    \right).
\end{equation*}
This proves \eqref{equ filling-power-bound}.

Let $p(t)=\sum_{j=0}^N a_jt^j$ be a polynomial.
By the proof of \Cref{lem finite-shift-blocks}, we have $|a_j|\|B^j\|\leqslant\|p(B)\|$ for every $0\leqslant j\leqslant N$.
Hence
\begin{equation*}
	\|p(T)-p(S)\|
    \leqslant\sum_{j=1}^N|a_j|\|T^j-S^j\|
    \leqslant\varepsilon\sum_{j=1}^N 2^{-j}|a_j|\|B^j\|
    \leqslant\varepsilon\|p(B)\|.
\end{equation*}
This proves $(3)$.
\end{proof}

We now use the preceding three lemmas to construct a homogeneous operator whose polynomial norms are comparable to those of a prescribed bilateral weighted shift.

\begin{lemma}\label{lem realization-B}
Let $B$ be an injective bilateral weighted shift operator.
Then for every $\varepsilon>0$, there exists $f\in L^\infty(\mathbb T,\mu)$, nonzero almost everywhere and not invertible in $L^\infty(\mathbb T,\mu)$, such that
\begin{equation*}
    \|p(B)\|\leqslant\|p(uf(v))\|
    \leqslant(1+\varepsilon)\|p(B)\|
\end{equation*}
for every polynomial $p$.
\end{lemma}

\begin{proof}
Without loss of generality, we assume that $\|B\|=1$.
Write $Be_n=w_ne_{n+1}$ for $n\in\mathbb Z$.
Choose a nowhere dense closed set $D\subseteq\mathbb T$ with $\mu(D)>0$, and let $D_m$ be as in
\Cref{lem return-time-partition}.
Using the partition in that lemma, define $g$ almost
everywhere by
\begin{equation*}
    g(z)=
    \begin{cases}
      0,&z\in D,\\[2pt]
      w_{j-1-\lfloor m/2\rfloor},
        &z\in\alpha^jD_m,\qquad 1\leqslant j<m,~m\geqslant 1.
    \end{cases}
\end{equation*}
Then $|g|\leqslant1$ and $Z(g)=D$ modulo a null set.
As in \eqref{equ B-f-z}, let $B_{g,z}$ be the bilateral weighted shift operator defined by
\begin{equation*}
    B_{g,z}e_n=g(\alpha^nz)e_{n+1},
    \qquad n\in\mathbb Z.
\end{equation*}
For almost every $z\in\mathbb T$, the set $\{n\in\mathbb Z\colon\alpha^nz\in D\}$ is unbounded in both directions, as shown in the proof of \Cref{lem return-time-partition}.
Enumerate this set as
\begin{equation*}
    \cdots<n_{-2}<n_{-1}<n_0<n_1<n_2<\cdots.
\end{equation*}
For each $r\in\mathbb Z$, let
\begin{equation*}
    m_r=n_{r+1}-n_r,
    \qquad
    \mathcal H_r
    =\operatorname{span}\{e_j:n_r<j\leqslant n_{r+1}\}.
\end{equation*}
Then $\alpha^{n_r}z\in D_{m_r}$ and $g(\alpha^{n_r}z)=g(\alpha^{n_{r+1}}z)=0$.
These zero weights imply that $\mathcal H_r$ reduces
$B_{g,z}$.
Moreover, we have
\begin{equation*}
    g(\alpha^{n_r+j}z)
    =w_{j-1-\lfloor m_r/2\rfloor},
    \qquad 1\leqslant j<m_r.
\end{equation*}
Thus, $B_{g,z}|_{\mathcal H_r}$ is unitarily equivalent to $B_{m_r}$ by \Cref{lem finite-shift-blocks}, and hence
\begin{equation*}
    \ell^2(\mathbb Z)=\bigoplus_{r\in\mathbb Z}\mathcal H_r,
    \qquad
    B_{g,z}\cong\bigoplus_{r\in\mathbb Z}B_{m_r}.
\end{equation*}
By \Cref{rem fiber-Laurent-norm,lem finite-shift-blocks}, we obtain
\begin{equation*}
    \|p(ug(v))\|
    =\|p(B_{g,z})\|
    =\sup_{r\in\mathbb Z}\|p(B_{m_r})\|
    \leqslant\|p(B)\|
\end{equation*}
for every polynomial $p$.
Let $f\in L^\infty(\mathbb T,\mu)$ be given by \Cref{lem filling-zero-weights}.
Then
\begin{equation*}
    \|p(uf(v))\|
    \leqslant
    \|p(ug(v))\|+\varepsilon\|p(B)\|
    \leqslant(1+\varepsilon)\|p(B)\|.
\end{equation*}

To prove the lower bound, observe that
$\{m\geqslant1\colon\mu(D_m)>0\}$ is unbounded.
Otherwise, \Cref{lem return-time-partition} would imply
\begin{equation*}
    \mu\left(\bigcup_{j=0}^{N-1}\alpha^jD\right)=1
\end{equation*}
for some $N\geqslant 1$.
This is impossible because this finite union is nowhere dense and closed, and hence its complement is a nonempty open set of positive measure.

Fix $m\geqslant 1$ with $\mu(D_m)>0$, and let $P_m$ be the
projection onto $\operatorname{span}\{e_1,\ldots,e_m\}$.
For almost every $z\in D_m$, we have
\begin{equation*}
    f(\alpha^jz)=g(\alpha^jz)
    =w_{j-1-\lfloor m/2\rfloor},
    \qquad 1\leqslant j<m.
\end{equation*}
Consequently, $P_mB_{f,z}|_{P_m\ell^2(\mathbb Z)}\cong B_m$.
As in the proof of \Cref{lem finite-shift-blocks}, we have
\begin{equation*}
    P_mp(B_{f,z})|_{P_m\ell^2(\mathbb Z)}
    \cong p(B_m)
\end{equation*}
for every polynomial $p$.
By \Cref{rem fiber-Laurent-norm} and $\mu(D_m)>0$, we may choose $z\in D_m$ such that
\begin{equation*}
    \|p(uf(v))\|
    =\|p(B_{f,z})\|
    \geqslant\|p(B_m)\|.
\end{equation*}
Letting $m\to\infty$ along the set $\{m\geqslant1\colon\mu(D_m)>0\}$ and applying
\Cref{lem finite-shift-blocks}, we obtain that $\|p(uf(v))\|\geqslant\|p(B)\|$.
This completes the proof.
\end{proof}

We now prove the main result of this subsection: the commutant of every injective noninvertible bilateral weighted shift admits a nearly isometric realization as the relative commutant of a homogeneous operator in $\mathcal R$.

\begin{proposition}\label{prop noninvertible-bilateral-realization}
Let $B$ be an injective noninvertible bilateral weighted shift operator.
Then for every $\varepsilon>0$, there exist $T=uf(v)\in\mathcal R$, where $f$ is nonzero almost everywhere and noninvertible, and a unital normal algebra isomorphism
\begin{equation*}
    \Phi\colon\{T\}'\cap\mathcal R\to\{B\}',
    \qquad \Phi(T)=B,
\end{equation*}
such that
\begin{equation*}
    \|\Phi\|=1,
    \qquad
    \|\Phi^{-1}\|\leqslant1+\varepsilon.
\end{equation*}
Moreover, $\Phi^{-1}$ is normal.
\end{proposition}

\begin{proof}
Let $f$ be given by \Cref{lem realization-B} and $T=uf(v)$.
Then $f\ne 0$ almost everywhere, $f$ is not invertible in $L^\infty(\mathbb T,\mu)$, and
\begin{equation*}
    \|p(B)\|\leqslant\|p(T)\|
    \leqslant(1+\varepsilon)\|p(B)\|
\end{equation*}
for every polynomial $p$.
Since $B$ is injective and noninvertible, \Cref{rem polynomial-comparison-commutant} gives a unital normal
algebra isomorphism
\begin{equation*}
    \Phi\colon\{T\}'\cap\mathcal R\to\{B\}',
    \qquad \Phi(T)=B,
\end{equation*}
such that $\|\Phi\|=1$, $\|\Phi^{-1}\|\leqslant1+\varepsilon$, and $\Phi^{-1}$ is normal.
This completes the proof.
\end{proof}

Under the assumption $\|B\|=r(B)$, the realization in \Cref{prop noninvertible-bilateral-realization} can be chosen to be isometric.

\begin{corollary}\label{cor isometric-bilateral-realization}
Let $B$ be an injective noninvertible bilateral weighted shift operator with $r(B)=\|B\|$.
Then there exist $T=uf(v)\in\mathcal R$, where $f$ is nonzero almost everywhere and noninvertible, and a unital normal isometric algebra isomorphism
\begin{equation*}
    \Phi\colon\{T\}'\cap\mathcal R\to\{B\}',
    \qquad \Phi(T)=B.
\end{equation*}
Moreover, $\Phi^{-1}$ is normal.
\end{corollary}

\begin{proof}
By scaling, we assume that $r(B)=\|B\|=1$.
Choose a nowhere dense closed set $D\subseteq\mathbb T$ of positive measure, and partition $D$ into measurable sets $\{D_k\}_{k\geqslant1}$, each of positive measure.
Define
\begin{equation*}
    f(z)=
    \begin{cases}
        1,&z\in\mathbb T\setminus D,\\
        2^{-k},&z\in D_k.
    \end{cases}
\end{equation*}
Then $f$ is nonzero everywhere and noninvertible in $L^\infty(\mathbb T,\mu)$.
Let $T=uf(v)$.

For every $n\geqslant1$, we have $\prod_{j=0}^{n-1}f(\alpha^jz)=1$ on the set
\begin{equation*}
    \mathbb T\setminus
    \bigcup_{j=0}^{n-1}\alpha^{-j}D,
\end{equation*}
which is nonempty and open, and hence has positive measure.
Since $|f|\leqslant1$, we obtain that $r(T)=\|T\|=1$.
It follows that $\sigma(B)=\overline{\mathbb D}=\sigma(T)$.
Thus, the spectral mapping theorem and von Neumann's inequality \cite{vN-51} give
\begin{equation*}
    \|p(B)\|
    =\max_{|z|\leqslant1}|p(z)|
    =\|p(T)\|
\end{equation*}
for every polynomial $p$.
The conclusion follows from
\Cref{rem polynomial-comparison-commutant}.
\end{proof}

For invertible shift operators, the comparison must also include negative powers.
Under corresponding assumptions on the norms and spectral radii of $B$ and $B^{-1}$, we obtain the following realization.

\begin{corollary}\label{cor invertible-bilateral-realization}
Let $B$ be an invertible bilateral weighted shift operator satisfying
\begin{equation*}
    \|B\|=r(B),
    \qquad
    \|B^{-1}\|=r(B^{-1}).
\end{equation*}
Then there exist an invertible $T=uf(v)\in\mathcal R$ and a unital normal algebra isomorphism
\begin{equation*}
    \Phi\colon\{T\}'\cap\mathcal R\to\{B\}',
    \qquad \Phi(T)=B.
\end{equation*}
Moreover, $\Phi^{-1}$ is normal.
\end{corollary}

\begin{proof}
Let $a=\|B^{-1}\|^{-1}$ and $b=\|B\|$.
If $a=b$, then $b^{-1}B$ is a unitary bilateral weighted shift, and the conclusion follows by taking $T=bu$.
In the following, we assume that $a<b$.
Choose pairwise disjoint nowhere dense closed sets
\begin{equation*}
    \{\alpha^jE_n,\alpha^jF_n\colon 0\leqslant j<n,~n\geqslant1\}
\end{equation*}
with $\mu(E_n),\mu(F_n)>0$.
We define
\begin{equation*}
    f(z)=
    \begin{cases}
        a, & z\in\bigsqcup_{n\geqslant1}
        \bigsqcup_{j=0}^{n-1}\alpha^jE_n,\\
        b, & \mbox{otherwise}.
    \end{cases}
\end{equation*}
Let $T=uf(v)$.
Then $\|T^n\|=b^n$ and $\|T^{-n}\|=a^{-n}$ for every $n\geqslant1$.
By \Cref{cor spectrum} and \cite[Theorem~5, p.~67]{Shields-74}, we have $\sigma(B)=\sigma(T)=\{z\in\mathbb C\colon a\leqslant|z|\leqslant b\}$.
For every Laurent polynomial $p$, let
\begin{equation*}
    M_p=\max_{a\leqslant|z|\leqslant b}|p(z)|.
\end{equation*}
By the spectral mapping theorem and \cite[Proposition~23, p.~82]{Shields-74}, there exists $C_{a,b}\geqslant2$ such that
\begin{equation*}
    M_p\leqslant\|p(X)\|\leqslant C_{a,b}M_p,
    \qquad X\in\{B,T\},
\end{equation*}
for every Laurent polynomial $p$.
Consequently,
\begin{equation*}
    C_{a,b}^{-1}\|p(B)\|
    \leqslant\|p(T)\|
    \leqslant C_{a,b}\|p(B)\|
\end{equation*}
for every Laurent polynomial $p$.
The proof of \Cref{lem Laurent-comparison-commutant} applies to this two-sided norm comparison and yields the required unital normal algebra isomorphism with $\|\Phi\|,\|\Phi^{-1}\|\leqslant C_{a,b}$.
Moreover, $\Phi^{-1}$ is also normal.
\end{proof}

The following example shows that \Cref{prop noninvertible-bilateral-realization} cannot be extended to all invertible bilateral weighted shift operators.

\begin{example}\label{eg no-realization}
Let $B$ be an invertible bilateral weighted shift operator defined by
\begin{equation*}
    Be_n=w_ne_{n+1},
    \qquad
    w_n=
    \begin{cases}
        1,&n\leqslant 0,\\
        \dfrac{n}{n+1},& n\geqslant 1.
    \end{cases}
\end{equation*}
Then $\{B\}'$ is not boundedly isomorphic to any unital Banach subalgebra of $\mathcal R$.
In particular, $\{B\}'$ is not boundedly isomorphic to $\{uf(v)\}'\cap\mathcal R$ for any $f\in L^\infty(\mathbb T,\mu)$.
\end{example}

\begin{proof}
Since $1/2\leqslant w_n\leqslant1$, the operator $B$ is invertible.
Moreover, a direct computation shows that
\begin{equation*}
    \|B^n\|=1,
    \qquad
    \|B^{-n}\|=n+1,\qquad n\geqslant 0.
\end{equation*}
Suppose that $\Phi$ is a bounded algebra isomorphism from $\{B\}'$ onto a norm-closed unital subalgebra of $\mathcal R$.
Then $\Phi$ is unital and $\Phi^{-1}$ is bounded.
Let $X=\Phi(B)$.
Since $X^{-1}=\Phi(B^{-1})$, we have
\begin{equation*}
    \|X^n\|\leqslant\|\Phi\|,
    \qquad
    \|X^{-n}\|\leqslant\|\Phi\|(n+1),
    \qquad n\geqslant0.
\end{equation*}
It follows that $r(X)\leqslant1$ and $r(X^{-1})\leqslant1$.
Thus, $\sigma(X)\subseteq\mathbb T$.
Since $X$ is power bounded, the Cassier-Fack theorem \cite[Main Theorem]{CF-96} implies that $X$ is similar to a unitary in $\mathcal R$.
Hence $\sup_{n\geqslant0}\|X^{-n}\|<\infty$.
However,
\begin{equation*}
    n+1=\|B^{-n}\|
    =\|\Phi^{-1}(X^{-n})\|
    \leqslant\|\Phi^{-1}\|\,\|X^{-n}\|,
\end{equation*}
which is a contradiction.
The last assertion follows because $\{uf(v)\}'\cap\mathcal R$ is a unital Banach subalgebra of $\mathcal R$.
\end{proof}

\subsection{Spectrum and point spectrum}\label{subsec spectrum}

As an application of \Cref{thm bilateral-model}, we compute the spectrum of a bounded injective homogeneous operator as follows.
For operators of degree one with continuous weights, see \cite[Theorem~3.3]{ZFS-17}.

\begin{corollary}\label{cor spectrum}
Let $T=u^mf(v)\in\mathcal R$, where $m\ne0$ and $f\ne0$ almost everywhere.
\begin{enumerate}
\item[$(1)$] If $T$ is not invertible, then $\sigma(T)=\{z\in\mathbb C\colon |z|\leqslant r(T)\}$;
\item[$(2)$] If $T$ is invertible, then $\sigma(T)=\{z\in\mathbb C\colon r(T^{-1})^{-1}\leqslant |z|\leqslant r(T)\}$.
\end{enumerate}
\end{corollary}

\begin{proof}
Let $\mathcal R_0=W^*(u^m,v)$.
It is clear that $\mathcal R_0$ is a hyperfinite type $\mathrm{II}_1$ factor and $\sigma(T)=\sigma_{\mathcal R_0}(T)$.
Since $vu^m=\alpha^mu^mv$ and $\alpha^m$ is an irrational rotation parameter, $T=u^mf(v)$ is $1$-homogeneous with respect to the generators $u^m,v$ of $\mathcal R_0$.
Without loss of generality, we assume that $m=1$.
Then $\{T\}'\cap\mathcal R\cong\{B\}'$ for some injective bilateral weighted shift operator as in \Cref{thm bilateral-model}.
It follows that $\sigma(T)=\sigma(B)$ and the conclusion follows from \cite[Theorem~5, p.~67]{Shields-74}.
\end{proof}

Let $T$ be a measurable operator in $L^0(\mathcal R,\tau)$.
A scalar $\lambda\in\mathbb C$ is called an \emph{eigenvalue} of $T$ if $T-\lambda$ has nontrivial kernel.
The \emph{point spectrum} $\sigma_p(T)$ is the set of all eigenvalues of $T$.

\begin{corollary}\label{cor point-spectrum}
Let $T=u^mf(v)\in L^0(\mathcal R,\tau)$ with $m\ne 0$.
Then
\begin{equation*}
    \sigma_p(T)=
    \begin{cases}
        \varnothing, & f\ne0~\text{almost everywhere},\\
        \{0\}, & \mbox{otherwise}.
    \end{cases}
\end{equation*}
\end{corollary}

\begin{proof}
Since $\ker T=\ker f(v)$, it follows that $0\notin\sigma_p(T)$ if and only if $f\ne 0$ almost everywhere.
Thus, it remains to exclude nonzero eigenvalues.
By scaling, it suffices to show that $1\notin\sigma_p(T)$.
Let $p\in\mathcal R$ be the projection onto $\ker(T-1)$.
Then $Tp=p$, and we shall prove that $p=0$.

Write $p=\sum_{n\in\mathbb Z}u^ng_n(v)\in\mathcal R$.
Since $g_0(v)=E_{W^*(v)}(p)$, we have $0\leqslant g_0\leqslant1$.
Since $Tp=p$ and $p=p^*$, it follows from \Cref{lem frequently} that
\begin{equation*}
	g_{m+n}(z)=f(\alpha^nz)g_n(z),\qquad g_n(z)=\overline{g_{-n}(\alpha^nz)},\qquad n\in\mathbb Z.
\end{equation*}
A direct computation shows that
\begin{equation*}
	g_0(\alpha^mz)=f(z)g_{-m}(\alpha^mz)=f(z)\overline{g_m(z)}=|f(z)|^2g_0(z).
\end{equation*}
Inductively, we have
\begin{equation*}
	g_0(\alpha^{mn}z)=|f(z)f(\alpha^mz)\cdots f(\alpha^{m(n-1)}z)|^2g_0(z),\qquad n\geqslant 1.
\end{equation*}
Similarly, by $g_{m+n}(z)=f(\alpha^nz)g_n(z)$, we have
\begin{equation*}
	g_{mn}(z)=f(z)f(\alpha^mz)\cdots f(\alpha^{m(n-1)}z)g_0(z),\qquad n\geqslant 1.
\end{equation*}
It follows that
\begin{equation*}
    |g_{mn}(z)|^2=g_0(z)g_0(\alpha^{mn}z),
    \qquad n\geqslant1.
\end{equation*}
Let $\{n_k\}_{k=1}^\infty$ be an increasing sequence of positive integers such that $\alpha^{mn_k}\to1$ as $k\to\infty$.
The Cauchy-Schwarz inequality shows that
\begin{align*}
	\left|\|g_{mn_k}\|_2^2-\|g_0\|_2^2\right|
	&=\left|\int_{\mathbb T}g_0(z)\bigl(g_0(\alpha^{mn_k}z)-g_0(z)\bigr)d\mu(z)\right|\\
    &\leqslant\|g_0\|_2\|g_0(\alpha^{mn_k}\cdot)-g_0\|_2
    \to0.
\end{align*}
Since $\sum_{n\in\mathbb Z}\|g_n\|_2^2=\|p\|_2^2<\infty$, we have $g_0=0$ and $\tau(p)=0$.
Thus, $p=0$.
\end{proof}

Recall that, by \cite[Proposition~31, p.~94]{Shields-74}, an injective unilateral weighted shift $S$ is \emph{strictly cyclic} if and only if there exists $C>0$ such that $\|p(S)\|\leqslant C\|p(S)e_0\|$ for every polynomial $p$.
We show that the analogous estimate $\|p(T)\|\leqslant C\|p(T)\|_2$ cannot hold for $T=uf(v)\in\mathcal R$.

\begin{corollary}\label{cor no-strictly-cyclic}
Let $T=uf(v)\in\mathcal R$, where $f\ne 0$ almost everywhere.
Then there exists a sequence of polynomials $\{p_j\}_{j=1}^\infty$ such that
\begin{equation*}
    \|p_j(T)\|=1
    \qquad\text{and}\qquad
    \lim_{j\to\infty}\|p_j(T)\|_2=0.
\end{equation*}
\end{corollary}

\begin{proof}
Choose $\lambda\in\sigma(T)$ with $|\lambda|=r(T)$.
Let $\{\lambda_j\}_{j=1}^{\infty}$ be a sequence of scalars such that $|\lambda_j|>r(T)$ for each $j\geqslant 1$ and $\lambda_j\to\lambda$ as $j\to\infty$.
Let
\begin{equation*}
    X_j=(\lambda_j-T)^{-1},
    \qquad
    Y_j=\frac{X_j}{\|X_j\|}.
\end{equation*}
Since $\lambda\in\sigma(T)$, we have $\|X_j\|\geqslant|\lambda_j-\lambda|^{-1}$, and hence
\begin{equation*}
    \|(\lambda-T)Y_j\|
    \leqslant|\lambda-\lambda_j|+\|X_j\|^{-1}
    \leqslant 2|\lambda-\lambda_j|
    \to0.
\end{equation*}
By \Cref{cor point-spectrum}, $\lambda-T$ has trivial kernel.
For $\varepsilon>0$, let
\begin{equation*}
    E_\varepsilon=\chi_{[0,\varepsilon]}(|\lambda-T|).
\end{equation*}
Then $\tau(E_\varepsilon)\to0$ as $\varepsilon\to0$.
Since $\|Y_j\|=1$, we have
\begin{equation*}
	\|Y_j\|_2^2=\|E_\varepsilon Y_j\|_2^2+\|(1-E_\varepsilon)Y_j\|_2^2
	\leqslant\tau(E_\varepsilon)+\varepsilon^{-2}\|(\lambda-T)Y_j\|_2^2.
\end{equation*}
It follows that $\limsup_{j\to\infty}\|Y_j\|_2^2\leqslant\tau(E_\varepsilon)$.
Therefore, $\|Y_j\|_2\to0$.
Since $|\lambda_j|>r(T)$, the series
\begin{equation*}
    X_j=\sum_{n=0}^{\infty}\lambda_j^{-n-1}T^n
\end{equation*}
converges in operator norm.
Consequently, there exists a polynomial $q_j$ such that
$\|q_j(T)-Y_j\|<2^{-j}$.
Let $p_j=\|q_j(T)\|^{-1}q_j$.
Then $\|p_j(T)\|=1$ and
\begin{equation*}
    \|p_j(T)\|_2
    \leqslant
    \frac{\|Y_j\|_2+2^{-j}}{\|Y_j\|-2^{-j}}
    \leqslant
    \frac{\|Y_j\|_2+2^{-j}}{1-2^{-j}}
    \to0.
\end{equation*}
This completes the proof.
\end{proof}

\section{Unilateral weighted shift operators}\label{sec unilateral}

For an injective noninvertible homogeneous operator of degree one in $\mathcal R$, each operator in its relative commutant has an $L^2$-Fourier expansion involving only nonnegative degrees.
This allows us to replace the bilateral weighted shift model of the preceding section by unilateral weighted shift operators.
We also establish a converse realization result.

\subsection{Unilateral weighted shift model}

A \emph{unilateral weighted shift operator} $S$ on $\ell^2(\mathbb N)$ is defined by
\begin{equation*}
    Se_n=w_ne_{n+1},\qquad n\geqslant0,
\end{equation*}
where $\{e_n\}_{n=0}^\infty$ is the standard orthonormal basis of $\ell^2(\mathbb N)$ and the weight sequence $\{w_n\}_{n=0}^\infty$ is bounded.
In this subsection, we establish the unilateral weighted shift model for $\{uf(v)\}'\cap\mathcal R$, where $f\in L^\infty(\mathbb T,\mu)$ is not invertible and $f\ne 0$ almost everywhere.

For $f\in L^\infty(\mathbb T,\mu)$ and $z\in\mathbb T$, let $S_{f,z}$ be the unilateral weighted shift operator defined by
\begin{equation}\label{equ S-f-z}
    S_{f,z}e_n=f(\alpha^nz)e_{n+1},
    \qquad n\geqslant 0.
\end{equation}
The operators $S_{f,z}$ are restrictions of the bilateral fibers introduced in \eqref{equ B-f-z}, that is,
\begin{equation*}
    S_{f,z}=B_{f,z}|_{\ell^2(\mathbb N)}.
\end{equation*}
Let $T=uf(v)$, $\mathcal V=W^*(v)$, and $\mathcal H_+=\bigoplus_{n=0}^{\infty}u^nL^2(\mathcal V,\tau)\subseteq L^2(\mathcal R,\tau)$.
Then the subspace $\mathcal H_+$ is invariant under left multiplication by $T$ and there exists a unitary operator
\begin{equation*}
	W_+\colon\mathcal H_+\to L^2(\mathbb T,\ell^2(\mathbb N))
\end{equation*}
defined by $\bigl(W_+u^ng_n(v)\bigr)(z)=g_n(z)e_n$ for $n\geqslant 0$.
Let $T_+=L_T|_{\mathcal H_+}$.
Similar to \eqref{equ T-direct-integral}, we have
\begin{equation}\label{equ unilateral-direct-integral}
	L^2(\mathbb T,\ell^2(\mathbb N))
    =
    \int_{\mathbb T}^{\oplus}\ell^2(\mathbb N)\,d\mu(z),\qquad
    W_+T_+W_+^*
    =
    \int_{\mathbb T}^{\oplus}S_{f,z}\,d\mu(z).
\end{equation}
We first show that passing from bilateral fibers to unilateral restrictions preserves polynomial norms almost everywhere as in \Cref{lem fiber-Laurent-norm}.

\begin{lemma}\label{lem fiber-polynomial-norm}
With the notation of \eqref{equ S-f-z} and $T=uf(v)\in\mathcal R$, there exists
a null set $E_f\subseteq\mathbb T$ such that
\begin{equation*}
    \|p(uf(v))\|=\|p(S_{f,z})\|,
    \qquad z\in\mathbb T\setminus E_f,
\end{equation*}
for every polynomial $p$.
\end{lemma}

\begin{proof}
Let $h(z)=\|p(S_{f,z})\|$.
Since $S_{f,\alpha z}$ is unitarily equivalent to the restriction of $S_{f,z}$ to $\overline{\operatorname{span}}\{e_n\colon n\geqslant1\}$, we have $h(\alpha z)\leqslant h(z)$.
Since
\begin{equation*}
    \int_{\mathbb T}\bigl(h(z)-h(\alpha z)\bigr)\,d\mu(z)=0,
\end{equation*}
we see that $h(\alpha z)=h(z)$ almost everywhere.
Thus, $h$ is essentially constant.

For each $k\geqslant0$, the restriction of $B_{f,z}$ to
$\overline{\operatorname{span}}\{e_n\colon n\geqslant-k\}$
is unitarily equivalent to $S_{f,\alpha^{-k}z}$.
Thus, we have
\begin{equation*}
    \|p(B_{f,z})\|
    =
    \sup_{k\geqslant0}\|p(S_{f,\alpha^{-k}z})\|
    =
    \|p(S_{f,z})\|
\end{equation*}
almost everywhere.
By applying \eqref{equ unilateral-direct-integral}, the conclusion follows from the bilateral fiber-norm formula and the same countable approximation argument as in \Cref{lem fiber-Laurent-norm}.
\end{proof}

As in the bilateral case \Cref{lem Laurent-comparison-commutant}, a comparison of polynomial norms gives an isomorphism between the corresponding commutants.

\begin{lemma}\label{lem unilateral-comparison-commutant}
Let $T=uf(v)\in\mathcal R$, where $f\in L^\infty(\mathbb T,\mu)$ is not invertible and $f\ne 0$ almost everywhere, and let $S$ be an injective unilateral weighted shift operator.
Suppose that there exists $C\geqslant1$ such that
\begin{equation*}
    \|p(S)\|\leqslant\|p(T)\|
    \leqslant C\|p(S)\|
\end{equation*}
for every polynomial $p$.
Then there exists a unique unital normal algebra isomorphism
\begin{equation*}
    \Phi\colon\{T\}'\cap\mathcal R\to\{S\}',
    \qquad \Phi(T)=S.
\end{equation*}
Moreover, $\|\Phi\|=1$, $\|\Phi^{-1}\|\leqslant C$,
and $\Phi^{-1}$ is normal.
\end{lemma}

\begin{proof}
The proof of \Cref{lem Laurent-comparison-commutant} applies with nonnegative powers in place of integer powers.
\end{proof}

Combining the preceding two lemmas yields the desired unilateral weighted shift model.

\begin{theorem}\label{thm unilateral-model}
Let $T=uf(v)\in\mathcal R$, where $f\in L^\infty(\mathbb T,\mu)$ is not invertible and $f\ne 0$ almost everywhere.
Then there exists a null set $E_f\subseteq\mathbb T$
such that, for every $z\in\mathbb T\setminus E_f$, the unilateral weighted shift operator
\begin{equation*}
	S_{f,z}e_n=f(\alpha^nz)e_{n+1},\qquad n\geqslant 0
\end{equation*}
is bounded and injective, and
there exists a unique unital normal isometric
algebra isomorphism
\begin{equation*}
    \Phi\colon\{T\}'\cap\mathcal R
    \to\{S_{f,z}\}',
    \qquad \Phi(T)=S_{f,z}.
\end{equation*}
Moreover, $\Phi^{-1}$ is normal.
\end{theorem}

\begin{proof}
Similar to the proof of \Cref{thm bilateral-model}, the conclusion follows from \Cref{lem fiber-polynomial-norm,lem unilateral-comparison-commutant}.
\end{proof}

If the weight function $f$ is continuous in \Cref{thm unilateral-model}, then we have the following result analogous to \Cref{cor bilateral-model-continuous}.

\begin{corollary}\label{cor unilateral-model-continuous}
Let $T=uf(v)\in\mathcal R$, where $f\in C(\mathbb T)$, $f\ne 0$ almost everywhere, and $f$ has a zero.
Then for every $z\in\mathbb T$ with $f(\alpha^nz)\ne0$ for all $n\geqslant0$, there exists
a unique unital normal isometric algebra isomorphism
\begin{equation*}
    \Phi\colon\{T\}'\cap\mathcal R
    \to\{S_{f,z}\}',
    \qquad \Phi(T)=S_{f,z}.
\end{equation*}
Moreover, $\Phi^{-1}$ is normal.
\end{corollary}

\begin{proof}
As in \Cref{cor bilateral-model-continuous}, uniform continuity of $f$ implies that $z\mapsto S_{f,z}$ is operator-norm continuous.
Hence the polynomial norm equality in \Cref{lem fiber-polynomial-norm} holds for every $z\in\mathbb T$.
Apply \Cref{lem unilateral-comparison-commutant} with $C=1$.
\end{proof}

\begin{remark}
If $f$ is invertible, no unital algebra isomorphism from $\{uf(v)\}'\cap\mathcal R$ onto the commutant of a unilateral weighted shift $S$ can send $uf(v)$ to $S$.
Indeed, such an isomorphism would send $(uf(v))^{-1}$ to an inverse of $S$, whereas $S$ is not invertible.
\end{remark}

\subsection{Unilateral weighted shift realization}

We now consider the converse realization problem.
The construction from the bilateral case also applies to injective unilateral weighted shift operators.
We begin with the corresponding approximation by finite compressions as in \Cref{lem finite-shift-blocks}.

\begin{lemma}\label{lem unilateral-finite-shift-blocks}
Let $S$ be an injective unilateral weighted shift operator.
For every integer $m\geqslant1$, let $S_m$ be its compression to the $m$-dimensional subspace
\begin{equation*}
    H_m=\operatorname{span}\{e_0,\ldots,e_{m-1}\}.
\end{equation*}
Then for every polynomial $p$, we have
\begin{equation*}
    \|p(S_m)\|\leqslant\|p(S)\|,
    \qquad
    \lim_{m\to\infty}\|p(S_m)\|=\|p(S)\|.
\end{equation*}
Moreover, for any contractions $C$ and $D$ with $C^m=D^m=0$, we have
\begin{equation*}
    \|p(C)-p(D)\|
    \leqslant\kappa_m\|C-D\|\|p(S)\|,
    \qquad
    \kappa_m=\sum_{j=1}^{m-1}\frac{j}{\|S^j\|}.
\end{equation*}
\end{lemma}

\begin{proof}
The proof of \Cref{lem finite-shift-blocks} applies.
\end{proof}

Using these finite compressions and the construction from \Cref{lem filling-zero-weights}, we obtain the following comparison of polynomial norms as in \Cref{lem realization-B}.

\begin{lemma}\label{lem realization-S}
Let $S$ be an injective unilateral weighted shift operator.
Then for every $\varepsilon>0$, there exists $f\in L^\infty(\mathbb T,\mu)$, nonzero almost everywhere and not invertible in $L^\infty(\mathbb T,\mu)$, such that
\begin{equation*}
    \|p(S)\|\leqslant\|p(uf(v))\|
    \leqslant(1+\varepsilon)\|p(S)\|
\end{equation*}
for every polynomial $p$.
\end{lemma}

\begin{proof}
By scaling, we assume that $\|S\|=1$.
Write $Se_n=w_ne_{n+1}$ for $n\geqslant0$.
Choose a nowhere dense closed set $D\subseteq\mathbb T$ of positive measure, and let $D_m$ be as in \Cref{lem return-time-partition}.
Define $g$ almost everywhere by
\begin{equation*}
    g(z)=
    \begin{cases}
        0,&z\in D,\\[2pt]
        w_{j-1},&z\in\alpha^jD_m,\qquad
        1\leqslant j<m,\qquad m\geqslant1.
    \end{cases}
\end{equation*}
Then $|g|\leqslant1$ and $Z(g)=D$ modulo a null set.
As in the proof of \Cref{lem realization-B}, almost every $S_{g,z}$ is an orthogonal direct sum of finite blocks $S_m$, except possibly for its first block, which is a compression of some $S_m$ to an interval.
By \Cref{lem fiber-polynomial-norm,lem unilateral-finite-shift-blocks}, we obtain
\begin{equation*}
    \|p(ug(v))\|\leqslant\|p(S)\|.
\end{equation*}
The proof of \Cref{lem filling-zero-weights} applies unchanged with $S$ in place of $B$, and produces $f$, nonzero almost everywhere and noninvertible, satisfying the required upper bound.
The internal weights of the original blocks remain unchanged, giving
\begin{equation*}
    \|p(uf(v))\|\geqslant\|p(S_m)\|
    \qquad\text{whenever }\mu(D_m)>0.
\end{equation*}
These indices are unbounded, as shown in \Cref{lem realization-B}.
The lower bound follows from \Cref{lem unilateral-finite-shift-blocks}.
\end{proof}

We now prove the main realization result of this subsection: the commutant of every injective unilateral weighted shift admits a nearly isometric realization as the relative commutant of a homogeneous operator in $\mathcal R$ analogous to \Cref{prop noninvertible-bilateral-realization}.

\begin{proposition}\label{prop unilateral-realization}
Let $S$ be an injective unilateral weighted shift operator.
Then for every $\varepsilon>0$, there exist $T=uf(v)\in\mathcal R$, where $f$ is nonzero almost everywhere and noninvertible, and a unital normal algebra isomorphism
\begin{equation*}
    \Phi\colon\{T\}'\cap\mathcal R\to\{S\}',
    \qquad \Phi(T)=S,
\end{equation*}
such that
\begin{equation*}
    \|\Phi\|=1,
    \qquad
    \|\Phi^{-1}\|\leqslant1+\varepsilon.
\end{equation*}
Moreover, $\Phi^{-1}$ is normal.
\end{proposition}

\begin{proof}
Apply \Cref{lem realization-S} and \Cref{lem unilateral-comparison-commutant} with $C=1+\varepsilon$.
\end{proof}

Under the additional assumption $r(S)=\|S\|$, the isomorphism can be chosen to be isometric.

\begin{corollary}\label{cor isometric-unilateral-realization}
Let $S$ be an injective unilateral weighted shift operator such that $r(S)=\|S\|$.
Then there exist $T=uf(v)\in\mathcal R$, where $f$
is nonzero almost everywhere and noninvertible, and a unital normal isometric algebra isomorphism
\begin{equation*}
    \Phi\colon\{T\}'\cap\mathcal R\to\{S\}',
    \qquad \Phi(T)=S.
\end{equation*}
Moreover, $\Phi^{-1}$ is normal.
\end{corollary}

\begin{proof}
The conclusion follows from \Cref{lem unilateral-comparison-commutant} and the proof of \Cref{cor isometric-bilateral-realization}.
\end{proof}

\subsection{A version of the corona theorem}

The commutant of the unweighted unilateral shift is naturally identified with $H^\infty(\mathbb D)$.
Under the same norm and spectral radius condition, this gives an explicit description of the relative commutant of $T$.

\begin{corollary}\label{cor commutant-Hinfty}
Let $T=uf(v)\in\mathcal R$, where $f\ne0$ almost everywhere.
If $T$ is not invertible and $r(T)=\|T\|$, then there exists a unital normal isometric algebra isomorphism
\begin{equation*}
    \Phi\colon\{T\}'\cap\mathcal R
    \to H^\infty(\mathbb D),
    \qquad \Phi(T/\|T\|)=z.
\end{equation*}
Moreover, $\Phi^{-1}$ is normal.
\end{corollary}

\begin{proof}
By scaling, we assume that $r(T)=\|T\|=1$.
Let $S$ be the unweighted unilateral shift.
Then $\sigma(S)=\overline{\mathbb D}=\sigma(T)$ by \Cref{cor spectrum}.
Thus, the spectral mapping theorem and von Neumann's inequality give
\begin{equation*}
    \|p(S)\|
    =\max_{|z|\leqslant1}|p(z)|
    =\|p(T)\|
\end{equation*}
for every polynomial $p$.
The proof is completed by \Cref{lem unilateral-comparison-commutant} and identifying the commutant $\{S\}'$ with $H^\infty(\mathbb D)$.
\end{proof}

\begin{remark}\label{rem commutant-Hinfty}
Let $K\subseteq\mathbb T$ be a nowhere dense closed set of positive measure, and partition $K$ into measurable sets $\{K_j\}_{j=1}^\infty$ of positive measure.
Define
\begin{equation*}
    f(z)=
    \begin{cases}
        2^{-j}, & z\in K_j,\\
        1, & z\notin K.
    \end{cases}
\end{equation*}
Then $T=uf(v)$ is injective and noninvertible, and $\|T^n\|=1$ for every $n\geqslant 1$.
In particular, $\|T\|=r(T)=1$.
Thus, $T$ satisfies the assumption in \Cref{cor commutant-Hinfty}.
\end{remark}

The identification with $H^\infty(\mathbb D)$ allows us to transfer Carleson's corona theorem to the relative commutant $\{T\}'\cap\mathcal R$.
Before stating the result, we record an elementary comparison of positive operators.
For finitely many operators $\{x_j\}_{j=1}^n$, let $C=\max\{\|x_j\|\colon 1\leqslant j\leqslant n\}$.
Since
\begin{equation*}
	n\sum_{j=1}^n|x_j|^2-\left(\sum_{j=1}^n|x_j|\right)^2=\sum_{i<j}(|x_i|-|x_j|)^2\geqslant0
\end{equation*}
and $|x_j|^2\leqslant C|x_j|$ for $1\leqslant j\leqslant n$, we have
\begin{equation*}
	\sum_{j=1}^n|x_j|\leqslant\left(n\sum_{j=1}^n|x_j|^2\right)^{1/2},\qquad
	\sum_{j=1}^n|x_j|^2\leqslant C\sum_{j=1}^n|x_j|.
\end{equation*}
Thus, $\sum_{j=1}^n|x_j|$ is bounded below if and only if $\sum_{j=1}^n|x_j|^2$ is.
We now establish the following version of the corona theorem.

\begin{proposition}\label{prop corona}
Let $T=uf(v)\in\mathcal R$, where $f\in L^\infty(\mathbb T,\mu)$ is not invertible and $f\ne 0$ almost everywhere.
Suppose that $r(T)=\|T\|=1$ and $x_1,\ldots,x_n\in\{T\}'\cap\mathcal R$.
Then the following statements are equivalent:
\begin{enumerate}
\item[$(1)$] there exists $c>0$ such that $\sum_{j=1}^nx_j^*x_j\geqslant c$;
\item[$(2)$] there exists $c>0$ such that $\sum_{j=1}^nx_jx_j^*\geqslant c$;
\item[$(3)$] there exist $y_1,\ldots,y_n\in\{T\}'\cap\mathcal R$ such that $\sum_{j=1}^nx_jy_j=1$.
\end{enumerate}
\end{proposition}

\begin{proof}
By \Cref{cor commutant-Hinfty}, there is a unital normal isometric algebra isomorphism
\begin{equation*}
    \Phi\colon\{T\}'\cap\mathcal R\to H^\infty(\mathbb D),
    \qquad \Phi(T)=z.
\end{equation*}
Let $\varphi_j=\Phi(x_j)\in H^\infty(\mathbb D)$ for $1\leqslant j\leqslant n$.

We first prove that either (1) or (2) implies (3).
Since $\{T\}'\cap\mathcal R$ is inverse closed in $\mathcal R$, we have
\begin{equation*}
    \sigma_{\mathcal R}(T)
    =\sigma_{\{T\}'\cap\mathcal R}(T)
    =\sigma_{H^\infty(\mathbb D)}(z)
    =\overline{\mathbb D}.
\end{equation*}
Fix $\lambda\in\mathbb D$.
Since both $T-\lambda$ and $T^*-\overline\lambda$ are noninvertible and every bounded below operator is invertible in $\mathcal R$, there are unit vectors $\{\xi_k\}_{k=1}^\infty,\{\eta_k\}_{k=1}^\infty$ such that
\begin{equation*}
	\|(T-\lambda)\xi_k\|\to 0,\qquad \|(T^*-\overline\lambda)\eta_k\|\to 0.
\end{equation*}
For each $1\leqslant j\leqslant n$, we have
\begin{equation*}
	\varphi_{j,\lambda}=\frac{\varphi_j(z)-\varphi_j(\lambda)}{z-\lambda}\in H^\infty(\mathbb D).
\end{equation*}
Let $x_{j,\lambda}=\Phi^{-1}(\varphi_{j,\lambda})\in\{T\}'\cap\mathcal R$.
Then
\begin{equation*}
    x_j-\varphi_j(\lambda)
    =x_{j,\lambda}(T-\lambda),
    \qquad
    x_j^*-\overline{\varphi_j(\lambda)}
    =x_{j,\lambda}^*(T^*-\overline\lambda).
\end{equation*}
It follows that
\begin{equation*}
    \|x_j\xi_k-\varphi_j(\lambda)\xi_k\|
    \to0,
    \qquad
    \|x_j^*\eta_k-\overline{\varphi_j(\lambda)}\eta_k\|
    \to0.
\end{equation*}
If (1) holds, then
\begin{equation*}
    \sum_{j=1}^n|\varphi_j(\lambda)|^2=
    \lim_{k\to\infty}\sum_{j=1}^n\|x_j\xi_k\|^2
    \geqslant c.
\end{equation*}
If (2) holds, the same conclusion follows from
\begin{equation*}
    \sum_{j=1}^n|\varphi_j(\lambda)|^2=
    \lim_{k\to\infty}\sum_{j=1}^n\|x_j^*\eta_k\|^2
    \geqslant c.
\end{equation*}
In either case, we have
\begin{equation*}
    \sum_{j=1}^n|\varphi_j(\lambda)|
    \geqslant
    \left(\sum_{j=1}^n|\varphi_j(\lambda)|^2\right)^{1/2}
    \geqslant\sqrt c
    \qquad(\lambda\in\mathbb D).
\end{equation*}
By Carleson's corona theorem \cite{Carleson-62}, there exist
$\psi_1,\ldots,\psi_n\in H^\infty(\mathbb D)$ such that
\begin{equation*}
    \sum_{j=1}^n\varphi_j\psi_j=1.
\end{equation*}
Taking $y_j=\Phi^{-1}(\psi_j)$ proves (3).

Conversely, suppose that (3) holds.
Since $\{T\}'\cap\mathcal R$ is commutative, we have
\begin{equation*}
    \sum_{j=1}^nx_jy_j
    =\sum_{j=1}^ny_jx_j
    =1.
\end{equation*}
Let $C=\sum_{j=1}^n\|y_j\|^2>0$.
For every vector $\xi$, we have
\begin{equation*}
    \|\xi\|^2
    =\left\|\sum_{j=1}^ny_jx_j\xi\right\|^2
    \leqslant\left(\sum_{j=1}^n\|y_j\|\|x_j\xi\|\right)^2
    \leqslant C\sum_{j=1}^n\|x_j\xi\|^2.
\end{equation*}
Taking adjoints of $\sum_{j=1}^nx_jy_j=1$ also gives
\begin{equation*}
    \|\xi\|^2
    =\left\|\sum_{j=1}^ny_j^*x_j^*\xi\right\|^2
    \leqslant C\sum_{j=1}^n\|x_j^*\xi\|^2.
\end{equation*}
Consequently, we have
\begin{equation*}
    \sum_{j=1}^nx_j^*x_j\geqslant C^{-1},
    \qquad
    \sum_{j=1}^nx_jx_j^*\geqslant C^{-1}.
\end{equation*}
Thus both (1) and (2) hold.
\end{proof}

\section{\texorpdfstring{$K_0$}{K0}-groups of relative commutants}\label{sec K0}

In this section, we study the $K_0$-groups of relative commutants of injective homogeneous operators and their direct sums.
We first identify a direct summand generated by the identity and compute the groups in two analytic models in \Cref{subsec K0=Z}.
Then we establish splitting criteria for direct sums in \Cref{prop different-degree-K0-split} and \Cref{thm K0-split}.

\subsection{\texorpdfstring{$K_0$}{K0}-groups of homogeneous operators}
\label{subsec K0=Z}

We begin by showing that the class of the identity always generates a direct summand isomorphic to $\mathbb Z$.

\begin{lemma}\label{lem commutant-K0-Z-summand}
Let $T=u^mf(v)\in L^0(\mathcal R,\tau)$, where $m\ne0$ and $f\ne0$ almost everywhere.
Then $\mathbb Z[1]$ is a direct summand of $K_0(\{T\}'\cap\mathcal R)$.
\end{lemma}

\begin{proof}
By \Cref{cor commutant-maximal-abelian}, $\{T\}'\cap\mathcal R$ is a commutative unital Banach algebra.
Choose a character $\varphi\colon\{T\}'\cap\mathcal R\to\mathbb C$, and let $\iota\colon\mathbb C\to\{T\}'\cap\mathcal R$ be the scalar inclusion.
Since $\varphi\iota=\operatorname{id}_{\mathbb C}$, the induced homomorphisms satisfy $\varphi_*\iota_*=\operatorname{id}_{K_0(\mathbb C)}$.
Identifying $K_0(\mathbb C)$ with $\mathbb Z$, we have
$\iota_*(n)=n[1]$.
Thus, $\iota_*$ is injective and
\begin{equation*}
    K_0(\{T\}'\cap\mathcal R)
    =
    \mathbb Z[1]\oplus\ker\varphi_*.
\end{equation*}
This completes the proof.
\end{proof}

\begin{remark}\label{rem commutant-K0-Z-summand}
Let $T=u^mf(v)\in\mathcal R$, where $m\ne0$ and $f\ne0$ almost everywhere.
If $f$ is not invertible in $L^{\infty}(\mathbb T,\mu)$, then $\tau$ is multiplicative on $\{T\}'\cap\mathcal R$.
In this case, we have the canonical splitting $K_0(\{T\}'\cap\mathcal R)=\mathbb Z[1]\oplus\ker\tau_*$.
\end{remark}

Using \Cref{cor commutant-Hinfty}, we obtain the following condition under which the $K_0$-group is generated by the class of the identity.

\begin{corollary}\label{cor H-infty-K0=Z}
Let $T=uf(v)\in\mathcal R$, where $f\ne0$ almost everywhere.
If $T$ is not invertible and $r(T)=\|T\|$, then
\begin{equation*}
	V(\{T\}'\cap\mathcal R)\cong\mathbb Z_{\geqslant 0},\qquad K_0(\{T\}'\cap\mathcal R)\cong\mathbb Z,
\end{equation*}
with $[1]$ corresponding to $1$.
\end{corollary}

\begin{proof}
Since $T\ne0$, we have $r(T)=\|T\|>0$.
Applying \Cref{cor commutant-Hinfty} to $\|T\|^{-1}T$,
we obtain a unital algebra isomorphism
\begin{equation*}
    \{T\}'\cap\mathcal R\cong H^\infty(\mathbb D).
\end{equation*}
Since $H^\infty(\mathbb D)$ is projective free by \cite[Theorem~1.5]{BS-09}, we have $V(\{T\}'\cap\mathcal R)\cong\mathbb Z_{\geqslant 0}$, with $[1]$ corresponding to $1$.
Taking the Grothendieck group gives the asserted
$K_0$-isomorphism.
\end{proof}

For invertible operators, an identification with bounded holomorphic functions on an annulus yields an analogue of \Cref{cor H-infty-K0=Z}.

\begin{corollary}\label{cor invertible-commutant-Hinfty}
Let $T=uf(v)\in\mathcal R$, where $f$ is invertible in
$L^\infty(\mathbb T,\mu)$.
Suppose that
\begin{equation*}
    r(T)=\|T\|=b,\qquad
    r(T^{-1})=\|T^{-1}\|=a^{-1},
    \qquad 0<a<b.
\end{equation*}
Then there exists a unique unital normal algebra isomorphism
\begin{equation*}
    \Phi\colon\{T\}'\cap\mathcal R\to H^\infty(\Omega),
    \qquad \Phi(T)=z,
\end{equation*}
where $\Omega=\{z\in\mathbb C\colon a<|z|<b\}$.
Moreover, $\Phi$ and $\Phi^{-1}$ are bounded, and
$\Phi^{-1}$ is normal.
In particular, we have
\begin{equation*}
    V(\{T\}'\cap\mathcal R)\cong\mathbb Z_{\geqslant 0},\qquad
    K_0(\{T\}'\cap\mathcal R)\cong\mathbb Z,
\end{equation*}
with $[1]$ corresponding to $1$.
\end{corollary}

\begin{proof}
By \Cref{cor spectrum}, $\sigma(T)=\overline\Omega$.
As in the proof of \Cref{cor invertible-bilateral-realization}, by the spectral mapping theorem and \cite[Proposition~23, p.~82]{Shields-74}, there exists $C_{a,b}\geqslant2$ such that
\begin{equation*}
    \|p\|_{H^\infty(\Omega)}
    \leqslant\|p(T)\|
    \leqslant C_{a,b}\|p\|_{H^\infty(\Omega)}
\end{equation*}
for every Laurent polynomial $p$.
The proof of \Cref{lem Laurent-comparison-commutant} applies here, using Fej\'er sums in $H^\infty(\Omega)$.
It yields a bounded unital normal algebra isomorphism
\begin{equation*}
    \Phi\colon\{T\}'\cap\mathcal R\to H^\infty(\Omega),
    \qquad \Phi(T)=z,
\end{equation*}
whose inverse is also normal.
Uniqueness follows from Fej\'er approximation.
Finally, $H^\infty(\Omega)$ is projective free by \cite[Theorem~1.5]{BS-09}.
Therefore, $V(\{T\}'\cap\mathcal R)\cong\mathbb Z_{\geqslant 0}$, with $[1]$ corresponding to $1$.
Taking the Grothendieck group gives $K_0(\{T\}'\cap\mathcal R)\cong\mathbb Z$.
\end{proof}

\subsection{A splitting principle}

We now turn to direct sums.
The following algebraic lemma allows us to remove a square-zero ideal without changing the $V$-monoid or the $K_0$-group.
Here $\mathcal J^2=\{0\}$ means that $xy=0$ for all $x,y\in\mathcal J$.

\begin{lemma}\label{lem square-zero-K0-splitting}
Let $\mathcal A$ be a unital algebra, $\mathcal B$ a unital subalgebra of $\mathcal A$, and $\mathcal J$ a two-sided ideal of $\mathcal A$ with $\mathcal J^2=\{0\}$ such that
\begin{equation*}
    \mathcal A=\mathcal B\oplus\mathcal J
\end{equation*}
as a vector space.
Then the inclusion $\iota\colon\mathcal B\to\mathcal A$ induces isomorphisms
\begin{equation*}
    V(\mathcal B)\cong V(\mathcal A),
    \qquad
    K_0(\mathcal B)\cong K_0(\mathcal A).
\end{equation*}
\end{lemma}

\begin{proof}
Let $\pi\colon\mathcal A\to\mathcal B$ be the projection along $\mathcal J$.
Since $\mathcal J$ is an ideal, $\pi$ is a unital homomorphism and $\pi\iota=\operatorname{id}$.
For an idempotent $P\in M_n(\mathcal A)$, let
\begin{equation*}
    P_0=\iota_n(\pi_n(P)),\qquad Q=P-P_0\in M_n(\mathcal J).
\end{equation*}
Then $P_0^2=P_0$, $Q^2=0$, and
\begin{equation*}
    P_0Q+QP_0=Q,\qquad P_0QP_0=0.
\end{equation*}
Let $Z=QP_0-P_0Q$. Since $Z\in M_n(\mathcal J)$, we have $Z^2=0$, and direct multiplication gives
\begin{equation*}
    (I_n+Z)P_0(I_n-Z)=P_0+Q=P.
\end{equation*}
Thus every matrix idempotent over $\mathcal A$ is similar to its image under $\iota_n\pi_n$.
Consequently, $\iota_*\pi_*=\operatorname{id}$ on $V(\mathcal A)$, while $\pi_*\iota_*=\operatorname{id}$ on $V(\mathcal B)$. Taking Grothendieck groups gives the asserted $K_0$-isomorphism.
\end{proof}

For the direct sum of two operators, the off-diagonal entries of its relative commutant are intertwiners between the two summands.
If products of intertwiners in opposite directions vanish, these entries form a square-zero ideal.
The preceding lemma therefore gives the following splitting principle.

\begin{proposition}\label{prop splitting-principle}
For $T_1,T_2\in L^0(\mathcal R,\tau)$, let
\begin{equation*}
    \mathcal A=\{T_1\oplus T_2\}'\cap M_2(\mathcal R),
    \qquad
    \mathcal A_j=\{T_j\}'\cap\mathcal R.
\end{equation*}
Suppose that $xy=yx=0$ for all $x,y\in\mathcal R$ satisfying
\begin{equation*}
	T_1x=xT_2,\qquad T_2y=yT_1.
\end{equation*}
Then the diagonal inclusion $\iota\colon\mathcal A_1\oplus\mathcal A_2\to\mathcal A$ induces isomorphisms
\begin{equation*}
    V(\mathcal A_1)\oplus V(\mathcal A_2)\cong V(\mathcal A),
    \qquad
    K_0(\mathcal A_1)\oplus K_0(\mathcal A_2)\cong K_0(\mathcal A).
\end{equation*}
\end{proposition}

\begin{proof}
Let $\mathcal B=\mathcal A_1\oplus\mathcal A_2$ and
\begin{equation*}
    \mathcal J=
    \left\{
    \begin{pmatrix}0&x\\y&0\end{pmatrix}
    \colon
    x,y\in\mathcal R,\
    T_1x=xT_2,\
    T_2y=yT_1
    \right\}.
\end{equation*}
Then $\mathcal A=\mathcal B\oplus\mathcal J$ as a vector space.
By the hypothesis, $\mathcal J^2=\{0\}$.
Moreover, $\mathcal B\mathcal J\subseteq\mathcal J$ and $\mathcal J\mathcal B\subseteq\mathcal J$.
Hence $\mathcal J$ is a two-sided ideal of $\mathcal A$.
The conclusion follows from \Cref{lem square-zero-K0-splitting}.
\end{proof}

By \Cref{prop similarity-sum-powers}, $u^m+v^n$ corresponding to distinct exponent pairs are not similar in $\mathcal R$.
We now apply \Cref{prop splitting-principle} to show that the $K_0$-group of the relative commutant of their direct sum splits.

\begin{corollary}\label{cor monomial-sum-K0-split}
Let $T_j=u^{m_j}+v^{n_j}$, where $m_j,n_j\in\mathbb Z\setminus\{0\}$, for $j=1,2$.
Let
\begin{equation*}
    \mathcal A=\{T_1\oplus T_2\}'\cap M_2(\mathcal R),
    \qquad
    \mathcal A_j=\{T_j\}'\cap\mathcal R.
\end{equation*}
If $(m_1,n_1)\ne(m_2,n_2)$, then the diagonal inclusion
$\iota\colon\mathcal A_1\oplus\mathcal A_2\to\mathcal A$
induces isomorphisms
\begin{equation*}
    V(\mathcal A_1)\oplus V(\mathcal A_2)\cong V(\mathcal A),
    \qquad
    K_0(\mathcal A_1)\oplus K_0(\mathcal A_2)
    \cong K_0(\mathcal A).
\end{equation*}
\end{corollary}

\begin{proof}
Let $x,y\in\mathcal R$ such that $T_1x=xT_2$ and $T_2y=yT_1$.
We claim that $xy=yx=0$.
Otherwise, there is no loss of generality in assuming that $A=xy\ne0$. 
Clearly, $A$ commutes with $T_1$.
By \Cref{lem semigroup,thm commutant,cor SI-sum-power},
$A$ has a one-sided Fourier expansion in suitable homogeneous coordinates $U$ and $V$.
If $A_n$ is its first nonzero component, then the corresponding component of $A^2$ is $A_n^2\ne0$ by \Cref{lem commutant-degree}.
Thus, $A^2=x(yx)y\ne0$, and hence $yx\ne0$.

As in the proof of \Cref{prop similarity-sum-powers}, for $z\in\mathbb T$, let $\rho_z$ be the automorphism determined by
\begin{equation*}
    \rho_z(u)=z^{n_1}u,\qquad
    \rho_z(v)=z^{m_1}v,
\end{equation*}
and let
\begin{equation*}
    S_z=z^{-m_1n_1}\rho_z(T_2)
    =
    z^{n_1(m_2-m_1)}u^{m_2}
    +
    z^{m_1(n_2-n_1)}v^{n_2}.
\end{equation*}
Since $\rho_z(T_1)=z^{m_1n_1}T_1$, the operators $x_z=y\rho_z(x)$ and $y_z=\rho_z(y)x$ satisfy
\begin{equation*}
    T_2x_z=x_zS_z,\qquad S_zy_z=y_zT_2.
\end{equation*}
Since $x_1=y_1=yx\ne0$, we have $x_z,y_z\ne0$ for every $z$ in some neighborhood of $1$ by the $L^2$-continuity.
By \Cref{lem mutual-phase-changes-sum-powers}, we have
\begin{equation*}
    z^{n_1(m_2-m_1)}\in\alpha^{m_2\mathbb Z},
    \qquad
    z^{m_1(n_2-n_1)}\in\alpha^{n_2\mathbb Z}
\end{equation*}
throughout that neighborhood.
Countability of these sets forces both exponents to be zero.
Since $m_1n_1\ne0$, this implies
$(m_1,n_1)=(m_2,n_2)$, a contradiction.

Therefore, $xy=0$.
Interchanging the two summands also gives $yx=0$.
The conclusion follows from \Cref{prop splitting-principle}.
\end{proof}

\subsection{Homogeneous operators of different degrees}

In this subsection, we show that the $K_0$-group of the relative commutant of a direct sum splits when the summands are injective homogeneous operators of different degrees.
We first introduce the scalar intertwining set.

\begin{definition}\label{def Lambda-T}
For every measurable operator $T\in L^0(\mathcal R,\tau)$, its \emph{scalar intertwining set} $\Lambda(T)$ is defined as $\Lambda(T)=\{\lambda\in\mathbb C\colon\lambda Tx=xT~\text{for some nonzero}~x\in\mathcal R\}$.
\end{definition}

In the following lemma, we show that the scalar intertwining set of an injective homogeneous operator in $L^0(\mathcal R,\tau)$ is always countable.

\begin{lemma}\label{lem countable-twisted-commutation}
Let $T=u^mf(v)\in L^0(\mathcal R,\tau)$, where $m\in\mathbb Z$ and $f\ne0$ almost everywhere.
Then $\Lambda(T)$ is a countable subset of $\mathbb C\setminus\{0\}$.
Moreover, $\Lambda(T)$ is a subset of $\{\alpha^n\colon n\in\mathbb Z\}$ whenever $m\ne 0$.
\end{lemma}

\begin{proof}
Since $T$ is invertible in $L^0(\mathcal R,\tau)$, we have $0\notin\Lambda(T)$.
Let $\lambda\in\Lambda(T)$ and choose a nonzero $x\in\mathcal R$ such that $\lambda Tx=xT$.
By \Cref{lem frequently}, we may assume that $x=u^ng(v)$ is homogeneous.
Then
\begin{equation}\label{equ lambda-fg}
    \lambda f(\alpha^nz)g(z)
    =
    g(\alpha^mz)f(z)
    \qquad\text{almost everywhere}.
\end{equation}
If $m=0$, then $\bigl(\lambda f(\alpha^nz)-f(z)\bigr)g(z)=0$.
For each $n\in\mathbb Z$ and $\lambda\ne0$, let
\begin{equation*}
    \Omega_{n,\lambda}
    =
    \left\{z\in\mathbb T\colon
    \frac{f(z)}{f(\alpha^nz)}=\lambda\right\}.
\end{equation*}
It follows from \eqref{equ lambda-fg} that $\mu(\Omega_{n,\lambda})\geqslant\mu(\{z\in\mathbb T\colon g(z)\ne0\})>0$.
For each fixed $n$, the sets $\Omega_{n,\lambda}$ are pairwise disjoint as $\lambda$ varies.
Since $\mu(\mathbb T)=1$, only countably many of them can have positive measure.
Taking the union over $n\in\mathbb Z$ proves countability when $m=0$.

Suppose now that $m\ne0$.
Since $f\ne0$ almost everywhere and $\lambda\ne0$, it follows from \eqref{equ lambda-fg} that the zero set $Z(g)$ of $g$ satisfies $Z(g)=\alpha^mZ(g)$ modulo a null set.
Since $g$ is not identically zero, $g\ne0$ almost everywhere.
Thus, both $T$ and $x$ are invertible in $L^0(\mathcal R,\tau)$.
Let
\begin{equation*}
    d=\gcd(|m|,|n|)\geqslant1,\qquad m=m_1d,\qquad n=n_1d.
\end{equation*}
Since $x^{m_1}T^{-n_1}$ has degree zero, we can write $h(v)=x^{m_1}T^{-n_1}$, where $h\ne 0$ almost everywhere.
Since $\lambda Tx=xT$, we have $\lambda^{m_1}Th(v)=h(v)T$ and $\lambda^{n_1}xh(v)=h(v)x$.
It follows that
\begin{equation}\label{equ h-2-identity}
    h(\alpha^mz)=\lambda^{m_1}h(z),\qquad
    h(\alpha^nz)=\lambda^{n_1}h(z)
    \qquad\text{almost everywhere}.
\end{equation}
By the first identity of \eqref{equ h-2-identity}, we see that $|\lambda|=1$ and $|h|$ is essentially constant.
In particular, $h\in L^\infty(\mathbb T,\mu)$.
Comparing Fourier coefficients in the first identity of \eqref{equ h-2-identity}, we obtain that
\begin{equation*}
    h(z)=cz^j,\qquad \lambda^{m_1}=\alpha^{jm}
\end{equation*}
for some $j\in\mathbb Z$ and $c\ne0$.
The second identity of \eqref{equ h-2-identity} gives $\lambda^{n_1}=\alpha^{jn}$.
Since $\gcd(|m_1|,|n_1|)=1$, there are $p,q\in\mathbb Z$ such that $pm_1+qn_1=1$.
Then
\begin{equation*}
    \lambda
    =(\lambda^{m_1})^p(\lambda^{n_1})^q
    =\alpha^{j(pm+qn)}
    =\alpha^{jd}.
\end{equation*}
Therefore, $\Lambda(T)\subseteq\{\alpha^n\colon n\in\mathbb Z\}$, which completes the proof.
\end{proof}

The preceding countability result implies that products of intertwiners in opposite directions vanish when the degrees are different.
The following lemma complements \Cref{lem 2-intertwining}.

\begin{lemma}\label{lem different-degree-zero-products}
Let $T_j=u^{m_j}f_j(v)\in L^0(\mathcal R,\tau)$, where
$m_j\in\mathbb Z$ and $f_j\ne0$ almost everywhere, for $j=1,2$.
Suppose that $m_1\ne m_2$.
If $x,y\in\mathcal R$ satisfy
\begin{equation*}
    T_1x=xT_2,\qquad T_2y=yT_1,
\end{equation*}
then $xy=yx=0$.
\end{lemma}

\begin{proof}
Let $d=m_1-m_2\ne0$.
For every $z\in\mathbb T$, let $\pi_z$ be the automorphism of $\mathcal R$ such that $\pi_z(u)=zu$ and $\pi_z(v)=v$ as defined in \eqref{equ pi-z}.
Since $\pi_z$ is trace-preserving, it extends to a $*$-automorphism of $L^0(\mathcal R,\tau)$, which we also denote by $\pi_z$.
Applying $\pi_z$ to $T_1x=xT_2$, we obtain that $z^dT_1\pi_z(x)=\pi_z(x)T_2$.
It follows that
\begin{equation*}
    z^dT_1\pi_z(x)y
    =\pi_z(x)T_2y
    =\pi_z(x)yT_1.
\end{equation*}
Thus, $z^d\in\Lambda(T_1)$ whenever $\pi_z(x)y\ne0$.
Since $d\ne0$ and $\Lambda(T_1)$ is countable by \Cref{lem countable-twisted-commutation}, we have $\pi_z(x)y=0$ for all but countably many $z\in\mathbb T$.
Since the map $z\mapsto\pi_z(x)$ is continuous in $L^2(\mathcal R,\tau)$, we have $\pi_z(x)y=0$ for every $z\in\mathbb T$.
In particular, $xy=0$.
By a similar argument, we also have $yx=0$.
\end{proof}

Now we prove the splitting result for different degrees.

\begin{proposition}\label{prop different-degree-K0-split}
Let $T_j=u^{m_j}f_j(v)\in L^0(\mathcal R,\tau)$, where $m_j\in\mathbb Z$ and $f_j\ne0$ almost everywhere, for $j=1,2$.
Let
\begin{equation*}
    \mathcal A=\{T_1\oplus T_2\}'\cap M_2(\mathcal R),\qquad
    \mathcal A_j=\{T_j\}'\cap\mathcal R.
\end{equation*}
If $m_1\ne m_2$, then the diagonal inclusion $\iota\colon\mathcal A_1\oplus\mathcal A_2\to\mathcal A$ induces isomorphisms
\begin{equation*}
    V(\mathcal A_1)\oplus V(\mathcal A_2)\cong V(\mathcal A),
    \qquad
    K_0(\mathcal A_1)\oplus K_0(\mathcal A_2)\cong K_0(\mathcal A).
\end{equation*}
\end{proposition}

\begin{proof}
By \Cref{lem different-degree-zero-products}, we have
$xy=yx=0$ whenever
\begin{equation*}
    T_1x=xT_2,\qquad
    T_2y=yT_1.
\end{equation*}
The conclusion follows from \Cref{prop splitting-principle}.
\end{proof}

\subsection{Homogeneous operators of the same degree}

In this subsection, we consider injective homogeneous operators of the same nonzero degree.
For every $T\in L^0(\mathcal R,\tau)$, its \emph{relative invariant subspace lattice} is defined by
\begin{equation*}
    \operatorname{Lat}_{\mathcal R}(T)
    =
    \{p\in\operatorname{Proj}(\mathcal R)\colon
    (1-p)Tp=0\},
\end{equation*}
with the order inherited from $\operatorname{Proj}(\mathcal R)$.
Denote by $R(T)$ the \emph{range projection} of $T$.

\begin{lemma}\label{lem lattice-intertwining-L0}
Let $T,S,H\in L^0(\mathcal R,\tau)$, where $TH=HS$ and $H$ is injective.
Then the map
\begin{equation*}
    \varphi\colon
    \operatorname{Lat}_{\mathcal R}(S)
    \to\operatorname{Lat}_{\mathcal R}(T),
    \qquad
    \varphi(e)=R(He),
\end{equation*}
is a trace-preserving order isomorphism, with inverse $\psi(p)=R(H^{-1}p)$.
\end{lemma}

\begin{proof}
Let $e\in\operatorname{Lat}_{\mathcal R}(S)$, $p=R(He)$, and $q=R(H^{-1}p)$.
The right supports of $He$ and $H^{-1}p$ are $e$ and $p$, respectively.
It follows that
\begin{equation*}
    \tau(e)=\tau(p)=\tau(q).
\end{equation*}
Since $pHe=He$ and $qH^{-1}p=H^{-1}p$, we obtain
\begin{equation*}
    qe=qH^{-1}pHe=H^{-1}pHe=e.
\end{equation*}
Thus, $e\leqslant q$, and equality of traces implies that $q=e$.
In particular, we have
\begin{equation*}
    H^{-1}p=eH^{-1}p.
\end{equation*}
Using $T=HSH^{-1}$, $Se=eSe$, and $(1-p)He=0$, we obtain that
\begin{equation*}
    (1-p)Tp
    =(1-p)HSH^{-1}p
    =(1-p)HeSeH^{-1}p
    =0.
\end{equation*}
Therefore, $\varphi(e)\in\operatorname{Lat}_{\mathcal R}(T)$, $\tau(\varphi(e))=\tau(e)$, and $(\psi\varphi)(e)=e$.
The same argument shows that $\psi(p)\in\operatorname{Lat}_{\mathcal R}(S)$, $\tau(\psi(p))=\tau(p)$, and $(\varphi\psi)(p)=p$.
It is clear that both maps preserve order.
This completes the proof.
\end{proof}

The following theorem gives several conditions under which the $K_0$-group of the relative commutant of $T_1\oplus T_2$ splits.

\begin{theorem}\label{thm K0-split}
Let $T_j=u^mf_j(v)\in\mathcal R$, where $m\ne 0$ and $f_j\ne0$ almost everywhere, for $j=1,2$.
Let
\begin{equation*}
    \mathcal A=\{T_1\oplus T_2\}'\cap M_2(\mathcal R),\qquad
    \mathcal A_j=\{T_j\}'\cap\mathcal R.
\end{equation*}
Then the diagonal inclusion $\iota\colon\mathcal A_1\oplus\mathcal A_2\to\mathcal A$ induces isomorphisms
\begin{equation*}
    V(\mathcal A_1)\oplus V(\mathcal A_2)\cong V(\mathcal A),
    \qquad
    K_0(\mathcal A_1)\oplus K_0(\mathcal A_2)\cong K_0(\mathcal A),
\end{equation*}
whenever one of the following statements holds:
\begin{enumerate}
\item[$(1)$] $\Lambda(T_1)\ne\Lambda(T_2)$;
\item[$(2)$] $\sigma(T_1)\ne\sigma(T_2)$;
\item[$(3)$] $\Delta(T_1)\ne\Delta(T_2)$;
\item[$(4)$] $\Gamma(T_1)\ne\Gamma(T_2)$;
\item[$(5)$] $\gcd\Gamma_b(T_1)\ne\gcd\Gamma_b(T_2)$;
\item[$(6)$] $\operatorname{Lat}_{\mathcal R}(T_1) \not\cong\operatorname{Lat}_{\mathcal R}(T_2)$;
\item[$(7)$] exactly one of $T_1,T_2$ is power bounded;
\item[$(8)$] $T_1$ and $T_2$ are not similar in $\mathcal R$, and at least one of them is invertible.
\end{enumerate}
\end{theorem}

\begin{proof}
We assume that $m\geqslant1$, since the case $m\leqslant-1$ is analogous.
Let
\begin{equation*}
    \mathcal I_{12}
    =
    \{x\in\mathcal R\colon T_1x=xT_2\},
    \qquad
    \mathcal I_{21}
    =
    \{y\in\mathcal R\colon T_2y=yT_1\}.
\end{equation*}
We show that each of the stated conditions forces either $\mathcal I_{12}=0$ or $\mathcal I_{21}=0$.

Suppose that both spaces contain nonzero operators.
By \Cref{lem mutual-homogeneous-intertwining}, we can choose nonzero homogeneous operators
\begin{equation*}
    x=u^ng(v)\in\mathcal I_{12},\qquad
    y=u^kh(v)\in\mathcal I_{21}
\end{equation*}
and an integer $a\geqslant1$ such that
\begin{equation}\label{equ mutual-intertwining-powers}
    (xy)^m=T_1^a,\qquad (yx)^m=T_2^a.
\end{equation}
Comparing degrees gives $n+k=a$.
Moreover, $x,y$ are invertible in $L^0(\mathcal R,\tau)$,
as in the proof of that lemma.
Let $A=xy$ and $B=yx$.
We now show that none of $(1)$--$(8)$ can hold.

$(1)$
Let $\lambda\in\Lambda(T_1)$ and choose a nonzero $w\in\mathcal R$ such that $\lambda T_1w=wT_1$.
Then
\begin{equation*}
    \lambda T_2(ywx)
    =
    y(\lambda T_1w)x
    =
    ywT_1x
    =
    (ywx)T_2.
\end{equation*}
Since $x$ and $y$ are invertible in $L^0(\mathcal R,\tau)$, we have $ywx\ne0$.
Thus, $\Lambda(T_1)\subseteq\Lambda(T_2)$.
Interchanging $T_1$ and $T_2$ gives
$\Lambda(T_1)=\Lambda(T_2)$.

$(2)$
The nonzero spectra of $xy$ and $yx$ agree.
Since $\mathcal R$ is finite, $xy$ is invertible if and only
if $yx$ is invertible.
Thus, $\sigma(A)=\sigma(B)$.
Combining this with \eqref{equ mutual-intertwining-powers} and spectral mapping, we obtain that
\begin{equation*}
    \sigma(T_1^a)=\sigma(A^m)=\sigma(B^m)=\sigma(T_2^a).
\end{equation*}
Thus, $\sigma(T_1)=\sigma(T_2)$ by \Cref{cor spectrum}.

$(3)$
By multiplicativity of the Fuglede--Kadison determinant, we have
\begin{equation*}
    \Delta(A)=\Delta(x)\Delta(y)=\Delta(B).
\end{equation*}
Taking determinants in \eqref{equ mutual-intertwining-powers}
gives
\begin{equation*}
    \Delta(T_1)^a
    =
    \Delta(A)^m
    =
    \Delta(B)^m
    =
    \Delta(T_2)^a.
\end{equation*}
Since $a\geqslant1$, it follows that $\Delta(T_1)=\Delta(T_2)$.

$(4)$
Since $\mathcal C_j(T_2)=x^{-1}\mathcal C_j(T_1)x$ for every $j\in\mathbb Z$, we have $\Gamma(T_1)=\Gamma(T_2)$.

$(5)$
Since $A\in\mathcal C_a(T_1)\cap\mathcal R$ and $B\in\mathcal C_a(T_2)\cap\mathcal R$, we have $a\in\Gamma_b(T_1)\cap\Gamma_b(T_2)$.
For every nonzero $w\in\mathcal C_j(T_1)\cap\mathcal R$, the operator $ywx$ is nonzero, bounded, has degree $j+a$, and commutes with $T_2$.
Interchanging the two operators, we therefore obtain
\begin{equation*}
    \Gamma_b(T_1)+a\subseteq\Gamma_b(T_2),
    \qquad
    \Gamma_b(T_2)+a\subseteq\Gamma_b(T_1).
\end{equation*}
Let $d_j=\gcd\Gamma_b(T_j)$ for $j=1,2$.
For every $j\in\Gamma_b(T_1)$, the integer $d_2$ divides both $a$ and $j+a$, and hence divides $j$.
Thus, $d_2\mid d_1$.
The reverse divisibility follows in the same way.
Hence $\gcd\Gamma_b(T_1)=\gcd\Gamma_b(T_2)$.

$(6)$
Since $x$ is injective and $T_1x=xT_2$, we have $\operatorname{Lat}_{\mathcal R}(T_1)\cong\operatorname{Lat}_{\mathcal R}(T_2)$ by \Cref{lem lattice-intertwining-L0}.

$(7)$
For every integer $N\geqslant0$, we have
\begin{equation*}
    T_1^{N+a}
    =
    T_1^NA^m
    =
    xT_2^NB^{m-1}y,
    \qquad
    T_2^{N+a}
    =
    yT_1^NA^{m-1}x.
\end{equation*}
Therefore, there exists $C>0$, independent of $N$, such that
\begin{equation*}
    \|T_1^{N+a}\|\leqslant C\|T_2^N\|,
    \qquad
    \|T_2^{N+a}\|\leqslant C\|T_1^N\|.
\end{equation*}
It follows that $T_1$ is power bounded if and only if $T_2$ is power bounded.

$(8)$
If $T_1$ is invertible in $\mathcal R$, then $A^m=T_1^a$ is invertible, and hence $A=xy$ is invertible.
Since $\mathcal R$ is finite, both $x$ and $y$ are invertible in $\mathcal R$.
Thus, $T_1=xT_2x^{-1}$.
The same argument applies if $T_2$ is invertible.
Consequently, if at least one of $T_1,T_2$ is invertible, then they are similar in $\mathcal R$.

We have shown that the existence of nonzero operators in both $\mathcal I_{12}$ and $\mathcal I_{21}$ excludes every stated condition.
Thus, under any of $(1)$--$(8)$, at least one of these intertwiner spaces is zero.
The conclusion follows from \Cref{prop splitting-principle}.
\end{proof}

The splitting results also yield a similarity criterion as follows.

\begin{corollary}\label{cor K0-similarity-criterion}
Let $T_j=u^{m_j}f_j(v)\in\mathcal R$, where $m_j\ne0$ and $f_j\ne0$ almost everywhere, for $j=1,2$, and let
\begin{equation*}
    \mathcal A=\{T_1\oplus T_2\}'\cap M_2(\mathcal R).
\end{equation*}
If $K_0(\mathcal A)\cong\mathbb Z$, then $T_1$ and $T_2$ are similar in $\mathcal R$.
\end{corollary}

\begin{proof}
Suppose, towards a contradiction, that $T_1$ and $T_2$
are not similar in $\mathcal R$, and let
\begin{equation*}
    \mathcal A_j=\{T_j\}'\cap\mathcal R,\qquad j=1,2.
\end{equation*}
Since each $K_0(\mathcal A_j)$ contains a direct summand isomorphic to $\mathbb Z$ by \Cref{lem commutant-K0-Z-summand}, a splitting
\begin{equation*}
    K_0(\mathcal A)
    \cong K_0(\mathcal A_1)\oplus K_0(\mathcal A_2)
\end{equation*}
would give a direct summand isomorphic to $\mathbb Z^2$
in $K_0(\mathcal A)$, contradicting
$K_0(\mathcal A)\cong\mathbb Z$.

If $m_1\ne m_2$, such a splitting follows from
\Cref{prop different-degree-K0-split}.
Thus, $m_1=m_2$.
If at least one of $T_1,T_2$ is invertible, the same
splitting follows from condition $(8)$ of
\Cref{thm K0-split}.
Therefore, both $T_1$ and $T_2$ are noninvertible.

In this case, $\tau$ is multiplicative on both $\mathcal A_1$ and $\mathcal A_2$.
Moreover, \Cref{lem 2-intertwining} and nonsimilarity imply
that
\begin{equation*}
    \tau(xy)=\tau(yx)=0
\end{equation*}
whenever $x,y\in\mathcal R$ satisfy
$T_1x=xT_2$ and $T_2y=yT_1$.
Hence the map
\begin{equation*}
    \rho\colon\mathcal A\to\mathbb C\oplus\mathbb C,
    \qquad
    \begin{pmatrix}a&x\\y&b\end{pmatrix}
    \mapsto(\tau(a),\tau(b))
\end{equation*}
is a unital algebra homomorphism.
Let $\iota\colon\mathbb C^2\to\mathcal A$ be the diagonal scalar inclusion.
Since $\rho\iota=\operatorname{id}_{\mathbb C^2}$, the induced homomorphisms satisfy $\rho_*\iota_*=\operatorname{id}_{\mathbb Z^2}$.
Thus, $\mathbb Z^2$ is a direct summand of $K_0(\mathcal A)$, again contradicting $K_0(\mathcal A)\cong\mathbb Z$.
Therefore, $T_1$ and $T_2$ are similar in $\mathcal R$.
\end{proof}

\subsection{Rotation pairs with different parameters}

The preceding results concern homogeneous operators associated with a fixed irrational rotation pair.
We now consider operators $T_j=u_j(1+v_j)$ associated with possibly different rotation parameters.
We first record a consequence of the Fourier expansion of the relative commutant of $u(1+v)$.

\begin{lemma}\label{lem rotation-sum-nonzero-products}
Let $T=u(1+v)\in\mathcal R$.
Then every nonzero operator $D\in\{T\}'\cap\mathcal R$ satisfies $D^2\ne 0$ and $DvD\ne 0$.
\end{lemma}

\begin{proof}
By \Cref{thm commutant}, $D$ has a unique $L^2$-Fourier expansion
\begin{equation*}
    D=\sum_{n=0}^{\infty}a_nT^n.
\end{equation*}
Let $n\geqslant0$ be the smallest integer such that $a_n\ne0$.
Then
\begin{equation*}
    \Pi_{2n}(D^2)=a_n^2T^{2n},
    \qquad
    \Pi_{2n}(DvD)=a_n^2T^nvT^n.
\end{equation*}
Since $T$ is invertible in $L^0(\mathcal R,\tau)$, we have $D^2\ne0$ and $DvD\ne0$.
\end{proof}

Using the preceding lemma and \Cref{lem countable-twisted-commutation}, we show that a nonzero product of intertwiners in opposite directions forces $\theta_2\equiv\pm\theta_1\pmod{\mathbb Z}$.

\begin{lemma}\label{lem rotation-sum-intertwining-products}
For $j=1,2$, let $u_j,v_j\in\mathcal R$ be unitary operators satisfying
\begin{equation*}
    \mathcal R=W^*(u_j,v_j),\qquad
    v_ju_j=\alpha_ju_jv_j,\qquad
    \alpha_j=e^{2\pi i\theta_j},\qquad
    \theta_j\in\mathbb R\setminus\mathbb Q.
\end{equation*}
Let $T_j=u_j(1+v_j)$.
If $x,y\in\mathcal R$ satisfy
\begin{equation*}
    T_1x=xT_2,\qquad T_2y=yT_1,\qquad xy\ne0,
\end{equation*}
then $\theta_2\equiv\pm\theta_1\pmod{\mathbb Z}$.
\end{lemma}

\begin{proof}
Let $A_1=xy\ne 0$ and $A_2=yx$.
Then $A_j$ commutes with $T_j$ for $j=1,2$.
By \Cref{lem rotation-sum-nonzero-products}, we have $xA_2y=A_1^2\ne0$, and hence $A_2\ne 0$.
By \Cref{lem rotation-sum-nonzero-products} again, we obtain that
\begin{equation*}
	x(yv_1x)y=A_1v_1A_1\ne 0,\qquad y(xv_2y)x=A_2v_2A_2\ne 0.
\end{equation*}
Thus, $B_1=xv_2y$ and $B_2=yv_1x$ are nonzero.
Since $T_jv_j=\alpha_j^{-1}v_jT_j$, we have
\begin{equation*}
    \alpha_2T_1B_1=B_1T_1,
    \qquad
    \alpha_1T_2B_2=B_2T_2.
\end{equation*}
It follows from \Cref{lem countable-twisted-commutation} that
\begin{equation*}
    \alpha_2\in\Lambda(T_1)\subseteq\alpha_1^{\mathbb Z},
    \qquad
    \alpha_1\in\Lambda(T_2)\subseteq\alpha_2^{\mathbb Z}.
\end{equation*}
Write $\alpha_2=\alpha_1^p$ and $\alpha_1=\alpha_2^q$ for some integers $p,q$.
Then $\alpha_1^{pq-1}=1$ and hence $pq=1$.
Consequently, $p=q=1$ or $p=q=-1$.
Therefore, $\theta_2\equiv\pm\theta_1\pmod{\mathbb Z}$.
\end{proof}

When the rotation parameters satisfy $\theta_2\not\equiv\pm\theta_1\pmod{\mathbb Z}$, the preceding lemma shows that all products of intertwiners in opposite directions vanish.
The splitting principle therefore applies.

\begin{proposition}\label{prop rotation-sums-K0-split}
For $j=1,2$, let $u_j,v_j\in\mathcal R$ be unitary operators satisfying
\begin{equation*}
    \mathcal R=W^*(u_j,v_j),\qquad
    v_ju_j=\alpha_ju_jv_j,\qquad
    \alpha_j=e^{2\pi i\theta_j},\qquad
    \theta_j\in\mathbb R\setminus\mathbb Q.
\end{equation*}
Let
\begin{equation*}
    T_j=u_j(1+v_j),\qquad
    \mathcal A_j=\{T_j\}'\cap\mathcal R,\qquad
    \mathcal A=\{T_1\oplus T_2\}'\cap M_2(\mathcal R).
\end{equation*}
If $\theta_2\not\equiv\pm\theta_1\pmod{\mathbb Z}$, then
the diagonal inclusion
$\iota\colon\mathcal A_1\oplus\mathcal A_2\to\mathcal A$
induces isomorphisms
\begin{equation*}
    V(\mathcal A_1)\oplus V(\mathcal A_2)\cong V(\mathcal A),
    \qquad
    K_0(\mathcal A_1)\oplus K_0(\mathcal A_2)
    \cong K_0(\mathcal A).
\end{equation*}
\end{proposition}

\begin{proof}
By \Cref{lem rotation-sum-intertwining-products}, $xy=yx=0$ for all $x,y\in\mathcal R$ satisfying
\begin{equation*}
    T_1x=xT_2,\qquad T_2y=yT_1.
\end{equation*}
The asserted isomorphisms follow from \Cref{prop splitting-principle}.
\end{proof}

\section{Strong generators and strongly irreducible indices}\label{sec SG}

In this section, we construct homogeneous operators of degree one in $\mathcal R$ such that every operator similar to any one of them generates $\mathcal R$ as a von Neumann algebra.
The operators we construct also remain strongly irreducible under every irreducible inclusion of $\mathcal R$ into a type $\mathrm{II}_1$ factor.
As applications, we obtain strongly irreducible operators in every separable type $\mathrm{II}_1$ factor and show that the strongly irreducible indices of $\mathcal R$ coincide with its irreducible subfactor indices.

\subsection{A criterion for strong generation}

An operator $T$ in $\mathcal R$ is strongly irreducible relative to $\mathcal R$ if and only if $XTX^{-1}$ is irreducible relative to $\mathcal R$ for every invertible operator $X$ in $\mathcal R$.
Recall that an operator $T$ in $\mathcal R$ is called a \emph{generator} if $W^*(T)=\mathcal R$.
This motivates the following definition.

\begin{definition}\label{def SG}
Let $T$ be an operator in $\mathcal R$.
We say that $T$ is a \emph{strong generator} of $\mathcal R$ if $XTX^{-1}$ is a generator of $\mathcal R$ for every invertible operator $X$ in $\mathcal R$.
\end{definition}

In this subsection, we establish a sufficient condition for an operator $T=uf(v)$ to be a strong generator of $\mathcal R$.
The proof begins with a Fourier-analytic lemma.
Recall that a real-valued function $\psi$ on $\mathbb R$ is \emph{sub-additive} if $\psi(s+t)\leqslant\psi(s)+\psi(t)$ for all $s,t\in\mathbb R$.

\begin{lemma}\label{lem subadditive-Fourier}
Suppose that $\{a_n\}_{n=1}^\infty$ is a nonnegative nondecreasing sequence and $\{q_n\}_{n=1}^\infty$ is a nonnegative sequence such that
\begin{equation*}
	\sum_{n=1}^\infty\frac{a_n}{n^2}=\infty,\qquad
	\sum_{n=1}^\infty a_nq_n<\infty,\qquad
	\sum_{n=1}^\infty q_n<\infty.
\end{equation*}
Let
\begin{equation*}
	\psi(t)=\sum_{n=1}^\infty
	q_n\bigl(1-\cos(2\pi nt)\bigr),
	\qquad t\in\mathbb R.
\end{equation*}
If $\psi$ is sub-additive, then $\psi=0$.
\end{lemma}

\begin{proof}
For every $N\geqslant 1$, let
\begin{equation*}
	\psi_N(t)=\sum_{n=1}^N
	q_n\bigl(1-\cos(2\pi nt)\bigr),
	\qquad
	Q_N=\sum_{n=N+1}^\infty q_n.
\end{equation*}
Then $\psi_N(0)=\psi_N'(0)=0$ and $0\leqslant\psi-\psi_N\leqslant2Q_N$.
By Bernstein's inequality for trigonometric polynomials \cite[Theorem~1.1]{QZ-19}, we have
\begin{equation*}
	\|\psi_N''\|_\infty\leqslant2\pi N\|\psi_N'\|_\infty.
\end{equation*}
Since $\psi_N'$ is an odd $1$-periodic function, there exists $s\in\mathbb R$ such that
$\psi_N'(s)=\|\psi_N'\|_\infty$.
By Taylor's mean value theorem, for every $t>0$, we have
\begin{equation*}
	\psi_N(s+t)-\psi_N(s)
	\geqslant \|\psi_N'\|_\infty t-\pi N\|\psi_N'\|_\infty t^2,
	\qquad
	\psi_N(t)\leqslant\pi N\|\psi_N'\|_\infty t^2.
\end{equation*}
Since $\psi$ is sub-additive and $0\leqslant\psi-\psi_N\leqslant 2Q_N$, we obtain
\begin{equation*}
	\psi_N(s+t)-\psi_N(s)
	\leqslant\psi(s+t)-\psi(s)+2Q_N
	\leqslant\psi(t)+2Q_N
	\leqslant \psi_N(t)+4Q_N.
\end{equation*}
It follows that $\|\psi_N'\|_\infty t-\pi N\|\psi_N'\|_\infty t^2\leqslant\pi N\|\psi_N'\|_\infty t^2+4Q_N$.
Taking $t=(4\pi N)^{-1}$, we have
\begin{equation*}
	\|\psi_N'\|_\infty\leqslant32\pi NQ_N.
\end{equation*}
Suppose to the contrary that $\psi\ne0$.
Since $\psi$ is nonnegative and $1$-periodic with $\psi(0)=0$, there exists $t_0\in(0,1)$ such that $\psi(t_0)>0$.
Since $\psi_N\to\psi$ uniformly, there exists $N_0\geqslant 1$ such that
\begin{equation*}
	\frac{\psi(t_0)}2
	\leqslant\psi_N(t_0)
	=\int_0^{t_0}\psi_N'(t)\,dt
	\leqslant\|\psi_N'\|_\infty
	\leqslant32\pi NQ_N,\qquad N\geqslant N_0.
\end{equation*}
It follows that $Q_N\geqslant\frac{c}{N}$ for $N\geqslant N_0$, where $c=\frac{\psi(t_0)}{64\pi}>0$.
By Abel's summation formula, we have
\begin{equation*}
	\sum_{n=1}^\infty a_nq_n=a_1\sum_{n=1}^\infty q_n
	+\sum_{n=1}^\infty(a_{n+1}-a_n)Q_n.
\end{equation*}
It follows that
\begin{equation*}
	\sum_{n=N_0}^\infty\frac{a_{n+1}-a_n}{n}
	\leqslant\frac1c
	\sum_{n=N_0}^\infty(a_{n+1}-a_n)Q_n
	\leqslant\frac1c\sum_{n=1}^\infty a_nq_n
	<\infty.
\end{equation*}
By Abel's summation formula again, we have
\begin{equation*}
	\sum_{n=1}^\infty\frac{a_n}{n(n+1)}
	=a_1+\sum_{k=1}^\infty
	\frac{a_{k+1}-a_k}{k+1}
	<\infty.
\end{equation*}
This contradicts $\sum_{n=1}^\infty\frac{a_n}{n^2}=\infty$.
Therefore, $\psi=0$.
\end{proof}

We now apply the preceding lemma to obtain a sufficient condition for $\mathcal R$ to be generated by $u$ and a positive operator $H$ in $\mathcal R$ comparable to $h(v)^2$, where $h$ is a nonnegative function in $L^\infty(\mathbb T,\mu)$.
The operator $H$ need not belong to $W^*(v)$.

\begin{lemma}\label{lem SG-rigid}
For a nonnegative function $h\in L^\infty(\mathbb T,\mu)$, $r>0$, and $z\in\mathbb T$, let
\begin{equation*}
	\Psi_r(z)=\int_{-r}^r
	h^2(ze^{2\pi it})\,dt,
	\qquad
	\Lambda=\{\alpha^n\colon n\in\mathbb Z\}.
\end{equation*}
Suppose that $\Psi_r(1)>0$ for every $r>0$ and
\begin{equation}\label{equ SG-rigid-1}
	\lim_{r\to0}\frac{\Psi_r(z)+\Psi_r(z^{-1})}{2\Psi_r(1)}=\infty
	\qquad(z\in\mathbb T\setminus\Lambda).
\end{equation}
Suppose also that there exists a nonnegative nondecreasing sequence
$\{a_n\}_{n=1}^\infty$ such that
\begin{equation}\label{equ SG-rigid-2}
	\sum_{n=1}^\infty\frac{a_n}{n^2}=\infty,
	\qquad
	\liminf_{r\to0}\frac{\Psi_r(\alpha^n)+\Psi_r(\alpha^{-n})}{2\Psi_r(1)}
	\geqslant a_{|n|}
	\qquad(n\in\mathbb Z\setminus\{0\}).
\end{equation}
Then $W^*(u,H)=\mathcal R$ for every positive operator $H\in\mathcal R$
satisfying
\begin{equation*}
	C_1h(v)^2\leqslant H\leqslant C_2h(v)^2
\end{equation*}
for some $C_1,C_2>0$.
\end{lemma}

\begin{proof}
Let $\mathcal U=W^*(u)$, $\mathcal V=W^*(v)$, and $\mathcal P=W^*(u,H)$.
Since $u\in\mathcal P$, we have
\begin{equation*}
	E_{\mathcal P}(u^*xu)=u^*E_{\mathcal P}(x)u
	\qquad(x\in\mathcal R).
\end{equation*}
It follows from $u^*E_{\mathcal P}(v^n)u=E_{\mathcal P}(u^*v^nu)=\alpha^nE_{\mathcal P}(v^n)$ that 
\begin{equation*}
	E_{\mathcal P}(v^n)=e_nv^n\qquad\text{for some}~e_n\in\mathcal U.
\end{equation*}
Since $e_n^2v^n=e_nE_{\mathcal P}(v^n)=E_{\mathcal P}(e_nv^n)=e_nv^n$, we have $e_n^2=e_n$, and hence $e_n$ is a projection in the abelian von Neumann algebra $\mathcal U$.
Let
\begin{equation*}
	\gamma_n=\tau(e_n)
	=\tau(v^{-n}E_{\mathcal P}(v^n))
	=\|E_{\mathcal P}(v^n)\|_2^2,\qquad (n\in\mathbb Z).
\end{equation*}
Then $\gamma_0=1$ and $\gamma_{-n}=\gamma_n\in[0,1]$.
Since both $e_m$ and $v^me_nv^{*m}$ are projections in $\mathcal U$, their product $q=e_mv^me_nv^{*m}$ is also a projection in $\mathcal U$.
Moreover,
\begin{equation*}
	qv^{m+n}=e_mv^me_nv^n=E_{\mathcal P}(v^m)E_{\mathcal P}(v^n)\in\mathcal P.
\end{equation*}
Applying $E_{\mathcal P}$ gives $qv^{m+n}=E_{\mathcal P}(qv^{m+n})=qe_{m+n}v^{m+n}$.
Hence $q=qe_{m+n}\leqslant e_{m+n}$.
It follows that
\begin{equation}\label{equ SG-defect-subadditive}
	\gamma_{m+n}
	=\tau(e_{m+n})
	\geqslant\tau(q)
	\geqslant\tau(e_m)+\tau(v^me_nv^{*m})-1
	=\gamma_m+\gamma_n-1.
\end{equation}
Let $K=E_{\mathcal V}E_{\mathcal P}|_{\mathcal V}$ and identify $\mathcal V$ with $L^\infty(\mathbb T,\mu)$.
Since $E_{\mathcal V}(e_nv^n)=\tau(e_n)v^n$, we have $K(z^n)=\gamma_nz^n$.
For $w\in\mathbb T$, we define rotation $(R_wg)(z)=g(wz)$ and reflection $(Jg)(z)=g(z^{-1})$.
Then
\begin{equation*}
	KR_w(z^n)=w^n\gamma_nz^n=R_wK(z^n),
	\qquad
	KJ(z^n)=\gamma_{-n}z^{-n}
	=\gamma_nz^{-n}=JK(z^n).
\end{equation*}
Thus, $KR_w=R_wK$ and $KJ=JK$ hold on $L^2(\mathbb T,\mu)$.
Moreover, $K$ preserves $C(\mathbb T)$, since it preserves
trigonometric polynomials and is contractive in the uniform norm.
By the assumption $C_1h(v)^2\leqslant H\leqslant C_2h(v)^2$, we have
\begin{equation*}
	C_1K(h(v)^2)\leqslant E_{\mathcal V}E_{\mathcal P}(H)=E_{\mathcal V}(H)\leqslant E_{\mathcal V}(C_2h(v)^2)=C_2h(v)^2.
\end{equation*}
It follows that $K(h^2)\leqslant(C_2/C_1)h^2$.

For simplicity, let $\Theta_r(z)=\Psi_r(z)+\Psi_r(z^{-1})$.
Since $K$ is contractive on $L^2(\mathbb T,\mu)$ and commutes with $J$ and every $R_w$, we have
\begin{equation*}
	K\Theta_r\leqslant(C_2/C_1)\Theta_r
	\qquad\text{almost everywhere}.
\end{equation*}
Since $h$ is bounded, $\Theta_r$ and $K\Theta_r$ are continuous functions in $C(\mathbb T)$.
Hence the above inequality holds everywhere.
We define $F_r(z)=\Theta_r(z)/\Theta_r(1)$.
Then
\begin{equation}\label{equ SG-rigid-bound}
	F_r(1)=1,\qquad (KF_r)(1)\leqslant C_2/C_1.
\end{equation}
Let $\ell$ be a state on $C(\mathbb T)$ given by $\ell(g)=(Kg)(1)$ and $(\pi_\ell,\mathcal H_\ell,\xi_\ell)$ the corresponding GNS-construction.
Let $\nu$ be the spectral measure of the unitary operator $V=\pi_\ell(z)$ associated with $\xi_\ell$.
Then
\begin{equation*}
	\ell(g)=\langle g(V)\xi_\ell,\xi_\ell\rangle
	=\int_{\mathbb T}g(z)\,d\nu(z),
	\qquad
	\gamma_n=\int_{\mathbb T}z^n\,d\nu(z).
\end{equation*}
In particular, \eqref{equ SG-rigid-bound} becomes
\begin{equation*}
	\int_{\mathbb T}F_r(z)\,d\nu(z)=\ell(F_r)=(KF_r)(1)\leqslant C_2/C_1.
\end{equation*}
Applying Fatou's lemma along $r=1/k$, and using \eqref{equ SG-rigid-1} and \eqref{equ SG-rigid-2}, we obtain
\begin{equation*}
	\sum_{n=1}^{\infty}a_n\nu(\{\alpha^n,\alpha^{-n}\})\leqslant C_2/C_1,\qquad
	\nu(\mathbb T\setminus\Lambda)=0.
\end{equation*}
Let $q_0=\nu(\{1\})$, $q_n=\nu(\{\alpha^n,\alpha^{-n}\})$ for $n\geqslant 1$, and $\psi(t)=\sum_{n=1}^{\infty}q_n\bigl(1-\cos(2\pi nt)\bigr)$.
Then
\begin{equation*}
	\sum_{n=1}^\infty\frac{a_n}{n^2}=\infty,\qquad
	\sum_{n=1}^\infty a_nq_n\leqslant C_2/C_1<\infty,\qquad
	\sum_{n=1}^\infty q_n=1-q_0<\infty.
\end{equation*}
Since $\gamma_n$ is real, it is clear that
\begin{equation*}
	\gamma_n=\mathrm{Re}\int_{\mathbb T}z^n\,d\nu(z)=q_0+\sum_{j=1}^{\infty}q_j\cos(2\pi jn\theta)=1-\psi(n\theta).
\end{equation*}
Thus, $\psi((m+n)\theta)\leqslant\psi(m\theta)+\psi(n\theta)$ by \eqref{equ SG-defect-subadditive}.
Since $\theta$ is irrational, continuity and periodicity of $\psi$ imply
\begin{equation*}
	\psi(s+t)\leqslant\psi(s)+\psi(t)
	\qquad(s,t\in\mathbb R).
\end{equation*}
Thus, $\psi=0$ by \Cref{lem subadditive-Fourier}.
In particular, $\gamma_1=1$, and hence
\begin{equation*}
	\|v-E_{\mathcal P}(v)\|_2^2
	=1-\|E_{\mathcal P}(v)\|_2^2
	=1-\gamma_1=0.
\end{equation*}
Therefore, $v\in\mathcal P$.
Combining this with $u\in\mathcal P$, we obtain $\mathcal P=\mathcal R$.
\end{proof}

Recall from \cite{Nelson-74} that the \emph{measure topology} on $L^0(\mathcal R,\tau)$ is induced by the neighborhood base $\{U(\varepsilon,\delta)\colon\varepsilon,\delta>0\}$ at zero, where $U(\varepsilon,\delta)$ is the set of all operators $A$ in $L^0(\mathcal R,\tau)$ such that
\begin{equation*}
	\tau(1-e)<\delta,\qquad Ae\in\mathcal R,\qquad\|Ae\|<\varepsilon,\qquad\text{for some projection}~e\in\mathcal R.
\end{equation*}
The following lemma converts the generation property in \Cref{lem SG-rigid} into strong generation.

\begin{lemma}\label{lem SG-rigid-again}
Let $h$ be a nonnegative function in $L^\infty(\mathbb T,\mu)$ such that $W^*(u,H)=\mathcal R$ for every positive operator $H\in\mathcal R$
satisfying
\begin{equation*}
	C_1h(v)^2\leqslant H\leqslant C_2h(v)^2
\end{equation*}
for some $C_1,C_2>0$.
Suppose that $h>0$ almost everywhere and
\begin{equation*}
	f(z)=\frac{h(\alpha z)}{h(z)}\in L^\infty(\mathbb T,\mu).
\end{equation*}
Then $T=uf(v)$ is a strong generator of $\mathcal R$.
\end{lemma}

\begin{proof}
Let $X$ be an invertible operator in $\mathcal R$, and
\begin{equation*}
	S=XTX^{-1},\qquad
	\mathcal S=W^*(S),\qquad
	\mathcal U=W^*(u),\qquad
	B=Xh(v).
\end{equation*}
Choose $C_1,C_2>0$ such that $C_1\leqslant X^*X\leqslant C_2$.
Since $E_{\mathcal U}(x)=\tau(x)1$ for $x\in W^*(v)$, we have
\begin{equation}\label{equ B-comparison}
	C_1h(v)^2\leqslant B^*B\leqslant C_2h(v)^2,\qquad C_1\tau(h(v)^2)\leqslant E_{\mathcal U}(B^*B)\leqslant C_2\tau(h(v)^2).
\end{equation}
Let $B_0=E_{\mathcal U}(B^*B)\in\mathcal U$ and $A_n=\frac1n\sum_{k=0}^{n-1}S^{*k}S^k\in\mathcal S$.
Then $B_0$ is invertible and positive.
Since $A_n\geqslant 1/n$, it is invertible in $\mathcal S$, and $A_n^{-1}\in\mathcal S$.
Since $SB=Bu$ and $\{u\}'\cap\mathcal R=\mathcal U$, the mean ergodic theorem gives
\begin{equation*}
	B^*A_nB
	=\frac1n\sum_{k=0}^{n-1}u^{-k}B^*Bu^k\to B_0
	\qquad\text{in }L^2(\mathcal R,\tau).
\end{equation*}
By continuity of inversion and multiplication in the measure topology, we have
\begin{equation*}
	A_n^{-1}
	=B(B^*A_nB)^{-1}B^*\to BB_0^{-1}B^*
	\qquad\text{in measure}.
\end{equation*}
Let $Q=BB_0^{-1}B^*\in\mathcal R$.
Since $A_n^{-1}\in\mathcal S$ and $A_n^{-1}\to Q$ in measure, we have
\begin{equation}\label{equ Q-in-S}
	Q\in L^0(\mathcal S,\tau)\cap\mathcal R=\mathcal S.
\end{equation}
Let $B_0^{-1/2}B^*=WQ^{1/2}$ be the polar decomposition, where $W$ is a unitary operator in $\mathcal R$.
Then $B=Q^{1/2}W^*B_0^{1/2}$.
Since $B_0\in\mathcal U$ and $SB=Bu$, we have
\begin{equation*}
	SQ^{1/2}W^*B_0^{1/2}=SB=Bu=Q^{1/2}W^*B_0^{1/2}u=Q^{1/2}W^*uB_0^{1/2}.
\end{equation*}
Combining this with \eqref{equ Q-in-S}, we obtain
\begin{equation*}
	W^*uW=Q^{-1/2}SQ^{1/2}
	\in L^0(\mathcal S,\tau)\cap\mathcal R
	=\mathcal S.
\end{equation*}
Thus, $u\in W\mathcal SW^*$ and hence $B_0\in\mathcal U\subseteq W\mathcal SW^*$.
By \eqref{equ Q-in-S}, we have
\begin{equation*}
	B^*B=B_0^{1/2}WQW^*B_0^{1/2}\in W\mathcal SW^*.
\end{equation*}
By \eqref{equ B-comparison} and the hypothesis on $h$, we obtain
\begin{equation*}
	\mathcal R=W^*(u,B^*B)
	\subseteq W\mathcal SW^*
	\subseteq\mathcal R.
\end{equation*}
Therefore, $\mathcal S=\mathcal R$.
This completes the proof.
\end{proof}

\subsection{Relative commutants under irreducible inclusions}

We now seek a condition ensuring that $\{T\}'\cap\mathcal M=\{T\}'\cap\mathcal R$ whenever $\mathcal R\subseteq\mathcal M$ is an irreducible inclusion of type $\mathrm{II}_1$ factors.
We prove a boundedness lemma for measurable operators satisfying an intertwining relation.

\begin{lemma}\label{lem Hx-XH}
Let $(\mathcal M,\tau)$ be a tracial von Neumann algebra and $\mathcal N$ a von Neumann subalgebra of $\mathcal M$.
Suppose that $H$ is a nonzero operator in $\mathcal M$ such that
\begin{equation*}
    E_{\mathcal N}(H^*H)\in\mathbb C1.
\end{equation*}
If $x\in L^0(\mathcal N,\tau)$ and $X\in\mathcal M$ satisfy $Hx=XH$,
then $x\in\mathcal N$ and $\|x\|\leqslant\|X\|$.
\end{lemma}

\begin{proof}
Write $E_{\mathcal N}(H^*H)=c1$, where $c=\tau(H^*H)>0$.
Let $b>a>\|X\|$ and $e=\chi_{(a,b]}(|x|)\in\mathcal N$.
Then $xe\in\mathcal N$ is bounded and
\begin{equation*}
    \|Hxe\|_2^2
    =\tau\bigl((xe)^*H^*Hxe\bigr)
    =\tau\bigl((xe)^*E_{\mathcal N}(H^*H)xe\bigr)
    =c\|xe\|_2^2.
\end{equation*}
Similarly, $\|He\|_2^2=c\|e\|_2^2$.
It follows from $Hx=XH$ that
\begin{equation*}
    ca^2\|e\|_2^2
    \leqslant c\|xe\|_2^2
    =\|Hxe\|_2^2
    =\|XHe\|_2^2
    \leqslant\|X\|^2\|He\|_2^2
    =c\|X\|^2\|e\|_2^2.
\end{equation*}
Since $a>\|X\|$, we have $\chi_{(a,b]}(|x|)=e=0$.
By letting $a\to\|X\|$ and $b\to\infty$, we conclude that
\begin{equation*}
    \chi_{(\|X\|,\infty)}(|x|)=0.
\end{equation*}
Therefore, $x$ is bounded with $\|x\|\leqslant\|X\|$ and $x\in L^0(\mathcal N,\tau)\cap\mathcal M=\mathcal N$.
\end{proof}

Next we provide another condition on the function $\Psi_r$ defined in \Cref{lem SG-rigid}, which ensures that the relative commutant of $T=uh(\alpha v)h(v)^{-1}$ does not enlarge under an irreducible inclusion.

\begin{lemma}\label{lem relative-commutant-in-M}
For a nonnegative function $h\in L^\infty(\mathbb T,\mu)$, $r>0$, and $z\in\mathbb T$, let
\begin{equation*}
	\Psi_r(z)=\int_{-r}^r
	h^2(ze^{2\pi it})\,dt,
	\qquad
	\Lambda=\{\alpha^n\colon n\in\mathbb Z\}.
\end{equation*}
Suppose that $\Psi_r(1)>0$ for every $r>0$ and
\begin{equation}\label{equ Phi-1}
	\lim_{r\to0}\frac{\Psi_r(z)}{\Psi_r(1)}=\infty
	\qquad(z\in\mathbb T\setminus\Lambda).
\end{equation}
Suppose also that $h>0$ almost everywhere, $f(z)=\frac{h(\alpha z)}{h(z)}$ is a function in $L^\infty(\mathbb T,\mu)$, and $T=uf(v)\in\mathcal R$.
Then
\begin{equation*}
	\{T\}'\cap\mathcal M=\{T\}'\cap\mathcal R
\end{equation*}
for any irreducible inclusion $\mathcal R\subseteq\mathcal M$ of type $\mathrm{II}_1$ factors.
\end{lemma}

\begin{proof}
Let $X$ be an operator in $\{T\}'\cap\mathcal M$, $H=h(v)$, and $\mathcal N=\{u\}'\cap\mathcal M$.
Then $H$ is invertible in $L^0(\mathcal M,\tau)$ and
\begin{equation*}
    TH=uf(v)h(v)=uh(\alpha v)=h(v)u=Hu.
\end{equation*}
Let $x=H^{-1}XH\in L^0(\mathcal M,\tau)$.
It follows from $TX=XT$ that $ux=xu$, that is,
\begin{equation*}
	x\in\{u\}'\cap L^0(\mathcal M,\tau)=L^0(\mathcal N,\tau).
\end{equation*}
For every $n\in\mathbb Z$, we have $E_{\mathcal N}(v^n)=u^*E_{\mathcal N}(v^n)u=E_{\mathcal N}(u^*v^nu)=\alpha^nE_{\mathcal N}(v^n)$, and hence
\begin{equation*}
	E_{\mathcal N}(v^n)=
	\begin{cases}
		1,& n=0,\\
		0,& n\ne 0.
	\end{cases}
\end{equation*}
It follows that $E_{\mathcal N}(H^2)=\tau(H^2)1$.
Since $Hx=XH$, \Cref{lem Hx-XH} implies that
\begin{equation*}
    x\in\mathcal N,\qquad \|x\|\leqslant\|X\|.
\end{equation*}
Since $u^*(x^*v^nxv^{-m})u=\alpha^{n-m}(x^*v^nxv^{-m})$, we have
\begin{equation*}
	\tau(x^*v^nxv^{-m})=0\qquad\text{whenever}~m\ne n.
\end{equation*}
Let $V$ be the unitary operator on $L^2(\mathcal M,\tau)$ defined by $V\xi=v\xi v^*$, and let $\nu$ be its spectral measure associated
with the vector $x\in L^2(\mathcal M,\tau)$.
Then
\begin{equation*}
	\tau(x^*v^nxv^{-m})=
	\begin{cases}
		\int_{\mathbb T}z^n\,d\nu(z), & m=n,\\
		0, & m\ne n.
	\end{cases}
\end{equation*}
For $0<r<1/2$, let $e_r=\chi_{\{e^{2\pi it}\colon -r<t<r\}}(v)$.
We claim that
\begin{equation}\label{equ relative-commutant-scalar}
    \tau(x^*H^2xe_r)=\int_{\mathbb T}\Psi_r(z)\,d\nu(z).
\end{equation}
Expanding $e_r$ in its $L^2$-Fourier series gives
\begin{equation*}
	\tau(x^*v^nxe_r)
    =\left(\int_{-r}^{r}e^{2\pi int}\,dt\right)
      \tau(x^*v^nxv^{-n})
      =\int_{\mathbb T}\int_{-r}^{r}
      (ze^{2\pi it})^n\,dt\,d\nu(z).
\end{equation*}
Consequently, for every trigonometric polynomial $q$, we have
\begin{equation*}
    \tau(x^*q(v)xe_r)
    =\int_{\mathbb T}\int_{-r}^{r}
      q(ze^{2\pi it})\,dt\,d\nu(z).
\end{equation*}
Thus, \eqref{equ relative-commutant-scalar} holds by taking a sequence $\{q_j\}_{j=1}^\infty$ of trigonometric polynomials converging to $h^2$ in $L^1(\mathbb T,\mu)$.
Since $x^*H^2x=HX^*XH\leqslant\|X\|^2H^2$, we have
\begin{equation*}
	\int_{\mathbb T}\Psi_r(z)\,d\nu(z)=\tau(x^*H^2xe_r)
	\leqslant\|X\|^2\tau(H^2e_r)
    =\|X\|^2\Psi_r(1).
\end{equation*}
Applying Fatou's lemma together with \eqref{equ Phi-1}, we conclude that $\nu(\mathbb T\setminus\Lambda)=0$.
For every $n\in\mathbb Z$, let $\xi_n=\chi_{\{\alpha^n\}}(V)x\in L^2(\mathcal M,\tau)$.
Then
\begin{equation*}
    x=\sum_{n\in\mathbb Z}\xi_n
    \qquad\text{in }L^2(\mathcal M,\tau),
    \qquad
    v\xi_nv^*=\alpha^n\xi_n.
\end{equation*}
Let $U$ be the unitary operator on $L^2(\mathcal M,\tau)$ defined by $U\xi=u\xi u^*$.
Since $UV=VU$ and $Ux=x$, we have 
\begin{equation*}
	u\xi_nu^*=U\xi_n=U\chi_{\{\alpha^n\}}(V)x=\chi_{\{\alpha^n\}}(V)Ux=\chi_{\{\alpha^n\}}(V)x=\xi_n.
\end{equation*}
Thus, $u(u^{-n}\xi_n)=(u^{-n}\xi_n)u$ and $v(u^{-n}\xi_n)=(u^{-n}\xi_n)v$.
It follows that
\begin{equation*}
	u^{-n}\xi_n\in\mathcal R'\cap L^2(\mathcal M,\tau)=\mathbb C1.
\end{equation*}
Consequently, we have $\xi_n\in\mathbb Cu^n$ for every $n\in\mathbb Z$, and hence
\begin{equation*}
    x\in L^2(W^*(u),\tau)\cap\mathcal M
    =W^*(u)\subseteq\mathcal R.
\end{equation*}
Therefore, $X=HxH^{-1}\in L^0(\mathcal R,\tau)\cap\mathcal M=\mathcal R$.
This completes the proof.
\end{proof}

\subsection{Construction of strong generators}

If $T=u^mf(v)$ is $m$-homogeneous in $\mathcal R$, then $W^*(T)\subseteq W^*(u^m,v)$.
Since $W^*(u^m,v)$ is a proper von Neumann subalgebra of $\mathcal R$ whenever $|m|\ne 1$, such an operator cannot generate $\mathcal R$.
Therefore, we will focus on the case $m=1$.

We now construct a nonnegative function $h\in L^\infty(\mathbb T,\mu)$ satisfying the hypotheses of \Cref{lem SG-rigid,lem relative-commutant-in-M}.
By \Cref{lem SG-rigid-again}, the associated operator $T=u\,h(\alpha v)h(v)^{-1}$ is a strong generator of $\mathcal R$.
Moreover, its relative commutant is unchanged under every irreducible inclusion $\mathcal R\subseteq\mathcal M$ of type $\mathrm{II}_1$ factors by \Cref{lem relative-commutant-in-M}.

The following theorem provides a strong generator of $\mathcal R$ that remains strongly irreducible in every type $\mathrm{II}_1$ factor containing $\mathcal R$ as an irreducible subfactor.

\begin{theorem}\label{thm SG}
There exists a homogeneous operator $T=uf(v)$ in $\mathcal R$ such that $T$ is a strong generator of $\mathcal R$ and is strongly irreducible relative to $\mathcal M$ for every irreducible inclusion $\mathcal R\subseteq\mathcal M$ of type $\mathrm{II}_1$ factors.
\end{theorem}

\begin{proof}
By \Cref{lem SG-rigid}, \Cref{lem SG-rigid-again}, and \Cref{lem relative-commutant-in-M}, it suffices to construct a function $h$ in $L^\infty(\mathbb T,\mu)$ such that $h>0$ almost everywhere, $\frac{h(\alpha z)}{h(z)}\in L^\infty(\mathbb T,\mu)$, and $h$ satisfies \eqref{equ SG-rigid-1}, \eqref{equ SG-rigid-2}, and \eqref{equ Phi-1}.
Note that \eqref{equ Phi-1} is stronger than \eqref{equ SG-rigid-1}.

Let $0<b<1$, $0<p<1$, and
\begin{equation*}
	G(z)=\sum_{n\in\mathbb Z}b^{|n|}
	|z-\alpha^n|^{-p},
	\qquad
	h(z)=G(z)^{-1}.
\end{equation*}
By Fubini's theorem, we have
\begin{equation*}
	\int_{\mathbb T}G(z)d\mu(z)
	=\sum_{n\in\mathbb Z}b^{|n|}
	\int_{\mathbb T}|z-1|^{-p}\,d\mu(z)<\infty.
\end{equation*}
Thus, $G<\infty$ almost everywhere and $h(z)\leqslant|z-1|^p\leqslant2^p$.
In particular, $h\in L^\infty(\mathbb T,\mu)$ and $h>0$
almost everywhere.
For $\Psi_r$ defined in \Cref{lem SG-rigid}, we have
\begin{equation}\label{equ Phi-1-bound}
	0<\Psi_r(1)
	\leqslant \int_{-r}^r|e^{2\pi it}-1|^{2p}\,dt
	\leqslant C_pr^{2p+1}.
\end{equation}
Clearly, there exists $D_p>0$ such that
\begin{equation*}
	r^{p-1}\int_{-r}^r|ze^{2\pi it}-\zeta|^{-p}\,dt
	\leqslant D_p
	\qquad(z,\zeta\in\mathbb T,\ 0<r<1/2).
\end{equation*}
Moreover,
\begin{equation*}
	\lim_{r\to0}r^{p-1}\int_{-r}^r
	|ze^{2\pi it}-\zeta|^{-p}\,dt
	=
	\begin{cases}
		\kappa_p,&z=\zeta,\\
		0,&z\ne\zeta,
	\end{cases}
	\qquad
	\kappa_p=\frac{2(2\pi)^{-p}}{1-p}.
\end{equation*}
Let $\Lambda=\{\alpha^n\colon n\in\mathbb Z\}$.
Since $\sum_{n\in\mathbb Z}b^{|n|}<\infty$, Lebesgue's dominated convergence theorem gives
\begin{equation}\label{equ limit-G}
	\lim_{r\to0}r^{p-1}
	\int_{-r}^rG(ze^{2\pi it})\,dt
	=
	\begin{cases}
		\kappa_pb^{|n|},&z=\alpha^n,\ n\in\mathbb Z,\\
		0,&z\notin\Lambda.
	\end{cases}
\end{equation}
Since $hG=1$ almost everywhere, H\"older's inequality yields
\begin{equation*}
	2r\leqslant
	\left(\int_{-r}^rh^2(ze^{2\pi it})\,dt\right)^{1/3}
	\left(\int_{-r}^rG(ze^{2\pi it})\,dt\right)^{2/3}.
\end{equation*}
Consequently, we have
\begin{equation*}
	\frac1{r^{2p+1}}\int_{-r}^rh^2(ze^{2\pi it})\,dt
	\geqslant
	8\left(r^{p-1}
	\int_{-r}^rG(ze^{2\pi it})\,dt\right)^{-2}.
\end{equation*}
Combining this with \eqref{equ Phi-1-bound}, \eqref{equ limit-G}, and $G(z^{-1})=G(z)$, we obtain
\begin{equation*}
	\lim_{r\to0}\frac{\Psi_r(z)}{\Psi_r(1)}=\infty
	\qquad(z\notin\Lambda),
\end{equation*}
and
\begin{equation*}
	\liminf_{r\to0}
	\frac{\Psi_r(\alpha^n)+\Psi_r(\alpha^{-n})}{2\Psi_r(1)}
	=\liminf_{r\to0}
	\frac{\Psi_r(\alpha^n)}{\Psi_r(1)}
	\geqslant\frac{8}{C_p\kappa_p^2}b^{-2|n|}
	\qquad(n\ne0).
\end{equation*}
Thus, $h$ satisfies \eqref{equ SG-rigid-2} with $a_n=\frac{8}{C_p\kappa_p^2}b^{-2n}$, and also satisfies \eqref{equ Phi-1}.
It is clear that
\begin{equation*}
	G(\alpha z)=\sum_{n\in\mathbb Z}
	b^{|n+1|}|z-\alpha^n|^{-p}.
\end{equation*}
Since $b\leqslant b^{|n+1|-|n|}\leqslant b^{-1}$, we have $bG(z)\leqslant G(\alpha z)\leqslant b^{-1}G(z)$, that is
\begin{equation*}
	bh(z)\leqslant h(\alpha z)\leqslant b^{-1}h(z).
\end{equation*}
Let $f(z)=h(\alpha z)/h(z)\in L^\infty(\mathbb T,\mu)$ and $T=uf(v)$.
By \Cref{lem SG-rigid,lem SG-rigid-again}, $T$ is a strong generator of $\mathcal R$, and hence is strongly irreducible relative to $\mathcal R$.
By \Cref{lem relative-commutant-in-M}, its relative commutant is unchanged under every irreducible inclusion $\mathcal R\subseteq\mathcal M$.
Therefore, $T$ is strongly irreducible relative to $\mathcal M$.
This completes the proof.
\end{proof}

\begin{remark}\label{rem SG}
The operator $T$ constructed above is invertible in $\mathcal R$.
If instead we choose
\begin{equation*}
    G(z)=\sum_{n=0}^\infty
    b^n|z-\alpha^n|^{-p},
    \qquad
    h(z)=G(z)^{-1},
\end{equation*}
then we obtain a noninvertible operator in $\mathcal R$ with the same properties.
Note that $G(\alpha z)=bG(z)+|z-\alpha^{-1}|^{-p}$, and hence $f=h(\alpha z)/h(z)$ is bounded; the same estimates show that $f$ is not bounded below and that the required hypotheses remain satisfied.
We omit the details.
\end{remark}

\subsection{Applications}

We first apply the preceding construction to obtain strongly irreducible operators in arbitrary separable type $\mathrm{II}_1$ factors.
By Popa's theorem \cite{Popa-81}, every such factor contains an irreducible hyperfinite subfactor.
Combining this with \Cref{thm SG}, we obtain the following corollary, which answers the existence question raised in \cite[Question~1.6]{AMSSY-26} for type $\mathrm{II}_1$ factors with separable predual.

\begin{corollary}\label{cor SG}
There exist strongly irreducible operators in every type $\mathrm{II}_1$ factor with separable predual.
\end{corollary}

We next consider the indices of the subfactors generated by strongly irreducible operators.
For an irreducible subfactor $\mathcal R_0$ of $\mathcal R$, we call $[\mathcal R:\mathcal R_0]$ an \emph{irreducible index} of $\mathcal R$.
Let $\mathcal I_{\mathrm{irr}}(\mathcal R)$ denote the set of all irreducible indices.
Since every irreducible operator relative to $\mathcal R$ generates an irreducible subfactor, and every irreducible subfactor of $\mathcal R$ is singly generated, we have
\begin{equation*}
	\mathcal I_{\mathrm{irr}}(\mathcal R)=\{[\mathcal R:W^*(T)]\colon T~\text{is irreducible relative to}~\mathcal R\}.
\end{equation*}
This motivates the following definition.

\begin{definition}\label{def SI-index}
For every strongly irreducible operator $T$ relative to $\mathcal R$, we call $[\mathcal R:W^*(T)]$ a \emph{strongly irreducible index} of $\mathcal R$.
We denote by $\mathcal I_{\mathrm{SI}}(\mathcal R)$ the set of all strongly irreducible indices, that is,
\begin{equation*}
    \mathcal I_{\mathrm{SI}}(\mathcal R)
    =
    \{[\mathcal R:W^*(T)]\colon
    T~\text{is strongly irreducible relative to}~
    \mathcal R\}.
\end{equation*}
\end{definition}

Applying \Cref{thm SG} inside each irreducible subfactor of $\mathcal R$ gives the following equality of index sets.

\begin{corollary}\label{cor SI-index}
Every irreducible index of $\mathcal R$ is a strongly irreducible index, that is, $\mathcal I_{\mathrm{SI}}(\mathcal R)=\mathcal I_{\mathrm{irr}}(\mathcal R)$.
\end{corollary}

\begin{proof}
Clearly, $\mathcal I_{\mathrm{SI}}(\mathcal R) \subseteq\mathcal I_{\mathrm{irr}}(\mathcal R)$.
Conversely, let $\mathcal R_0\subseteq\mathcal R$ be an irreducible subfactor.
Since $\mathcal R_0$ is a hyperfinite type $\mathrm{II}_1$ factor by \cite{Connes-76}, \Cref{thm SG} provides $T\in\mathcal R_0$ such that $W^*(T)=\mathcal R_0$ and $T$ is strongly irreducible relative to $\mathcal R$.
It follows that
\begin{equation*}
	[\mathcal R:\mathcal R_0]=[\mathcal R:W^*(T)]\in\mathcal I_{\mathrm{SI}}(\mathcal R).
\end{equation*}
This completes the proof.
\end{proof}

For every von Neumann algebra $\mathcal N$, let $\operatorname{Proj}(\mathcal N)$ denote its projection lattice.
Applying \Cref{lem lattice-intertwining-L0} to the construction in \Cref{thm SG}, we obtain a strongly irreducible operator whose relative invariant subspace lattice is an atomless Boolean lattice.

\begin{corollary}\label{cor SI-Boolean-lattice}
There exists a strongly irreducible operator $T$ relative to $\mathcal R$ such that $\operatorname{Lat}_{\mathcal R}(T)\cong\operatorname{Proj}(L^\infty(\mathbb T,\mu))$.
\end{corollary}

\begin{proof}
Let $T$ be the homogeneous operator constructed in the proof of \Cref{thm SG}, and let $H=h(v)$.
Then $T$ is strongly irreducible relative to $\mathcal R$, $H$ is bounded and injective, and $TH=Hu$.
By \Cref{lem lattice-intertwining-L0}, we have
\begin{equation*}
    \operatorname{Lat}_{\mathcal R}(T)
    \cong\operatorname{Lat}_{\mathcal R}(u)
    =\operatorname{Proj}(W^*(u))
    \cong\operatorname{Proj}(L^\infty(\mathbb T,\mu)).
\end{equation*}
This completes the proof.
\end{proof}

\section{Strongly irreducible operators with noncyclic \texorpdfstring{$K_0$}{K0}-groups}\label{sec K0=Z2}

In this section, we construct an invertible homogeneous operator $T\in\mathcal R$ such that $T\otimes 1+1\otimes T$ is strongly irreducible relative to $\mathcal R\bar\otimes\mathcal R$ and the $K_0$-group of its relative commutant is isomorphic to $\mathbb Z^2$.

\subsection{Two-dimensional bilateral weighted shift operators}

We now briefly introduce a two-dimensional bilateral weighted shift model as in \Cref{sec bilateral}.
Let
\begin{equation*}
    \delta_1=(1,0),\qquad \delta_2=(0,1).
\end{equation*}
Suppose that $B_1$ and $B_2$ are two-dimensional weighted shift operators on $\ell^2(\mathbb Z^2)$ of the form
\begin{equation*}
    B_je_{\mathbf n}
    =w_{j,\mathbf n}e_{\mathbf n+\delta_j},\qquad j=1,2,
\end{equation*}
where $\{e_{\mathbf n}\}_{\mathbf n\in\mathbb Z^2}$ is the standard orthonormal basis of $\ell^2(\mathbb Z^2)$ and $\{w_{j,\mathbf n}\}_{\mathbf n\in\mathbb Z^2}$ is a family of
positive scalars for $j=1,2$.
If $B_1B_2=B_2B_1$, then $w_{1,\mathbf n}w_{2,\mathbf n+\delta_1} =w_{2,\mathbf n}w_{1,\mathbf n+\delta_2}$ for every $\mathbf n\in\mathbb Z^2$.
This compatibility allows us to choose a single sequence $\{\beta_{\mathbf n}\}_{\mathbf n\in\mathbb Z^2}$ of positive scalars such that $w_{j,\mathbf n}=\beta_{\mathbf n+\delta_j}/\beta_{\mathbf n}$ for $j=1,2$.
The reader is referred to \cite{JL-79} for related topics.

Let $\{\beta_{\mathbf n}\}_{\mathbf n\in\mathbb Z^2}$ be a family of positive scalars with $\beta_{\mathbf 0}=1$.
We define weighted shift operators on
$\ell^2(\mathbb Z^2)$ by
\begin{equation}\label{equ shift-2d}
    B_je_{\mathbf n}
    =\frac{\beta_{\mathbf n+\delta_j}}{\beta_{\mathbf n}}
       e_{\mathbf n+\delta_j},\qquad j=1,2.
\end{equation}
For our purpose, we assume that
\begin{equation}\label{equ b-2d}
    C_1\leqslant
    \frac{\beta_{\mathbf n+\delta_j}}{\beta_{\mathbf n}}
    \leqslant C_2,\qquad j=1,2,
\end{equation}
for some $C_1,C_2>0$.
Then $(B_1,B_2)$ is a commuting pair of bounded invertible operators.
Let $B^{\mathbf m}=B_1^{m_1}B_2^{m_2}$ for $\mathbf m=(m_1,m_2)\in\mathbb Z^2$.
Then
\begin{equation*}
    B^{\mathbf m}e_{\mathbf n}
    =\frac{\beta_{\mathbf n+\mathbf m}}{\beta_{\mathbf n}}
       e_{\mathbf n+\mathbf m},\qquad \mathbf n\in\mathbb Z^2.
\end{equation*}
For $Y\in\{B_1,B_2\}'$, define the matrix entries $\{Y_{\mathbf m,\mathbf n}\}_{\mathbf m,\mathbf n\in\mathbb Z^2}$ of $Y$ by
\begin{equation*}
    Ye_{\mathbf n}
    =\sum_{\mathbf m\in\mathbb Z^2}
       Y_{\mathbf m,\mathbf n}e_{\mathbf m}.
\end{equation*}
Since $YB_j=B_jY$, comparing the coefficients of $e_{\mathbf m+\delta_j}$ in $YB_je_{\mathbf n}$ and $B_jYe_{\mathbf n}$ gives that $\frac{\beta_{\mathbf n+\delta_j}}{\beta_{\mathbf n}}Y_{\mathbf m+\delta_j,\mathbf n+\delta_j}=\frac{\beta_{\mathbf m+\delta_j}}{\beta_{\mathbf m}}Y_{\mathbf m,\mathbf n}$, that is,
\begin{equation*}
    \frac{\beta_{\mathbf n+\delta_j}}
         {\beta_{\mathbf m+\delta_j}}
       Y_{\mathbf m+\delta_j,\mathbf n+\delta_j}
    =
    \frac{\beta_{\mathbf n}}{\beta_{\mathbf m}}
       Y_{\mathbf m,\mathbf n},\qquad j=1,2.
\end{equation*}
It follows that
\begin{equation*}
    Y_{\mathbf m,\mathbf n}
    =a_{\mathbf m-\mathbf n}
       \frac{\beta_{\mathbf m}}{\beta_{\mathbf n}},
    \qquad\text{where}~
    a_{\mathbf n}
    =\beta_{\mathbf n}^{-1}Y_{\mathbf n,\mathbf 0}.
\end{equation*}
In particular, the coefficients $a_{\mathbf n}$ are uniquely
determined by the column $Ye_{\mathbf 0}$.
Moreover, we have
\begin{equation*}
    Ye_{\mathbf n}
    =\sum_{\mathbf m\in\mathbb Z^2}
       a_{\mathbf m}
       \frac{\beta_{\mathbf n+\mathbf m}}{\beta_{\mathbf n}}
       e_{\mathbf n+\mathbf m}
    =\sum_{\mathbf m\in\mathbb Z^2}
       a_{\mathbf m}B^{\mathbf m}e_{\mathbf n}.
\end{equation*}
This series converges in $\ell^2(\mathbb Z^2)$ for each $\mathbf n\in\mathbb Z^2$.
This is the meaning of the formal notation $Y=\sum_{\mathbf n\in\mathbb Z^2}a_{\mathbf n}B^{\mathbf n}$.

Next we define a Hilbert space by
\begin{equation*}
	L^2(\beta)=\left\{
	\varphi(z)=\sum_{\mathbf n\in\mathbb Z^2}\widehat{\varphi}(\mathbf n)z^{\mathbf n}\colon
	\|\varphi\|^2=\sum_{\mathbf n\in\mathbb Z^2}|\widehat{\varphi}(\mathbf n)|^2\beta_{\mathbf n}^2<\infty
	\right\}.
\end{equation*}
Let $L^\infty(\beta)$ denote the set of all formal Laurent series $\varphi(z)=\sum_{\mathbf n\in\mathbb Z^2}\widehat{\varphi}(\mathbf n)z^{\mathbf n}$ such that $\varphi L^2(\beta)\subseteq L^2(\beta)$.
It follows from $1\in L^2(\beta)$ that $L^\infty(\beta)\subseteq L^2(\beta)$.
By the closed graph theorem, the multiplication operator $M_\varphi$ is bounded on $L^2(\beta)$ for every $\varphi\in L^\infty(\beta)$.
As usual, we identify $L^\infty(\beta)\cong\{M_\varphi\colon \varphi\in L^\infty(\beta)\}$.
The following proposition is analogous to Proposition 7 in \cite[p.~59]{Shields-74}.

\begin{proposition}\label{prop Shields-prop-7}
Let $J\colon\ell^2(\mathbb Z^2)\to L^2(\beta)$ be the unitary operator such that
\begin{equation*}
	Je_{\mathbf n}=\beta_{\mathbf n}^{-1}z^{\mathbf n},\qquad \mathbf n\in\mathbb Z^2.
\end{equation*}
With the notation of \eqref{equ shift-2d}, we have $JB_jJ^*=M_{z_j}$ for $j=1,2$.
\end{proposition}

\begin{proof}
By the definition of $L^2(\beta)$, the family $\{\beta_{\mathbf n}^{-1}z^{\mathbf n}\}_{\mathbf n\in\mathbb Z^2}$ is an orthonormal basis of $L^2(\beta)$.
Thus, $J$ is a unitary operator.
For $j=1,2$, we have
\begin{equation*}
    JB_je_{\mathbf n}
    =
    \frac{\beta_{\mathbf n+\delta_j}}{\beta_{\mathbf n}}
    \beta_{\mathbf n+\delta_j}^{-1}
    z^{\mathbf n+\delta_j}
    =
    \beta_{\mathbf n}^{-1}z^{\mathbf n+\delta_j}
    =
    M_{z_j}Je_{\mathbf n},\qquad \mathbf n\in\mathbb Z^2.
\end{equation*}
It follows that $JB_jJ^*=M_{z_j}$.
This completes the proof.
\end{proof}

In the following, we take
\begin{equation}\label{equ beta-n}
	\beta_{\mathbf n}=\sqrt{(1+n_1^2)(1+n_2^2)},\qquad\mathbf n=(n_1,n_2)\in\mathbb Z^2.
\end{equation}
Then \eqref{equ b-2d} holds.
Since $\sum_{n\in\mathbb Z}(1+n^2)^{-1}<4$, the Cauchy--Schwarz inequality shows that
\begin{equation}\label{equ 4F}
	\sum_{\mathbf n\in\mathbb Z^2}|\widehat{\varphi}(\mathbf n)|
	\leqslant\|\varphi\|
	\left(
	\sum_{\mathbf n\in\mathbb Z^2}\beta_{\mathbf n}^{-2}
	\right)^{1/2}
	=\|\varphi\|\sum_{n\in\mathbb Z}\frac{1}{1+n^2}
	\leqslant 4\|\varphi\|.
\end{equation}
Thus, the Fourier series $\varphi(z)=\sum_{\mathbf n\in\mathbb Z^2}\widehat{\varphi}(\mathbf n)z^{\mathbf n}$ converges absolutely and uniformly.
In particular, $L^2(\beta)$ embeds continuously into $C(\mathbb T^2)$.
We next show that it is closed under multiplication and that every element defines a bounded multiplier.

\begin{lemma}\label{lem L-2=infty}
With the notation of \eqref{equ beta-n}, we have
\begin{equation*}
	\|\varphi\|\leqslant\|M_\varphi\|\leqslant 16\|\varphi\|,\qquad \varphi\in L^2(\beta).
\end{equation*}
In particular, $L^\infty(\beta)=L^2(\beta)$.
Moreover, $L^\infty(\beta)$ is a maximal commutative subalgebra of $\mathcal B(L^2(\beta))$.
\end{lemma}

\begin{proof}
Let $c_n=\sqrt{1+n^2}$ for $n\in\mathbb Z$.
Then $\beta_{\mathbf n}=c_{n_1}c_{n_2}$ for $\mathbf n=(n_1,n_2)\in\mathbb Z^2$.
The elementary inequality $|c_m-c_n|\leqslant|m-n|\leqslant c_{m-n}$ gives
\begin{equation*}
	\left(\frac{c_n}{c_mc_{n-m}}\right)^2
	\leqslant\left(\frac{1}{c_m}+\frac{1}{c_{n-m}}\right)^2
	\leqslant\frac{2}{c_m^2}+\frac{2}{c_{n-m}^2},
\end{equation*}
and hence
\begin{equation*}
	\sum_{m\in\mathbb Z}\left(\frac{c_n}{c_mc_{n-m}}\right)^2
	\leqslant 4\sum_{m\in\mathbb Z}\frac{1}{1+m^2}<16,\qquad n\in\mathbb Z.
\end{equation*}
Consequently, for every $\mathbf n=(n_1,n_2)\in\mathbb Z^2$, we have
\begin{equation*}
	\sum_{\mathbf m\in\mathbb Z^2}\frac{\beta_{\mathbf n}^2}{\beta_{\mathbf m}^2\beta_{\mathbf n-\mathbf m}^2}
	=\prod_{j=1}^{2}
	\sum_{m\in\mathbb Z}\left(\frac{c_{n_j}}{c_mc_{n_j-m}}\right)^2\leqslant 16^2.
\end{equation*}
For $\varphi(z)=\sum_{\mathbf n\in\mathbb Z^2}\widehat{\varphi}(\mathbf n)z^{\mathbf n}\in L^2(\beta)$, it is clear that $\|M_\varphi\|\geqslant\|M_\varphi1\|=\|\varphi\|$.
By the absolute convergence \eqref{equ 4F}, for $\psi=\sum_{\mathbf n\in\mathbb Z^2}\widehat{\psi}(\mathbf n)z^{\mathbf n}\in L^2(\beta)$, we have
\begin{equation*}
	M_\varphi\psi=\varphi\psi=\sum_{\mathbf n\in\mathbb Z^2}c_{\mathbf n}z^{\mathbf n},\qquad
	c_{\mathbf n}=\sum_{\mathbf m\in\mathbb Z^2}\widehat{\varphi}(\mathbf m)\widehat{\psi}(\mathbf n-\mathbf m).
\end{equation*}
By applying the Cauchy--Schwarz inequality, we obtain
\begin{align*}
	|c_{\mathbf n}|^2\beta_{\mathbf n}^2
	&\leqslant\sum_{\mathbf m\in\mathbb Z^2}
	\frac{\beta_{\mathbf n}^2}{\beta_{\mathbf m}^2\beta_{\mathbf n-\mathbf m}^2}
	\sum_{\mathbf m\in\mathbb Z^2}
    |\widehat\varphi(\mathbf m)|^2\beta_{\mathbf m}^2
    |\widehat\psi(\mathbf n-\mathbf m)|^2\beta_{\mathbf n-\mathbf m}^2\\
    &\leqslant 16^2\sum_{\mathbf m\in\mathbb Z^2}
    |\widehat\varphi(\mathbf m)|^2\beta_{\mathbf m}^2
    |\widehat\psi(\mathbf n-\mathbf m)|^2\beta_{\mathbf n-\mathbf m}^2.
\end{align*}
It follows that
\begin{equation*}
	\|M_\varphi\psi\|^2
	\leqslant 16^2\sum_{\mathbf n\in\mathbb Z^2}\sum_{\mathbf m\in\mathbb Z^2}
    |\widehat\varphi(\mathbf m)|^2\beta_{\mathbf m}^2
    |\widehat\psi(\mathbf n-\mathbf m)|^2\beta_{\mathbf n-\mathbf m}^2
    =16^2\|\varphi\|^2\|\psi\|^2.
\end{equation*}
This proves $\|\varphi\|\leqslant\|M_\varphi\|\leqslant 16\|\varphi\|$.
In particular, $L^\infty(\beta)=L^2(\beta)$.

Since $M_\varphi M_\psi=M_{\varphi\psi}=M_\psi M_\varphi$ for all $\varphi,\psi\in L^\infty(\beta)$, we have $L^\infty(\beta)\subseteq L^\infty(\beta)'$.
Conversely, let $X\in L^\infty(\beta)'$ and $\varphi=X1\in L^2(\beta)$.
For every $\psi\in L^\infty(\beta)$, we have
\begin{equation*}
    X\psi=XM_\psi1=M_\psi X1=M_\psi\varphi=\psi\varphi=M_\varphi\psi.
\end{equation*}
Therefore, $X=M_\varphi\in L^\infty(\beta)$, and hence $L^\infty(\beta)'=L^\infty(\beta)$.
This completes the proof.
\end{proof}

By \eqref{equ 4F}, we may regard $L^\infty(\beta)$ as a Banach algebra of continuous functions on $\mathbb T^2$.
The following lemma establishes density in $C(\mathbb T^2)$ and inverse closedness at every matrix level.

\begin{lemma}\label{lem Fourier-inverse-closed}
With the notation of \eqref{equ beta-n}, the inclusion $\iota\colon L^\infty(\beta)\to C(\mathbb T^2)$ is continuous and unital with dense range such that
\begin{equation*}
    M_n(L^\infty(\beta))\cap\operatorname{GL}_n(C(\mathbb T^2))
    =
    \operatorname{GL}_n(L^\infty(\beta)),\qquad
    n\geqslant1.
\end{equation*}
\end{lemma}

\begin{proof}
We equip $L^\infty(\beta)$ with the norm
$\varphi\mapsto\|M_\varphi\|$.
By \Cref{lem L-2=infty}, this norm is equivalent to the Hilbert space
norm on $L^2(\beta)$ and makes $L^\infty(\beta)$ a unital commutative
Banach algebra.

The inclusion $\iota$ is unital, and \eqref{equ 4F} shows that
\begin{equation*}
    \|\varphi\|_\infty
    \leqslant4\|\varphi\|
    \leqslant4\|M_\varphi\|,
    \qquad \varphi\in L^\infty(\beta).
\end{equation*}
Its range is dense because it contains all Laurent polynomials,
which are uniformly dense in $C(\mathbb T^2)$.
Moreover, Laurent polynomials are dense in $L^2(\beta)$ by truncation.
The norm equivalence in \Cref{lem L-2=infty} therefore shows that
they are also dense in $L^\infty(\beta)$ in the multiplier norm.

We next identify the characters of $L^\infty(\beta)$.
Evaluation at any point of $\mathbb T^2$ is a continuous character
by the continuity of $\iota$.
Conversely, let $\chi$ be a character of $L^\infty(\beta)$ and let
\begin{equation*}
    \lambda_j=\chi(z_j),\qquad j=1,2.
\end{equation*}
Since $z_j^{-1}\in L^\infty(\beta)$, we have $\lambda_j\ne0$.
Characters of a unital Banach algebra are continuous.
Thus, for every integer $N\geqslant1$, \Cref{lem L-2=infty} gives
\begin{equation*}
    |\lambda_j|^{\pm N}
    =|\chi(z_j^{\pm N})|
    \leqslant\|\chi\|\|M_{z_j^{\pm N}}\|
    \leqslant16\|\chi\|\sqrt{1+N^2}.
\end{equation*}
Taking $N$-th roots and letting $N\to\infty$, we obtain
$|\lambda_j|\leqslant1$ and $|\lambda_j|^{-1}\leqslant1$.
Hence $|\lambda_j|=1$ for $j=1,2$, and
$\lambda=(\lambda_1,\lambda_2)\in\mathbb T^2$.
For every Laurent polynomial $p$, multiplicativity and linearity give
$\chi(p)=p(\lambda)$.
By the density of Laurent polynomials and the continuity of both
$\chi$ and evaluation at $\lambda$, we obtain
\begin{equation*}
    \chi(\varphi)=\varphi(\lambda),
    \qquad \varphi\in L^\infty(\beta).
\end{equation*}
Thus, every character is evaluation at a point of $\mathbb T^2$.
Therefore, $L^\infty(\beta)$ has maximal ideal space $\mathfrak M(L^\infty(\beta))=\mathbb T^2$.
Consequently, $L^\infty(\beta)$ is inverse closed in $C(\mathbb T^2)$.

Now let $A\in M_n(L^\infty(\beta))$ be invertible in $M_n(C(\mathbb T^2))$.
Then $\det A\in L^\infty(\beta)$ is invertible in $C(\mathbb T^2)$, and hence $(\det A)^{-1}\in L^\infty(\beta)$.
Since every entry of $\operatorname{adj}(A)$ also belongs to
$L^\infty(\beta)$, we have
\begin{equation*}
    A^{-1}=(\det A)^{-1}\operatorname{adj}(A)
    \in M_n(L^\infty(\beta)).
\end{equation*}
The reverse inclusion is immediate because $\iota$ is a unital homomorphism.
This completes the proof.
\end{proof}

The $K_0$-assertion in the following lemma follows from Karoubi's density theorem; see \cite[Theorem~3.2]{BH-21}.
For the reader's convenience, we include a direct proof which also establishes the corresponding isomorphism of $V$-monoids.

\begin{lemma}\label{lem dense-inverse-closed}
Let $\mathcal A$ and $\mathcal B$ be unital Banach algebras.
Suppose that $\iota\colon\mathcal B\to\mathcal A$ is an injective continuous unital homomorphism with dense range such that
\begin{equation}\label{equ inverse-closed}
    M_n(\mathcal B)\cap\operatorname{GL}_n(\mathcal A)
    =
    \operatorname{GL}_n(\mathcal B),\qquad
    n\geqslant1.
\end{equation}
Then inclusion $\iota\colon\mathcal B\to\mathcal A$ induces isomorphisms
\begin{equation*}
    V(\mathcal B)\cong V(\mathcal A),
    \qquad
    K_0(\mathcal B)\cong K_0(\mathcal A).
\end{equation*}
\end{lemma}

\begin{proof}
We identify $\mathcal B$ with $\iota(\mathcal B)$.

\emph{Surjectivity.}
Let $\gamma=\{z\in\mathbb C\colon|z-1|=1/2\}$, oriented counterclockwise.
For every idempotent $P$ in $M_n(\mathcal A)$, we have $\gamma\cap\sigma(P)=\varnothing$.
By the density of $\mathcal B$ in $\mathcal A$, we may choose $B$ in $M_n(\mathcal B)$ sufficiently close to $P$ in $M_n(\mathcal A)$ such that $\gamma\cap\sigma(B)=\varnothing$.
Note that the spectra of $B$ in $M_n(\mathcal A)$ and $M_n(\mathcal B)$ coincide by \eqref{equ inverse-closed}.
Let $Q$ be an idempotent in $M_n(\mathcal B)$ defined by
\begin{equation*}
    Q=
    \frac{1}{2\pi i}
    \int_{|z-1|=1/2}(z-B)^{-1}\,dz.
\end{equation*}
Continuity of $\iota$ identifies the contour integral defining $Q$ with the same integral in $M_n(\mathcal A)$.
Moreover, the idempotent $Q$ is arbitrarily close to $P$ in $M_n(\mathcal A)$.
Thus
\begin{equation*}
    X=QP+(1-Q)(1-P)=1+(Q-P)(2P-1)
\end{equation*}
is invertible in $M_n(\mathcal A)$ and satisfies $QX=XP$.
Therefore $[P]=[Q]$ in $V(\mathcal A)$.

\emph{Injectivity.}
Suppose that two idempotents $P,Q\in M_n(\mathcal B)$ are equivalent in $M_n(\mathcal A)$.
After adding zero blocks, we can choose $X_0,Y_0\in M_n(\mathcal A)$ such that
\begin{equation*}
	X_0=PX_0Q,\qquad Y_0=QY_0P,\qquad X_0Y_0=P,\qquad Y_0X_0=Q.
\end{equation*}
By the density of $\mathcal B$ in $\mathcal A$, we may choose $X,Y\in M_n(\mathcal B)$ close to $X_0, Y_0$ in $M_n(\mathcal A)$, such that $X=PXQ$, $Y=QYP$, and
\begin{equation*}
	A=XY+(1-P),\qquad B=YX+(1-Q)
\end{equation*}
are invertible in $M_n(\mathcal A)$.
Then $A^{-1},B^{-1}\in M_n(\mathcal B)$ by \eqref{equ inverse-closed}.
Let $\widetilde Y=YA^{-1}$.
A direct computation shows that
\begin{equation*}
	X\widetilde Y=XYA^{-1}=P,\qquad
    \widetilde YX=YXB^{-1}=Q.
\end{equation*}
Thus $P,Q$ are equivalent in $M_n(\mathcal B)$.

Therefore, $V(\mathcal B)\cong V(\mathcal A)$.
Passing to their Grothendieck groups gives the asserted isomorphism on $K_0$.
\end{proof}

\begin{proposition}\label{prop K0-sobolev}
With the notation of \eqref{equ shift-2d} and \eqref{equ beta-n}, we have
\begin{equation*}
	V(\{B_1,B_2\}')\cong V(C(\mathbb T^2)),
    \qquad
    K_0(\{B_1,B_2\}')\cong K_0(C(\mathbb T^2)).
\end{equation*}
\end{proposition}

\begin{proof}
By \Cref{lem L-2=infty}, Laurent polynomials are dense in $L^\infty(\beta)$.
It follows that
\begin{equation*}
	\{M_{z_1},M_{z_2}\}'=L^\infty(\beta)'=L^\infty(\beta).
\end{equation*}
Combining this with \Cref{prop Shields-prop-7}, we have $\{B_1,B_2\}'\cong L^\infty(\beta)$.
We complete the proof by applying \Cref{lem Fourier-inverse-closed} and \Cref{lem dense-inverse-closed}.
\end{proof}

\subsection{Two-dimensional bilateral weighted shift model}

We now develop the corresponding two-dimensional model as in \Cref{subsec bilateral-model}.
Let $(\mathcal R^{\otimes2},\tau^{\otimes2})=(\mathcal R\bar\otimes\mathcal R,\tau\otimes\tau)$, and
\begin{equation*}
	U_1=u\otimes 1,\qquad U_2=1\otimes u,\qquad
	V_1=v\otimes 1,\qquad V_2=1\otimes v.
\end{equation*}
Let $\mathcal U=W^*(u)$ and $\mathcal V=W^*(v)$.
For $\mathbf n=(n_1,n_2)\in\mathbb Z^2$ and $z=(z_1,z_2)\in\mathbb T^2$, write
\begin{equation*}
	U^{\mathbf n}=U_1^{n_1}U_2^{n_2},\qquad
	V^{\mathbf n}=V_1^{n_1}V_2^{n_2},\qquad
	z^{\mathbf n}=z_1^{n_1}z_2^{n_2},\qquad
	\alpha^{\mathbf n}z=(\alpha^{n_1}z_1,\alpha^{n_2}z_2).
\end{equation*}
Note that $\mathcal V^{\otimes2}=W^*(V_1,V_2)$.
With respect to the orthogonal decomposition
\begin{equation*}
	L^2(\mathcal R^{\otimes 2},\tau^{\otimes 2})=\bigoplus_{\mathbf n\in\mathbb Z^2}U^{\mathbf n}L^2(\mathcal V^{\otimes 2},\tau^{\otimes2}),
\end{equation*}
every operator $x$ in $\mathcal R^{\otimes2}$ has a unique $L^2$-Fourier expansion
\begin{equation}\label{equ Fourier-2d}
	x=\sum_{\mathbf n\in\mathbb Z^2}U^{\mathbf n}f_{\mathbf n}(V),\qquad f_{\mathbf n}(V)=E_{\mathcal V^{\otimes2}}(U^{-\mathbf n}x),\qquad \|f_{\mathbf n}\|_\infty\leqslant\|x\|,
\end{equation}
where $f_{\mathbf n}\in L^\infty(\mathbb T^2,\mu_2)$ and $\mu_2=\mu\otimes\mu$.
We first compute the relative commutant of the normal operator $U_1+U_2$ in $\mathcal R^{\otimes2}$.

\begin{lemma}\label{lem normal-commutant}
For $c_1,c_2\ne 0$, we have
\begin{equation*}
	\{c_1U_1+c_2U_2\}'\cap L^0(\mathcal R^{\otimes2},\tau^{\otimes2})=L^0(\mathcal U^{\otimes2},\tau^{\otimes2}).
\end{equation*}
\end{lemma}

\begin{proof}
As in \eqref{equ Fourier-2d}, every operator $x$ in $\mathcal R^{\otimes2}$ has a unique $L^2$-Fourier expansion
\begin{equation*}
	x=\sum_{\mathbf n\in\mathbb Z^2}f_{\mathbf n}(U)V^{\mathbf n}.
\end{equation*}
Similar to \Cref{lem frequently}, if $(c_1U_1+c_2U_2)x=x(c_1U_1+c_2U_2)$, then
\begin{equation*}
	(c_1U_1+c_2U_2)f_{\mathbf n}(U)V^{\mathbf n}=f_{\mathbf n}(U)V^{\mathbf n}(c_1U_1+c_2U_2),\qquad \mathbf n\in\mathbb Z^2.
\end{equation*}
It follows that
\begin{equation*}
	\bigl(c_1(1-\alpha^{n_1})z_1+c_2(1-\alpha^{n_2})z_2\bigr)f_{\mathbf n}(z_1,z_2)=0\qquad\text{almost everywhere}.
\end{equation*}
If $\mathbf n\ne\mathbf 0$, then $c_1(1-\alpha^{n_1})z_1+c_2(1-\alpha^{n_2})z_2\ne 0$ almost everywhere, and hence $f_{\mathbf n}=0$.
Therefore, $x=f_{\mathbf 0}(U)\in\mathcal U^{\otimes 2}$ and $\{c_1U_1+c_2U_2\}'\cap\mathcal R^{\otimes2}=\mathcal U^{\otimes2}$.

Now let $x\in\{c_1U_1+c_2U_2\}'\cap L^0(\mathcal R^{\otimes2},\tau^{\otimes2})$.
Since $c_1U_1+c_2U_2$ is a normal operator, the Fuglede-Putnam theorem \cite[Theorem 6]{BCSZ-18} for measurable operators shows that $x^*$ also commutes with $c_1U_1+c_2U_2$.
Let $y=x(1+x^*x)^{-1/2}$.
It follows that $y\in\{c_1U_1+c_2U_2\}'\cap\mathcal R^{\otimes2}=\mathcal U^{\otimes2}$ and therefore
\begin{equation*}
	x=y(1-y^*y)^{-1/2}\in L^0(\mathcal U^{\otimes2},\tau^{\otimes2}).
\end{equation*}
This completes the proof.
\end{proof}

The preceding lemma shows that commuting with $c_1U_1+c_2U_2$ is equivalent to commuting with both $U_1$ and $U_2$.
By conjugating with $H$, we obtain the corresponding statement for $T_1,T_2$ and a description of their joint commutant.
The following lemma can be viewed as a two-dimensional version of \Cref{thm commutant}.

\begin{lemma}\label{lem commutant-2d}
Let $H$ be invertible in $L^0(\mathcal V^{\otimes2},\tau^{\otimes2})$ and $T_j=HU_jH^{-1}$ for $j=1,2$.
Then for $c_1,c_2\ne 0$, we have
\begin{equation*}
	\{c_1T_1+c_2T_2\}'\cap\mathcal R^{\otimes 2}=\{T_1,T_2\}'\cap\mathcal R^{\otimes 2}.
\end{equation*}
Moreover, every operator $x\in\{T_1,T_2\}'\cap\mathcal R^{\otimes 2}$ has a unique $L^2$-Fourier expansion
\begin{equation*}
	x=\sum_{\mathbf n\in\mathbb Z^2}a_{\mathbf n}T^{\mathbf n},\qquad T^{\mathbf n}=HU^{\mathbf n}H^{-1}.
\end{equation*}
\end{lemma}

\begin{proof}
If $x$ commutes with $c_1T_1+c_2T_2$, then $H^{-1}xH$ commutes with $c_1U_1+c_2U_2$.
By \Cref{lem normal-commutant}, $H^{-1}xH$ belongs to $L^0(\mathcal U^{\otimes2},\tau^{\otimes2})$ and hence commutes with both $U_1$ and $U_2$.
Consequently, $x\in\{T_1,T_2\}'\cap\mathcal R^{\otimes 2}$.

Write $x=\sum_{\mathbf n\in\mathbb Z^2}U^{\mathbf n}f_{\mathbf n}(V)$ as in \eqref{equ Fourier-2d}.
Similar to \Cref{lem frequently}, for every $\mathbf n\in\mathbb Z^2$, we have
\begin{equation*}
	T_jU^{\mathbf n}f_{\mathbf n}(V)=U^{\mathbf n}f_{\mathbf n}(V)T_j,\qquad j=1,2.
\end{equation*}
It follows from \Cref{lem normal-commutant} that
\begin{equation*}
	H^{-1}U^{\mathbf n}f_{\mathbf n}(V)H\in\{U_1,U_2\}'\cap L^0(\mathcal R^{\otimes2},\tau^{\otimes2})
	=L^0(\mathcal U^{\otimes2},\tau^{\otimes2}).
\end{equation*}
Since $H\in L^0(\mathcal V^{\otimes2},\tau^{\otimes2})$, we can write $H^{-1}U^{\mathbf n}f_{\mathbf n}(V)H=U^{\mathbf n}g_{\mathbf n}(V)$, where $g_{\mathbf n}$ is a measurable function on $(\mathbb T^2,\mu_2)$.
It follows that
\begin{equation*}
	g_{\mathbf n}(V)\in U^{-\mathbf n}L^0(\mathcal U^{\otimes2},\tau^{\otimes2})
	=L^0(\mathcal U^{\otimes2},\tau^{\otimes2}),
\end{equation*}
Thus, $g_{\mathbf n}(V)=a_{\mathbf n}$ for some scalar $a_{\mathbf n}\in\mathbb C$.
Consequently, we have
\begin{equation*}
	U^{\mathbf n}f_{\mathbf n}(V)=a_{\mathbf n}HU^{\mathbf n}H^{-1}=a_{\mathbf n}T^{\mathbf n}.
\end{equation*}
This completes the proof.
\end{proof}

Analogous to the hypothesis \eqref{equ b-2d}, let $h$ be a measurable function on $(\mathbb T^2,\mu_2)$ such that $h>0$ almost everywhere and
\begin{equation}\label{equ h-2d}
	C_1\leqslant\frac{h(\alpha^{\delta_j}z)}{h(z)}\leqslant C_2,\qquad j=1,2,
\end{equation}
for some $C_1,C_2>0$.
Let $H=h(V)\in L^0(\mathcal V^{\otimes2},\tau^{\otimes2})$ and $T_j=HU_jH^{-1}$.
Then $T_1$ and $T_2$ are bounded invertible operators in $\mathcal R^{\otimes2}$ of the form
\begin{equation}\label{equ T-2d}
    T_j
    =
    U_j\frac{h(\alpha^{\delta_j}V)}{h(V)},\qquad j=1,2.
\end{equation}
Moreover, we have
\begin{equation*}
    T^{\mathbf m}
    =
    T_1^{m_1}T_2^{m_2}
    =
    U^{\mathbf m}
    \frac{h(\alpha^{\mathbf m}V)}{h(V)},
    \qquad \mathbf m\in\mathbb Z^2.
\end{equation*}
As in \eqref{equ shift-2d}, let $B_{j,z}$ be the shift operators on $\ell^2(\mathbb Z^2)$ defined by
\begin{equation*}
	B_{j,z}e_{\mathbf n}=\frac{h(\alpha^{\mathbf n+\delta_j}z)}{h(\alpha^{\mathbf n}z)}e_{\mathbf n+\delta_j},\qquad j=1,2.
\end{equation*}
Then
\begin{equation*}
	B_z^{\mathbf m}e_{\mathbf n}=B_{1,z}^{m_1}B_{2,z}^{m_2}e_{\mathbf n}=\frac{h(\alpha^{\mathbf n+\mathbf m}z)}{h(\alpha^{\mathbf n}z)}e_{\mathbf n+\mathbf m},\qquad\mathbf n\in\mathbb Z^2.
\end{equation*}
For $x\in\mathcal R^{\otimes2}$, let $L_x$ denote the left multiplication operator by $x$ on $L^2(\mathcal R^{\otimes2},\tau^{\otimes2})$.
Similar to \eqref{equ T-direct-integral}, there exists a unitary operator
\begin{equation*}
	W\colon L^2(\mathcal R^{\otimes 2},\tau^{\otimes 2})\to L^2(\mathbb T^2,\ell^2(\mathbb Z^2))
\end{equation*}
defined by $\bigl(WU^{\mathbf n}g_{\mathbf n}(V)\bigr)(z)=g_{\mathbf n}(z)e_{\mathbf n}$, such that
\begin{equation*}
	L^2(\mathbb T^2,\ell^2(\mathbb Z^2))
    =
    \int_{\mathbb T^2}^{\oplus}\ell^2(\mathbb Z^2)\,d\mu_2(z),\qquad
    WL_{T_j}W^*
    =
    \int_{\mathbb T^2}^{\oplus}B_{j,z}\,d\mu_2(z).
\end{equation*}
It follows that
\begin{equation*}
	WL_{T^{\mathbf m}}W^*
    =
    \int_{\mathbb T^2}^{\oplus}B_z^{\mathbf m}\,d\mu_2(z),\qquad \mathbf m\in\mathbb Z^2.
\end{equation*}
Similar to \Cref{lem fiber-Laurent-norm}, there exists a null set $E_h\subseteq\mathbb T^2$ such that
\begin{equation}\label{equ fiber-Laurent-norm-2d}
    \|p(T_1,T_2)\|=\|p(B_{1,z},B_{2,z})\|,
    \qquad z\in\mathbb T^2\setminus E_h,
\end{equation}
for every Laurent polynomial $p(z_1,z_2)$.
The following lemma is a two-dimensional version of \Cref{lem Laurent-comparison-commutant}.

\begin{lemma}\label{lem Laurent-comparison-commutant-2d}
With the preceding notation \eqref{equ shift-2d} and \eqref{equ T-2d}, suppose that there exists $C\geqslant1$ such that
\begin{equation*}
    \|p(B_1,B_2)\|\leqslant\|p(T_1,T_2)\|\leqslant C\|p(B_1,B_2)\|
\end{equation*}
for every Laurent polynomial $p$.
Then there exists a unique unital normal algebra isomorphism
\begin{equation*}
    \Phi\colon\{T_1,T_2\}'\cap\mathcal R^{\otimes2}\to\{B_1,B_2\}',
    \qquad \Phi(T^{\mathbf m})=B^{\mathbf m}.
\end{equation*}
Moreover, $\|\Phi\|=1$, $\|\Phi^{-1}\|\leqslant C$, and $\Phi^{-1}$ is also normal.
\end{lemma}

\begin{proof}
By the above argument, every $Y\in\{B_1,B_2\}'$ has a unique formal Fourier expansion $Y=\sum_{\mathbf n\in\mathbb Z^2}b_{\mathbf n}B^{\mathbf n}$.
By \Cref{lem commutant-2d}, every $X\in\{T_1,T_2\}'\cap\mathcal R^{\otimes2}$ has a unique $L^2$-Fourier expansion $X=\sum_{\mathbf n\in\mathbb Z^2}a_{\mathbf n}T^{\mathbf n}$.
Define the two-dimensional Fej\'er sums of $X$ and $Y$ by
\begin{equation*}
	\sigma_N(X)=\sum_{|n_1|\leqslant N,|n_2|\leqslant N}
	\left(1-\frac{|n_1|}{N+1}\right)\left(1-\frac{|n_2|}{N+1}\right)a_{\mathbf n}T^{\mathbf n},
\end{equation*}
and
\begin{equation*}
	\sigma_N(Y)=\sum_{|n_1|\leqslant N,|n_2|\leqslant N}
	\left(1-\frac{|n_1|}{N+1}\right)\left(1-\frac{|n_2|}{N+1}\right)b_{\mathbf n}B^{\mathbf n},
\end{equation*}
These sums converge ultraweakly to $X$ and $Y$, respectively, and satisfy
\begin{equation*}
    \|\sigma_N(X)\|\leqslant\|X\|,
    \qquad
    \|\sigma_N(Y)\|\leqslant\|Y\|.
\end{equation*}
The remainder of the proof is identical to that of \Cref{lem Laurent-comparison-commutant}.
\end{proof}

\subsection{Two-dimensional bilateral weighted shift realization}

We now realize the weighted shift model from the preceding subsection inside $\mathcal R^{\otimes2}$.
The construction gives a strongly irreducible operator whose relative commutant has a $K_0$-group isomorphic to $\mathbb Z^2$.
Throughout this subsection, $B_1,B_2$ denote the shift operators defined by \eqref{equ shift-2d} with the profile \eqref{equ beta-n}.
We first construct a measurable set with arbitrarily long isolated orbit visits.

\begin{lemma}\label{lem isolated-orbit-visits}
There exists a measurable set $E\subseteq\mathbb T$ such that
\begin{equation*}
    E_N
    =
    \left\{
        z\in E\colon
        \alpha^jz\notin E
        ~\text{whenever}~0<|j|\leqslant N
    \right\}
\end{equation*}
has positive measure for every integer $N\geqslant1$.
\end{lemma}

\begin{proof}
We inductively choose nowhere dense closed sets $D_m\subseteq\mathbb T$ of positive measure such that all the sets
\begin{equation*}
    \alpha^jD_m,\qquad |j|\leqslant m,\quad m\geqslant1,
\end{equation*}
are pairwise disjoint.

Since $\alpha^{-1}$, $1$, and $\alpha$ are distinct, there exists a compact nowhere dense set $D_1\subseteq\mathbb T$ of positive measure such that $\alpha^{-1}D_1$, $D_1$, and $\alpha D_1$ are pairwise disjoint.
Suppose that $m\geqslant2$ and $D_1,\ldots,D_{m-1}$ have been chosen.
Let
\begin{equation*}
    K_{m-1}
    =
    \bigcup_{n=1}^{m-1}
    \bigcup_{|j|\leqslant n}\alpha^jD_n.
\end{equation*}
Since $\bigcup_{|j|\leqslant m}\alpha^jK_{m-1}$ is nowhere dense and closed, we can choose
\begin{equation*}
    z_m\notin
    \bigcup_{|j|\leqslant m}\alpha^jK_{m-1}.
\end{equation*}
The points $\alpha^jz_m$, $|j|\leqslant m$, are distinct and lie outside $K_{m-1}$.
Hence there exists an open arc $I_m$ containing $z_m$ such that the arcs $\alpha^jI_m$, $|j|\leqslant m$, are pairwise disjoint and avoid $K_{m-1}$.
Choose a nowhere dense closed set $D_m\subseteq I_m$ of positive measure.
This completes the induction.

Let $E=\bigcup_{m=1}^{\infty}D_m$ and $N\geqslant 1$.
For $0<|j|\leqslant N$, pairwise disjointness implies that $\alpha^jD_N$ is disjoint from every $D_m$
and hence from $E$.
Since $D_N\subseteq E$, we have
\begin{equation*}
    D_N\subseteq E_N,
    \qquad
    \mu(E_N)\geqslant\mu(D_N)>0.
\end{equation*}
This completes the proof.
\end{proof}

Fix the set $E$ from \Cref{lem isolated-orbit-visits}.
To distinguish the one-dimensional coefficients from the two-dimensional profile $\beta_{\mathbf n}$, write $c_n=\sqrt{1+n^2}$ for $n\in\mathbb Z$ and define
\begin{equation}\label{equ g-definition}
    g(z)
    =
    \sum_{n\in\mathbb Z}c_n^{-2}\chi_E(\alpha^nz)
    =\sum_{n\in\mathbb Z}\frac{\chi_E(\alpha^nz)}{1+n^2},
    \qquad z\in\mathbb T.
\end{equation}
The summability of $\{c_n^{-2}\}_{n\in\mathbb Z}$ ensures that $g$ is bounded.
The following lemma provides an approximation on isolated orbit segments.

\begin{lemma}\label{lem orbit-profile}
With the notation of \eqref{equ g-definition}, the function $g\in L^\infty(\mathbb T,\mu)$ is not invertible and
\begin{equation}\label{equ adjacent-profile}
	0<g(z)<4,\qquad
	\frac13g(z)
    \leqslant g(\alpha z)
    \leqslant 3g(z),\qquad\text{almost everywhere}.
\end{equation}
Moreover, for integers $0\leqslant M<N$, we have
\begin{equation}\label{equ isolated-profile}
    0\leqslant
    g(\alpha^mz)-c_m^{-2}
    \leqslant\frac2{N-M},
    \qquad z\in E_N,\quad |m|\leqslant M.
\end{equation}
\end{lemma}

\begin{proof}
Since $\sum_{n\in\mathbb Z}c_n^{-2}<4$, we have $0\leqslant g<4$ everywhere.
Moreover, $g>0$ precisely on $\bigcup_{n\in\mathbb Z}\alpha^nE$, which has full measure by ergodicity.
Thus $g>0$ almost everywhere.
The elementary inequality $1+(n\pm1)^2\leqslant3(1+n^2)$ gives
\begin{equation*}
    \frac13c_n^{-2}
    \leqslant c_{n-1}^{-2}
    \leqslant3c_n^{-2}.
\end{equation*}
Since $g(\alpha z)=\sum_{n\in\mathbb Z}c_{n-1}^{-2}\chi_E(\alpha^nz)$, comparison term by term proves \eqref{equ adjacent-profile}.

For $z\in E_N$ and $|m|\leqslant M<N$, the definition of $E_N$ in \Cref{lem isolated-orbit-visits} gives
\begin{equation*}
	0\leqslant g(\alpha^mz)-c_m^{-2}
	=\sum_{|n|>N}c_{n-m}^{-2}\chi_E(\alpha^nz)
	\leqslant\sum_{|n|>N-M}c_n^{-2}
	\leqslant\frac2{N-M}.
\end{equation*}
This proves \eqref{equ isolated-profile}.

Taking $N=2M$ and $m=M\geqslant1$, we obtain
\begin{equation*}
    g(\alpha^Mz)
    \leqslant c_M^{-2}+\frac2M,
    \qquad z\in E_{2M}.
\end{equation*}
Since $\alpha^ME_{2M}$ has positive measure and the right-hand side tends to zero, $g$ is not invertible in $L^\infty(\mathbb T,\mu)$.
\end{proof}

With the notation of \eqref{equ g-definition}, we define
\begin{equation*}
	h(z_1,z_2)
    =
    \frac1{\sqrt{g(z_1)g(z_2)}},
    \qquad
    H=h(V).
\end{equation*}
By \eqref{equ adjacent-profile}, $h$ satisfies \eqref{equ h-2d} with $C_1=1/\sqrt3$ and $C_2=\sqrt3$.
Thus, $T_1,T_2$ defined in \eqref{equ T-2d} are bounded invertible operators in $\mathcal R^{\otimes2}$ of the form
\begin{equation}\label{equ f-definition}
    T_1=uf(v)\otimes1,\qquad T_2=1\otimes uf(v),\qquad
    f(z)=\sqrt{\frac{g(z)}{g(\alpha z)}}.
\end{equation}
The next lemma supplies the norm comparison required by \Cref{lem Laurent-comparison-commutant-2d}.

\begin{lemma}\label{lem polynomial-norm-equality}
With the notation of \eqref{equ shift-2d}, \eqref{equ beta-n}, \eqref{equ g-definition}, and \eqref{equ f-definition}, we have
\begin{equation*}
    \|p(T_1,T_2)\|=\|p(B_1,B_2)\|
\end{equation*}
for every Laurent polynomial $p$.
\end{lemma}

\begin{proof}
For $j=1,2$, let $S_j=(T_j^{-1})^*$ and $C_j=(B_j^{-1})^*$.
Then
\begin{equation*}
    S_j=H^{-1}U_jH,\qquad
    C_je_{\mathbf n}
    =
    \frac{\gamma_{\mathbf n+\delta_j}}{\gamma_{\mathbf n}}
    e_{\mathbf n+\delta_j},\qquad \gamma_{\mathbf n}=\beta_{\mathbf n}^{-1}.
\end{equation*}
Let $k(z_1,z_2)=\sqrt{g(z_1)g(z_2)}$.
The direct-integral model from the preceding subsection
represents $S_j$ by the fiber shift operators
\begin{equation*}
    C_{j,z}e_{\mathbf n}
    =
    \frac{k(\alpha^{\mathbf n+\delta_j}z)}
         {k(\alpha^{\mathbf n}z)}
    e_{\mathbf n+\delta_j}.
\end{equation*}
We work on the full measure set on which $k(\alpha^{\mathbf n}z)>0$ for every $\mathbf n\in\mathbb Z^2$.
Write
\begin{equation*}
    k_{\mathbf n}(z)=k(\alpha^{\mathbf n}z),\qquad
    \chi_{\mathbf n}(z)=\chi_E(\alpha^{n_1}z_1)\chi_E(\alpha^{n_2}z_2).
\end{equation*}
Next we show that
\begin{equation}\label{equ norm-p-S-p-C}
	\|p(S_1,S_2)\|=\|p(C_1,C_2)\|,
\end{equation}
where $p$ is a Laurent polynomial of the form $p(z)=\sum_{\mathbf m\in\mathbb Z^2}a_{\mathbf m}z^{\mathbf m}$.
For simplicity, let $K=\|p(C_1,C_2)\|$.

\noindent\textbf{Upper bound.}
By \eqref{equ g-definition} and the symmetry $\gamma_{-\mathbf n}=\gamma_{\mathbf n}$, we have
\begin{equation}\label{equ k-n-z}
    k_{\mathbf n}(z)^2
    =
    \sum_{\mathbf r\in\mathbb Z^2}
        \chi_{\mathbf r}(z)
        \gamma_{\mathbf n-\mathbf r}^2.
\end{equation}
For a finitely supported vector $\sum_{\mathbf n\in\mathbb Z^2}\xi_{\mathbf n}e_{\mathbf n}\in\ell^2(\mathbb Z^2)$, we have
\begin{equation*}
    p(C_1,C_2)
    \sum_{\mathbf n\in\mathbb Z^2}\gamma_{\mathbf n}\xi_{\mathbf n}
    e_{\mathbf n}
    =
    \sum_{\mathbf n\in\mathbb Z^2}\gamma_{\mathbf n}q_{\mathbf n}e_{\mathbf n},\qquad
    q_{\mathbf n}
    =
    \sum_{\mathbf m\in\mathbb Z^2}a_{\mathbf m}
    \xi_{\mathbf n-\mathbf m}.
\end{equation*}
Consequently, $\sum_{\mathbf n\in\mathbb Z^2}\gamma_{\mathbf n}^2|q_{\mathbf n}|^2\leqslant K^2\sum_{\mathbf n\in\mathbb Z^2}\gamma_{\mathbf n}^2|\xi_{\mathbf n}|^2$.
Applying this inequality to the translated vector $\sum_{\mathbf n\in\mathbb Z^2}\xi_{\mathbf n+\mathbf r}e_{\mathbf n}\in\ell^2(\mathbb Z^2)$ yields
\begin{equation*}
    \sum_{\mathbf n\in\mathbb Z^2}\gamma_{\mathbf n-\mathbf r}^2
        |q_{\mathbf n}|^2
    \leqslant
    K^2\sum_{\mathbf n\in\mathbb Z^2}\gamma_{\mathbf n-\mathbf r}^2
        |\xi_{\mathbf n}|^2,
    \qquad \mathbf r\in\mathbb Z^2.
\end{equation*}
Multiplying by $\chi_{\mathbf r}(z)$, summing over $\mathbf r\in\mathbb Z^2$, and using \eqref{equ k-n-z}, we obtain that
\begin{equation*}
    \sum_{\mathbf n\in\mathbb Z^2}k_{\mathbf n}(z)^2|q_{\mathbf n}|^2
    \leqslant
    K^2\sum_{\mathbf n\in\mathbb Z^2}k_{\mathbf n}(z)^2
        |\xi_{\mathbf n}|^2.
\end{equation*}
On the other hand, we have
\begin{equation*}
    p(C_{1,z},C_{2,z})
    \sum_{\mathbf n\in\mathbb Z^2}k_{\mathbf n}(z)\xi_{\mathbf n}e_{\mathbf n}
    =
    \sum_{\mathbf n\in\mathbb Z^2}k_{\mathbf n}(z)q_{\mathbf n}e_{\mathbf n}.
\end{equation*}
Since each $k_{\mathbf n}(z)>0$, every finitely supported vector has the form appearing on the left.
It follows that
\begin{equation*}
    \|p(C_{1,z},C_{2,z})\|\leqslant K
    \quad\text{almost everywhere}.
\end{equation*}
Applying \eqref{equ fiber-Laurent-norm-2d} to $S_1,S_2$, we have
\begin{equation*}
    \|p(S_1,S_2)\|
    =
    \operatorname*{ess\,sup}_{z\in\mathbb T^2}
        \|p(C_{1,z},C_{2,z})\|
    \leqslant K.
\end{equation*}

\noindent\textbf{Lower bound.}
Let $\Lambda_1=\{\mathbf m\in\mathbb Z^2\colon a_{\mathbf m}\ne0\}$.
Given $\varepsilon>0$, choose a finitely supported unit vector $\xi=\sum_{\mathbf n\in\Lambda_2}\xi_{\mathbf n}e_{\mathbf n}\in\ell^2(\mathbb Z^2)$ such that
\begin{equation*}
    \|p(C_1,C_2)\xi\|>K-\varepsilon.
\end{equation*}
Choose an integer $M\geqslant0$ such that
\begin{equation*}
    \Lambda_2\cup(\Lambda_2+\Lambda_1)\subseteq[-M,M]^2.
\end{equation*}
Let $E_N$ be given by \Cref{lem isolated-orbit-visits}.
For $N>M$ and $z\in E_N\times E_N$,
\eqref{equ isolated-profile} implies that
\begin{equation*}
    \frac1{1+n_j^2}
    \leqslant g(\alpha^{n_j}z_j)
    \leqslant
    \frac1{1+n_j^2}+\frac2{N-M},
    \qquad |n_j|\leqslant M,\quad j=1,2.
\end{equation*}
Since
\begin{equation*}
    k_{\mathbf n}(z)^2
    =
    g(\alpha^{n_1}z_1)g(\alpha^{n_2}z_2),
    \qquad
    \gamma_{\mathbf n}^2
    =
    \frac1{(1+n_1^2)(1+n_2^2)},
\end{equation*}
we obtain
\begin{equation*}
    \sup_{z\in E_N\times E_N}
    \max_{\mathbf n\in[-M,M]^2\cap\mathbb Z^2}
        |k_{\mathbf n}(z)-\gamma_{\mathbf n}|
    \to0
    \qquad\text{as}~N\to\infty.
\end{equation*}
Combining this with $\min_{\mathbf n\in[-M,M]^2\cap\mathbb Z^2}\gamma_{\mathbf n}>0$, it follows that
\begin{equation*}
    \sup_{z\in E_N\times E_N}
    \max_{\substack{\mathbf m\in\Lambda_1\\
                    \mathbf n\in\Lambda_2}}
    \left|
        \frac{k_{\mathbf n+\mathbf m}(z)}{k_{\mathbf n}(z)}
        -
        \frac{\gamma_{\mathbf n+\mathbf m}}{\gamma_{\mathbf n}}
    \right|
    \to0.
\end{equation*}
The shift formulas give
\begin{align*}
    p(C_{1,z},C_{2,z})\xi
    &=
    \sum_{\mathbf m\in\Lambda_1}
    \sum_{\mathbf n\in\Lambda_2}
        a_{\mathbf m}\xi_{\mathbf n}
        \frac{k_{\mathbf n+\mathbf m}(z)}{k_{\mathbf n}(z)}
        e_{\mathbf n+\mathbf m},\\
    p(C_1,C_2)\xi
    &=
    \sum_{\mathbf m\in\Lambda_1}
    \sum_{\mathbf n\in\Lambda_2}
        a_{\mathbf m}\xi_{\mathbf n}
        \frac{\gamma_{\mathbf n+\mathbf m}}{\gamma_{\mathbf n}}
        e_{\mathbf n+\mathbf m}.
\end{align*}
Since these sums contain only finitely many terms, for
sufficiently large $N$ we have
\begin{equation*}
    \|p(C_{1,z},C_{2,z})\xi-p(C_1,C_2)\xi\|
    <\varepsilon
\end{equation*}
for almost every $z\in E_N\times E_N$.
It follows that
\begin{equation*}
    \|p(C_{1,z},C_{2,z})\|
    \geqslant\|p(C_{1,z},C_{2,z})\xi\|
    >K-2\varepsilon
\end{equation*}
on a set of positive measure, because $\mu_2(E_N\times E_N)=\mu(E_N)^2>0$.
Applying \eqref{equ fiber-Laurent-norm-2d} to $S_1,S_2$, we have
\begin{equation*}
    \|p(S_1,S_2)\|
    =
    \operatorname*{ess\,sup}_{z\in\mathbb T^2}
        \|p(C_{1,z},C_{2,z})\|
    \geqslant K-2\varepsilon.
\end{equation*}
It follows that $\|p(S_1,S_2)\|\geqslant K$ by letting $\varepsilon\to0$.
This proves \eqref{equ norm-p-S-p-C}.

Finally, for $p(z)=\sum_{\mathbf m\in\mathbb Z^2}a_{\mathbf m}z^{\mathbf m}$,
define
\begin{equation*}
    p^\sharp(z)
    =
    \sum_{\mathbf m\in\mathbb Z^2}
        \overline{a_{\mathbf m}}z^{-\mathbf m}.
\end{equation*}
Since $S_j=(T_j^{-1})^*$ and $C_j=(B_j^{-1})^*$,
we have
\begin{equation*}
    p(T_1,T_2)=p^\sharp(S_1,S_2)^*,
    \qquad
    p(B_1,B_2)=p^\sharp(C_1,C_2)^*.
\end{equation*}
Applying \eqref{equ norm-p-S-p-C} to $p^\sharp$ shows that
\begin{equation*}
    \|p(T_1,T_2)\|
    =
    \|p^\sharp(S_1,S_2)\|
    =
    \|p^\sharp(C_1,C_2)\|
    =
    \|p(B_1,B_2)\|.
\end{equation*}
This completes the proof.
\end{proof}

The norm equality in \Cref{lem polynomial-norm-equality} allows us to transfer the model commutant to $\mathcal R^{\otimes2}$.
Now we present our main result in this section.

\begin{theorem}\label{thm SI-K0-Z2}
There exists an invertible homogeneous operator $T=uf(v)$ such that for any $c_1,c_2\in\mathbb C\setminus\{0\}$, the operator $T_{c_1,c_2}=c_1(T\otimes 1)+c_2(1\otimes T)$ is strongly irreducible relative to $\mathcal R^{\otimes2}$ and satisfies
\begin{equation*}
    V(\{T_{c_1,c_2}\}'\cap\mathcal R^{\otimes2})
    \cong V(C(\mathbb T^2)),
    \qquad
    K_0(\{T_{c_1,c_2}\}'\cap\mathcal R^{\otimes2})
    \cong K_0(C(\mathbb T^2))
    \cong\mathbb Z^2.
\end{equation*}
\end{theorem}

\begin{proof}
Let $g$ be defined by \eqref{equ g-definition} and
\begin{equation*}
    f(z)=\sqrt{\frac{g(z)}{g(\alpha z)}},
    \qquad T=uf(v).
\end{equation*}
By \eqref{equ adjacent-profile},
$1/\sqrt3\leqslant f\leqslant\sqrt3$ almost everywhere.
Thus $T$ is a bounded invertible homogeneous operator in $\mathcal R$.
Let $T_1=T\otimes1$ and $T_2=1\otimes T$.
It follows from \Cref{lem commutant-2d} that
\begin{equation*}
    \{T_{c_1,c_2}\}'\cap\mathcal R^{\otimes2}
    =
    \{T_1,T_2\}'\cap\mathcal R^{\otimes2}.
\end{equation*}
By \Cref{lem polynomial-norm-equality}, the hypothesis of \Cref{lem Laurent-comparison-commutant-2d} holds with $C=1$.
Together with the proof of \Cref{prop K0-sobolev}, we have
\begin{equation*}
	\{T_1,T_2\}'\cap\mathcal R^{\otimes2}\cong\{B_1,B_2\}'\cong L^\infty(\beta).
\end{equation*}
Since $L^\infty(\beta)\subseteq C(\mathbb T^2)$, it contains no nontrivial idempotent.
Consequently, $T_{c_1,c_2}$ is strongly irreducible relative to $\mathcal R^{\otimes2}$.
Thus, \Cref{prop K0-sobolev} gives
\begin{equation*}
    V(\{T_1,T_2\}'\cap\mathcal R^{\otimes2})\cong V(C(\mathbb T^2)),
    \qquad
    K_0(\{T_1,T_2\}'\cap\mathcal R^{\otimes2})
    \cong K_0(C(\mathbb T^2)).
\end{equation*}
By \cite[Chapter~9, Theorem~1.2; Chapter~17, Theorem~3.4]{Husemoller-book-94}, we have
\begin{equation*}
	V(C(\mathbb T^2))\cong\{0\}\cup(\mathbb Z_{\geqslant 1}\times\mathbb Z),\qquad
	K_0(C(\mathbb T^2))\cong\mathbb Z^2.
\end{equation*}
This completes the proof.
\end{proof}

As shown in \Cref{eg CFJ-applications}~$(3)$, the computation of the $V$-monoid also implies that $T_{c_1,c_2}\oplus T_{c_1,c_2}$ admits infinitely many inequivalent strongly irreducible decompositions, each consisting of two summands.

Let $T$ be the operator constructed in \Cref{thm SI-K0-Z2}.
By \Cref{lem polynomial-norm-equality}, we have
\begin{equation*}
    \|T^n\|=\|B_1^n\|=\sup_{k\in\mathbb Z}\frac{\sqrt{1+(k+n)^2}}{\sqrt{1+k^2}}\leqslant1+|n|,
    \qquad n\in\mathbb Z.
\end{equation*}
Hence $r(T)\leqslant1$ and $r(T^{-1})\leqslant1$.
It follows from \Cref{cor spectrum} that $\sigma(T)=\mathbb T$.
Fix $c_1,c_2\ne0$ and let
\begin{equation*}
    \mathcal A=\{T_{c_1,c_2}\}'\cap
    (\mathcal R\bar\otimes\mathcal R).
\end{equation*}
Then $\sigma(T_{c_1,c_2})=\sigma_{\mathcal A}(T_{c_1,c_2})$.
Under the established $\mathcal A\cong L^\infty(\beta)$ and $\mathfrak M(L^\infty(\beta))\cong\mathbb T^2$, we have
\begin{align*}
    \sigma(T_{c_1,c_2})
    &=\{c_1z_1+c_2z_2\colon(z_1,z_2)\in\mathbb T^2\}\\
    &=\left\{\lambda\in\mathbb C\colon
    \bigl||c_1|-|c_2|\bigr|
    \leqslant|\lambda|\leqslant|c_1|+|c_2|\right\}.
\end{align*}
In particular, $T_{c_1,c_2}$ is invertible if and only if $|c_1|\ne|c_2|$.

\appendix
\section{Finite strongly irreducible decompositions}
\label[appendix]{app CFJ}

We present a short algebraic formulation of the criterion of Cao, Fang, and Jiang \cite[Definitions~0.2--0.3 and Theorem~1.1]{CFJ-02}, using the Krull--Schmidt characterization recorded in \cite[Introduction, item~(1)]{Facchini-06}.
Throughout, $\mathcal A$ is a nonzero unital complex algebra with identity $I$, and $\mathbb N=\{0,1,2,\ldots\}$.
We write $I_n$ for the identity of $M_n(\mathcal A)$.

A nonzero idempotent $P\in\mathcal A$ is \emph{primitive} if $P\mathcal AP$ contains no idempotents other than $0$ and $P$.
A \emph{primitive decomposition} of an idempotent is a finite sum of pairwise orthogonal primitive idempotents, where orthogonality means that both products are zero.
The zero idempotent has the empty decomposition.
The same terminology applies to matrix algebras and corners.

For idempotents $P\in M_m(\mathcal A)$ and $Q\in M_n(\mathcal A)$, write $P\sim_{\mathrm a}Q$ if there are $X\in P M_{m,n}(\mathcal A)Q$ and $Y\in Q M_{n,m}(\mathcal A)P$ such that
\begin{equation*}
    XY=P,\qquad YX=Q.
\end{equation*}
Let $V(\mathcal A)$ denote the commutative monoid of algebraic-equivalence classes of matrix idempotents over $\mathcal A$, with addition $[P]+[Q]=[P\oplus Q]$ and zero element $[0]$.
For every finite orthogonal sum, $[\sum_jP_j]=\sum_j[P_j]$.
At a fixed matrix size, equality of classes is precisely algebraic equivalence at that size.
The group $K_0(\mathcal A)$ is the algebraic Grothendieck group of $V(\mathcal A)$.
We use the same notation for other unital algebras, with algebraic equivalence always taken relative to the algebra under consideration.

We first recall the Krull--Schmidt criterion in this language.

\begin{lemma}\label{lem CFJ-Facchini}
For a nonzero unital complex algebra $\mathcal A$, the following are equivalent.
\begin{enumerate}
\item[$(1)$] Every matrix idempotent has a primitive decomposition, and any two primitive decompositions of the same idempotent have the same algebraic-equivalence types, counted with multiplicity.
\item[$(2)$] $V(\mathcal A)$ is a free commutative monoid.
\end{enumerate}
In this case, $V(\mathcal A)\cong\mathbb N^k$ for some finite $k\geqslant 1$, and primitive idempotents represent exactly the free generators.
\end{lemma}

\begin{proof}
This is the matrix-idempotent form of the Krull--Schmidt characterization in \cite[Introduction, item~(1)]{Facchini-06}.
Every finitely generated projective right $\mathcal A$-module is isomorphic to $P\mathcal A^n$ for some idempotent $P\in M_n(\mathcal A)$, and
\begin{equation*}
    \operatorname{End}_{\mathcal A}(P\mathcal A^n)
    \cong P M_n(\mathcal A)P.
\end{equation*}
Algebraic equivalence corresponds to module isomorphism, primitivity to indecomposability, and orthogonal splittings to direct-sum decompositions.
Thus, $(1)$ is exactly the cited Krull--Schmidt property.
The cited characterization gives the equivalence and identifies the indecomposable classes with the free generators.
The identity $[P]+[I_n-P]=n[I]$ for every idempotent $P\in M_n(\mathcal A)$ shows that each free generator lies in the finite support of $[I]$.
Hence the generating set is finite and, since $\mathcal A\ne0$, nonempty.
\end{proof}

The next lemma assembles equivalences between corresponding summands into a single conjugacy.

\begin{lemma}\label{lem CFJ-gluing}
Suppose that $P=\sum_{j=1}^sP_j=\sum_{j=1}^sQ_j$ are two orthogonal
idempotent splittings of a nonzero idempotent $P\in M_n(\mathcal A)$.
If $P_j\sim_{\mathrm a}Q_j$ for every $1\leqslant j\leqslant s$,
the two families are simultaneously conjugate by an invertible
element of $P M_n(\mathcal A)P$.
\end{lemma}

\begin{proof}
Choose $X_j\in Q_jM_n(\mathcal A)P_j$ and $Y_j\in P_jM_n(\mathcal A)Q_j$ such that
\begin{equation*}
    X_jY_j=Q_j,\qquad Y_jX_j=P_j.
\end{equation*}
Let $X=\sum_{j=1}^sX_j$ and $Y=\sum_{j=1}^sY_j$.
Then $X,Y\in P M_n(\mathcal A)P$, and orthogonality gives $XY=YX=P$ and $XP_j=X_j=Q_jX$.
Thus, $X$ is invertible in the corner, with inverse $Y$, and $XP_jX^{-1}=Q_j$ for $1\leqslant j\leqslant s$.
\end{proof}

We now apply these facts to strongly irreducible decompositions.
Fix $T\in\mathcal A$ and write $T^{(n)}=\operatorname{diag}(T,\ldots,T)$ for $n\geqslant 1$.
Let
\begin{equation*}
    \mathcal B=\{T\}'\cap\mathcal A
    =\{X\in\mathcal A\colon XT=TX\}.
\end{equation*}
Then $\{T^{(n)}\}'\cap M_n(\mathcal A)=M_n(\mathcal B)$.
The element $T$ is \emph{strongly irreducible} relative to $\mathcal A$ if $\mathcal B$ contains no idempotents other than $0$ and $I$.
For every idempotent $P\in\mathcal B$, we have
\begin{equation*}
    \{TP\}'\cap P\mathcal AP=P\mathcal BP.
\end{equation*}
Thus, for $P\ne0$, $TP$ is strongly irreducible in $P\mathcal AP$ exactly when $P$ is primitive in $\mathcal B$.
Accordingly, a \emph{finite strongly irreducible decomposition} of $T^{(n)}$ means a primitive decomposition of $I_n$ in $M_n(\mathcal B)$.

With this terminology, the CFJ criterion \cite{CFJ-02} takes the following algebraic form.

\begin{theorem}\label{thm CFJ-algebraic}
The following statements are equivalent.
\begin{enumerate}
\item[$(1)$] For every $n\geqslant 1$, every idempotent of $M_n(\mathcal B)$ has a primitive decomposition, and any two primitive decompositions of $I_n$ are simultaneously conjugate, up to a permutation, by an invertible element of $M_n(\mathcal B)$.
\item[$(2)$] For some integer $k\geqslant 1$, there is a monoid isomorphism $h\colon V(\mathcal B)\to\mathbb N^k$.
\end{enumerate}
In this case, write $h([I])=(n_1,\ldots,n_k)$.
Then $n_i\geqslant 1$ for every $1\leqslant i\leqslant k$ and $T$ has a finite strongly irreducible decomposition $\{E_{ij}\colon 1\leqslant i\leqslant k,\ 1\leqslant j\leqslant n_i\}$ such that
\begin{equation}\label{eq CFJ-types}
    h([E_{ij}])=e_i,\qquad
    E_{ij}\sim_{\mathrm a}E_{i_1j_1}
    \quad\Longleftrightarrow\quad i=i_1,
\end{equation}
where $e_i$ are the coordinate vectors.
Every finite strongly irreducible decomposition of $T^{(n)}$ has exactly $nn_i$ summands of type $i$.
\end{theorem}

\begin{proof}
Assume $(1)$, and let $P\in M_n(\mathcal B)$ be an idempotent.
Choose any two primitive decompositions
\begin{equation*}
    P=\sum_{i=1}^rP_i=\sum_{j=1}^sQ_j,
\end{equation*}
and fix a primitive decomposition $I_n-P=\sum_{\ell=1}^tR_\ell$.
These give two primitive decompositions of $I_n$:
\begin{equation*}
    I_n=\sum_{i=1}^rP_i+\sum_{\ell=1}^tR_\ell
       =\sum_{j=1}^sQ_j+\sum_{\ell=1}^tR_\ell.
\end{equation*}
By $(1)$, the two decompositions have the same number of summands of each algebraic-equivalence type.
For any fixed type, the common summands $R_\ell$ contribute the same multiplicity to both decompositions. Subtracting this common integer shows that the families $\{P_i\}$ and $\{Q_j\}$ also have the same number of summands of that type.
Thus primitive decompositions of $P$ are unique by type, counted with multiplicity. Since this holds for every $n$ and $P$, \Cref{lem CFJ-Facchini}, applied to $\mathcal B$, gives $(2)$.

Conversely, $(2)$ and \Cref{lem CFJ-Facchini} give existence and uniqueness by type of primitive decompositions of every matrix idempotent.
Given two primitive decompositions of $I_n$, match corresponding types and apply \Cref{lem CFJ-gluing} in $M_n(\mathcal B)$.
This proves $(1)$.

Finally, choose an idempotent $R_i\in M_{d_i}(\mathcal B)$ with $h([R_i])=e_i$.
Then
\begin{equation*}
    e_i+h([I_{d_i}-R_i])=d_i h([I]).
\end{equation*}
Taking the $i$-th coordinate gives $d_i n_i\geqslant 1$, and hence $n_i\geqslant 1$.
Primitive classes correspond precisely to the coordinate vectors.
Applying $h$ to a primitive decomposition of $I$ therefore shows that it has exactly $n_i$ summands of type $i$.
Indexing these summands as $E_{ij}$ gives \eqref{eq CFJ-types}.
Since $h([I_n])=n h([I])$, the same argument gives the stated multiplicities for every amplification.
\end{proof}

The coordinates supplied by $h$ also classify matrix idempotents and determine the ordered $K_0$-group.

\begin{corollary}\label{cor CFJ-classification}
Assume the equivalent conditions in \Cref{thm CFJ-algebraic}.
For any idempotents $P,Q\in M_n(\mathcal B)$, we have
\begin{equation*}
    [P]=[Q]\quad\Longleftrightarrow\quad
    P\sim_{\mathrm a}Q
    \quad\Longleftrightarrow\quad
    P\text{ is similar to }Q\text{ in }M_n(\mathcal B).
\end{equation*}
For every nonzero idempotent $P\in M_n(\mathcal B)$, any two primitive decompositions of $P$ are simultaneously conjugate inside $P M_n(\mathcal B)P$, up to a permutation.
Moreover, if we write $h([P])=(m_1,\ldots,m_k)$, then each primitive decomposition of $P$ has exactly $m_i$ summands of type $i$.
The vectors realized by idempotents of $M_n(\mathcal B)$ are exactly
\begin{equation*}
    \{(m_1,\ldots,m_k)\in\mathbb N^k
      \colon 0\leqslant m_i\leqslant nn_i\}.
\end{equation*}
Furthermore, we have
\begin{equation}\label{eq CFJ-ordered-K0}
    \bigl(K_0(\mathcal B),K_0(\mathcal B)^+,[I]\bigr)
    \cong\bigl(\mathbb Z^k,\mathbb N^k,(n_1,\ldots,n_k)\bigr),
\end{equation}
where the positive cone is the image of $V(\mathcal B)$ in its Grothendieck group.
\end{corollary}

\begin{proof}
Equality of classes is algebraic equivalence at the fixed matrix size.
If $[P]=[Q]$, cancellation in $V(\mathcal B)\cong\mathbb N^k$ gives
\begin{equation*}
    [I_n-P]=[I_n-Q].
\end{equation*}
Apply \Cref{lem CFJ-gluing} to the splittings $I_n=P+(I_n-P)=Q+(I_n-Q)$ to obtain similarity.
Conversely, similarity implies algebraic equivalence and hence equality of classes.
For primitive decompositions of $P$, \Cref{lem CFJ-Facchini} gives matching types, and \Cref{lem CFJ-gluing} gives simultaneous conjugacy inside its corner.
The multiplicities are the coordinates of $h([P])$.
Their bounds follow from $[P]+[I_n-P]=n[I]$.
Conversely, choosing the prescribed number of summands of each type from a primitive decomposition of $I_n$ and summing them realizes every vector within these bounds.
Finally, group completion gives \eqref{eq CFJ-ordered-K0}.
\end{proof}

\begin{remark}\label{rem CFJ-scope}
The requirement that every idempotent admit a primitive decomposition is the hereditary condition in \cite[Definition~0.2]{CFJ-02}; here it is imposed at every matrix level.
The monoid formulation retains information beyond the abstract Grothendieck group alone.
If $T$ is strongly irreducible and satisfies the equivalent conditions of \Cref{thm CFJ-algebraic}, then $I$ is primitive in $\mathcal B$, and hence $[I]$ is a free generator.
Since every coordinate of $h([I])$ is positive, this forces $k=1$ and $h([I])=1$.
Consequently, $V(\mathcal B)\cong\mathbb N$ and $K_0(\mathcal B)\cong\mathbb Z$.
\end{remark}

Suppose that $\mathcal B=\{T\}'\cap\mathcal A$ admits a monoid isomorphism $V(\mathcal B)\cong\mathbb N$ sending $[I]$ to $1$.
By \Cref{thm CFJ-algebraic,cor CFJ-classification}, every idempotent in $M_n(\mathcal B)$ is similar in that algebra to $\operatorname{diag}(I_r,0_{n-r})$ for a unique integer $0\leqslant r\leqslant n$.
In particular, every finite strongly irreducible decomposition of $T^{(n)}$ is simultaneously conjugate, up to a permutation, to the canonical decomposition into $n$ copies of $T$.

\begin{example}\label{eg CFJ-applications}
The following applications illustrate the decomposition properties of the operators studied in this paper.
\begin{enumerate}
\item[$(1)$]
The preceding arguments apply relative to $C^*(u,W^*(v))$ to the bounded operators in \Cref{prop norm-commutant-K0}, relative to $C^*(u,v)$ under \Cref{cor rotation-commutant-K0}, and relative to $\mathcal R$ under \Cref{cor H-infty-K0=Z,cor invertible-commutant-Hinfty}.

\item[$(2)$]
Let $T$ satisfy the hypotheses of either \Cref{cor H-infty-K0=Z} or \Cref{cor invertible-commutant-Hinfty}, let $c>1$, and let
\begin{equation*}
    S=T\oplus cT,\qquad
    \mathcal B_S=\{S\}'\cap M_2(\mathcal R).
\end{equation*}
The same hypotheses hold for $cT$.
Moreover, $r(cT)=cr(T)>r(T)>0$, and hence $T$ and $cT$ have different spectra.
Thus, \Cref{thm K0-split} and the corresponding computation of the individual $V$-monoids give
\begin{equation*}
    V(\mathcal B_S)\cong\mathbb N^2,
    \qquad [I_2]\mapsto(1,1),
\end{equation*}
with the coordinate vectors represented by the canonical block projections.
By \Cref{thm CFJ-algebraic}, every finite strongly irreducible decomposition of $S^{(n)}$ is simultaneously conjugate, up to a permutation, to the canonical decomposition into $n$ copies of $T$ and $n$ copies of $cT$.
The conjugating element may be chosen invertible in $M_n(\mathcal B_S)$.

\item[$(3)$]
Let $S=T_{c_1,c_2}$ be the strongly irreducible operator in \Cref{thm SI-K0-Z2}, where $c_1,c_2\ne0$, and let $\mathcal B_S=\{S\}'\cap\mathcal R^{\otimes2}$.
Then
\begin{equation*}
    V(\mathcal B_S)\cong
    \{0\}\cup(\mathbb Z_{\geqslant1}\times\mathbb Z),
    \qquad [I]\longmapsto(1,0).
\end{equation*}
It follows from this monoid identification that the primitive
classes are precisely $(1,\ell)$, $\ell\in\mathbb Z$, and that
every matrix idempotent over $\mathcal B_S$ admits a finite
primitive decomposition.

For every $j\in\mathbb N$, the equality
\begin{equation*}
    [I_2]=(2,0)=(1,j)+(1,-j)
\end{equation*}
therefore yields an orthogonal primitive decomposition
$I_2=P_j+Q_j$ in $M_2(\mathcal B_S)$ satisfying
\begin{equation*}
    [P_j]=(1,j),\qquad [Q_j]=(1,-j).
\end{equation*}
As $j$ varies over $\mathbb N$, these unordered pairs
of classes are distinct.
Hence the resulting decompositions are pairwise nonconjugate,
even up to a permutation, by invertible elements of
$M_2(\mathcal B_S)$.
Thus $S^{(2)}$ has infinitely many inequivalent finite
strongly irreducible decompositions, each consisting
of two summands.
\end{enumerate}
\end{example}

\end{document}